\documentclass[10pt,a4paper]{article}
\usepackage[latin9]{inputenc}
\usepackage{xcolor}
\usepackage{amsmath}
\usepackage{amssymb}
\usepackage{graphics}

\makeatletter

\pdfpageheight\paperheight
\pdfpagewidth\paperwidth

\newtheorem{theo}{Theorem}
\newtheorem{fact}{Fact}
\newtheorem{prop}{Proposition}
\newtheorem{coro}{Corollary}

\newtheorem{rema}{Remark}
\newtheorem{lema}{Lemma}
\newtheorem{defi}{Definition}\def\be#1\ee{\begin{equation}#1\end{equation}}

\numberwithin{theo}{section}
\numberwithin{coro}{section}
\numberwithin{rema}{section}
\numberwithin{fact}{section}

\newcommand{\ba}{\begin{eqnarray} }
\newcommand{\ea}{\end{eqnarray} }

\makeatother

\begin{document}

\global\long\def\intr{\int_{R}}
 \global\long\def\sbr#1{\left[ #1\right] }
 \global\long\def\cbr#1{\left\{  #1\right\}  }
 \global\long\def\rbr#1{\left(#1\right)}
 \global\long\def\ev#1{\mathbb{E}{#1}}
 \global\long\def\E{\mathbb{E}}
 \global\long\def\P{\mathbb{P}}
 \global\long\def\R{\mathbb{R}}
 \global\long\def\Z{\mathbb{Z}}
 \global\long\def\N{\mathbb{N}}
 \global\long\def\Q{\mathbb{Q}}
 \global\long\def\norm#1#2#3{\Vert#1\Vert_{#2}^{#3}}
 \global\long\def\pr#1{\mathbb{P}\rbr{#1}}
 \global\long\def\cleq{\lssssim}
 \global\long\def\ceq{\eqsim}
 \global\long\def\conv{\rightarrow}
 \global\long\def\Var#1{\text{Var}(#1)}
 \global\long\def\TDD{}
 \global\long\def\dd#1{\textnormal{d}#1}
 \global\long\def\inti{\int_{0}^{\infty}}
 \global\long\def\crr{\mathcal{C}([0;\infty),\R)}
 \global\long\def\sb#1{\langle#1\rangle}
 %\global\long\def\pm#1{d_{P}\rbr{#1}}
 \global\long\def\crt{\mathcal{C}([0;T],\R)}
 \global\long\def\nuu{\nu_{n;\lambda}}
 \global\long\def\ZZ{Z_{\Lambda_{n}}}
 \global\long\def\PP{\mathbb{P}_{\Lambda_{n}}}
 \global\long\def\EE{\mathbb{E}_{\Lambda_{n}}}
 \global\long\def\LL{\Lambda_{n}}
 \global\long\def\AA{\mathcal{A}}
 \global\long\def\evx{\mathbb{E}_{x}}
 \global\long\def\pin#1{1_{\cbr{#1\in\mathcal{A}}}}
 \global\long\def\Zd{\mathbb{Z}^{d}}
\global\long\def\ls{\ls}
 \global\long\def\gs{\gs}
 \global\long\def\ddp#1#2{\langle#1,#2\rangle}
 \global\long\def\intc#1{\int_{0}^{#1}}
 \global\long\def\T#1{\mathcal{P}_{#1}}
 \global\long\def\ii{\mathbf{i}}
 \global\long\def\star#1{\left.#1^{*}\right.}
 \global\long\def\pspace{\mathcal{C}}
 \global\long\def\eq{\varphi}
 \global\long\def\grad{\text{grad}}
 \global\long\def\var{\text{var}}
 \global\long\def\ab{[a,b]}
 \global\long\def\ra{\rightarrow}
 \global\long\def\TTV#1#2#3{\text{TV}^{#3}\!\rbr{#1,#2}}
 \global\long\def\V#1#2#3{\text{V}^{#3}\!\rbr{#1,#2}}
 \global\long\def\vi#1#2{\text{vi}\rbr{#1,#2}}
 \global\long\def\eqdef{:=}
 \global\long\def\UTV#1#2#3{\text{UTV}^{#3}\!\rbr{#1,#2}}
 \global\long\def\DTV#1#2#3{\text{DTV}^{#3}\!\rbr{#1,#2}}
 \global\long\def\ns{\infty}
 \global\long\def\f{:\left[a,b\right]\ra\R}
 \global\long\def\TV{\text{TV}}
 \global\long\def\osc{\text{osc}}
 \global\long\def\calu{{\cal U}^{p}\left[a,b\right]}
 \global\long\def\cont{\text{cont}}
\global\long\def\ccN{{\cal N}}
\global\long\def\Cov{\text{Cov}}
\def\bt#1\et{\begin{theo}#1\end{theo}}
\def\bl#1\el{\begin{lema}#1\end{lema}}
\def\bp#1\ep{\begin{prop}#1\end{prop}}
\def\bd#1\ed{\begin{defi}#1\end{defi}}
\global\long\def\sbr#1{\left[ #1\right] }
 \global\long\def\cbr#1{\left\{  #1\right\}  }
 \global\long\def\rbr#1{\left(#1\right)}
 \global\long\def\tbr#1{\left\langle #1\right\rangle }
%litery 'krecone'
\def\ccA{{\cal A}}
\def\ccB{{\cal B}}
\def\ccC{{\cal C}}
\def\ccD{{\cal D}}
\def\ccE{{\cal E}}
\def\ccF{{\cal F}}
\def\ccG{{\cal G}}
\def\ccH{{\cal H}}
\def\ccI{{\cal I}}
\def\ccJ{{\cal J}}
\def\ccK{{\cal K}}
\def\ccL{{\cal L}}
\def\ccM{{\cal M}}
\def\ccN{{\cal N}}
\def\ccO{{\cal O}}
\def\ccP{{\cal P}}
\def\ccQ{{\cal Q}}
\def\ccR{{\cal R}}
\def\ccS{{\cal S}}
\def\ccT{{\cal T}}
\def\ccY{{\cal Y}}
\def\ccX{{\cal X}}
\def\ccZ{{\cal Z}}
%strzalki
\def\va{\varepsilon}
\def\ra{\rightarrow}
\def\lra{\longrightarrow}
\def\la{\leftarrow}
\def\da{\downarrow}
\def\al{\langle}
\def\ar{\rangle}
 %standardowa notacja algebraiczna: ciala liczbowe, grupy, pierscienie
\def\pint{-\hspace{-11pt}\int}
\def\C{{\mathbb C}}
\def\E{\mathbf{E}}
\def\Cov{\mathbf{Cov}}
\def\P{\mathbf{P}}
\def\Var{\mathbf{Var}}
\def\EM{\mathbb{E}_\mu}
\def\F{{\mathbb F}}
\def\N{{\mathbb N}}
\def\Q{{\mathbb Q}}
\def\R{{\mathbb R}}
\def\Z{{\mathbb Z}}
\def\c{{\hat{c}}}
\def\ls{\leqslant}
\def\gs{\geqslant}
\def\for{\mbox{for}}
\def\For{\mbox{for all}}
\def\Bc{{B}}
\def\Bo{{B^{\circ}}}

\title{\textbf{Concentration of the truncated variation of fractional Brownian motions of any Hurst index, their $1/H$-variations and local times}}

\author{{Witold M. Bednorz and Rafa{\l} M. {\L}ochowski} \thanks{Research partially funded by National Science Centre, Poland, Grant no. 2022/47/B/ST1/02114.}}

\maketitle
\begin{abstract}
We obtain bounds for probabilities of deviations of the truncated variation functional
of fractional Brownian motions (fBm) of any Hurst index $H \in (0,1)$ from their expected values. Obtained bounds are optimal for large values of deviations up to multiplicative constants depending on the parameter $H$ only. As an application, we give tight bounds for tails of $1/H$-variations of fBm along Lebesgue partitions and establish the a.s. weak convergence (in $L^1$) of normalized numbers of strip crossings by the trajectories of fBm to their local times for any Hurst parameter $H \in (0,1)$.
\end{abstract}

\section{Introduction}

Let $f:[0,+\ns) \ra \R$ be a c\`adl\`ag path (right-continuous with left limits).  For $s$, $t$ such that $t>s\ge 0$, {\em the truncated variation  of $f$ on the time interval $[s,t]$ with the truncation parameter $c \gs 0$} is defined as 
\[
\TTV{f}{[s,t]}{c} : = \sup_{n \gs 1}\sup_{s \ls s_1 < t_1 \ls s_2 < t_2 < \dots \ls s_n <t_n \ls t} \sum_{i=1}^n \rbr{\left| f\rbr{t_{i}}-f\rbr{s_{i}} \right|- c}_+,
\]
where $(x)_+ = \max\rbr{x,0}$.

For $S>0$, $\TTV{f}{[0,S]}{c}$ will be shortly denoted by $\TTV{f}{S}{c}$. $\TTV{f}{[s,t]}{c}$ is always finite when $c>0$ (see \cite[Fact 2.2]{LochowskiColloquium:2013}) and for $c=0$ it coincides with the total variation of $f$ on the interval $[s,t]$, which we will denote by $\TTV{f}{[s,t]}{}$ and for $[s,t] = [0,S]$ we will denote by $\TTV{f}{S}{}$.

Truncated variation is the solution of the following optimization problem (see \cite{LochowskiColloquium:2013}): 
\[
\TTV{f}{S}{c} = \inf \cbr{ \TTV{g}{S}{}:   \Vert g - f \Vert_{\infty, [0,S]} \ls c/2},
\]
where the infimum is over all c\`adl\`ag  functions $g:[0, + \ns) \ra \R$ and $\Vert g - f \Vert_{\infty, [0,S]}$ denotes the supremum norm on the interval $[0, S]$:
\[
\Vert g - f \Vert_{\infty, [0,S]} := \sup_{s \in [0,S]} \left| g\rbr{s}-f\rbr{s} \right|.
\]
It is also related to the numbers of interval crossings and $\psi$-variation of continuous $f$ along the Lebesgue partitions, see \cite{LochowskiColloquium:2017}, \cite{AIHPlocaltimes:2024}. ($\psi$-variation of $f$ on the time interval $[0,S]$ along some partition $\Pi = \cbr{0=t_0 < t_1 \ldots < t_n =S}$ may be defined as
$
V_{\psi}^{\Pi}(f,S) := \sum_{i=1}^n \psi\rbr{\left| f\rbr{t_i}-f\rbr{t_{i-1}} \right|}
$. The Lebesgue partition for the path $f$ is obtained as  consecutive times at which $f$ attains values from the grid $c \Z + r = \cbr{r, \pm c+r, \pm 2c+r, \ldots}$.)
The convergence of $1/H$-variation (corresponding to $\psi$-variation with $\psi(x) = x^{1/H}$) of fractional Brownian motions along sequence of  Lebesgue partitions is very much related to the concentration properties of the truncated variation functional of trajectories of fractional Brownian motion, see \cite{Toyomu:2023}. Let us recall that the fractional Brownian motion $\rbr{W^H_t}_{ t\gs 0}$ (fBm in short) with the Hurst index $H \in (0,1)$ is a centered Gaussian process with the following covariance structure 
\[
\E \rbr{W^H_sW^H_t} = \frac{1}{2}\cbr{|s|^{2H} + |t|^{2H} - |t-s|^{2H}},  \ s,t\gs 0.
\]

Due to the above mentioned properties of the truncated variation, taking also advantage of the fact that it may be written as a maximum of shifted suprema of mean zero Gaussian processes, in the first part of this paper we study the concentration properties of the truncated variation functional of trajectories of fractional Brownian motion of any Hurst index $H$.
These properties were already investigated in our previous paper \cite{BednorzLochowski:2013}. However, they were formulated in much more general setting -- for all processes $X_s$, $s \in [0,S]$, satisfying the following condition:
\be \label{growth}
\E \varphi \rbr{ \frac{X_t-X_s}{|t-s|^q} } \ls 1,
\ee
where $q \in (0,1)$ and $\varphi:[0,+\ns) \ra [0,+\ns)$ is strictly increasing and satisfying a growth condition of exponential type, and few other technical conditions. Typical example of $\varphi$ is $\varphi_p(x) = \exp\rbr{x^p}-1$, $p>0$. As a result one gets the following concentration estimates for $X$ satisfying \eqref{growth} with $\varphi = \varphi_p$:
\be \label{Bernoulli}
\P \rbr{ \TTV{X}{S}{c} \gs S c^{1-1/q} \sbr{A_{p,q} + B_{p,q}u}} \ls D_{p,q} \exp\rbr{ - u^{pq} }, \ u >0.
\ee
Here $A_{p,q}$, $B_{p,q}$ and $D_{p,q}$ are positive constants depending on $p$ and $q$ only. Since the fractional Brownian motion $\rbr{W^H_t}_{ t\gs 0}$ has the following variance of increments
\[
\E \rbr{W^H_t - W^H_s}^2 = |t-s|^{2H},
\]
the process $\rbr{\tfrac{1}{2}W^H_t}_{ t\gs 0}$ satisfies condition \eqref{growth} with $q=H$ and $p=2$, thus \eqref{Bernoulli} takes the form
\be \label{BernoullifBm}
\P \rbr{ \TTV{ \tfrac{1}{2} W^H}{S}{c} \gs S c^{1-1/H} \sbr{A_{2,H} + B_{2,H}u}} \ls D_{2,H} \exp\rbr{ - u^{2H} }, \ u >0.
\ee
It is worth to mention that for $c \ls S^H$, $S c^{1-1/H} $ is comparable (up to a multiplicative bounds depending on $H$ only) with $\E \TTV{ W^H}{S}{c}$, see Remark \ref{rema1}.
In this paper we managed to obtain much better concentration, namely, if $H \in (0,1/2)$ and $c \ls S^H$ then
\[
\P \rbr{ \left| \TTV{W^H}{S}{c}-\E {\TTV{W^H}{S}{c}}\right| > Sc^{1-{1}/{H}}v} \ls {\bar A}_H \exp\rbr{ - {\bar B}_H S{c}^{-{1}/{H}} \min\rbr{v^{1+2H} , v^2}}, \ v>0,
\]
while for $H \in [1/2,1)$ and $c \ls S^H$, 
\be
\P \rbr{ \left| \TTV{W^H}{S}{c}-\E {\TTV{W^H}{S}{c}}\right| > Sc^{1-\frac{1}{H}}v} \ls {\bar A} \exp\rbr{ - {\bar B} \rbr{S{c}^{-\frac{1}{H}}}^{2-2H}  v^2}, \nonumber
\ee
see Theorem \ref{theo-1} and Theorem  \ref{theo-2}. Moreover, both estimates can not be improved for large $v$ (opposite estimates hold with other positive constants put in place of ${\bar A}$, ${\bar A}_H$, ${\bar B}$, ${\bar B}_H$), see Remarks \ref{rema2} and \ref{rema3}. Thus, in the case $H \in (0, 1/2)$ we obtained a somewhat surprising result that the truncated variation functional of a Gaussian process does not reveal Gaussian concentration phenomenon.  

Together with the  truncated variation we will consider two related functionals -- {\em upward} and {\em downward  truncated variations} -- which are defined as 
\[
\UTV{f}{[s,t]}{c} : = \sup_{n \gs 1}\sup_{s \ls s_1 < t_1 < s_2 < t_2 < \dots < s_n <t_n \ls t} \sum_{i=1}^n \rbr{ f\rbr{t_{i}}-f\rbr{s_{i}} - c}_+
\]
and
\[
\DTV{f}{[s,t]}{c} : = \sup_{n \gs 1}\sup_{s \ls s_1 < t_1 < s_2 < t_2 < \dots < s_n <t_n \ls t} \sum_{i=1}^n \rbr{ f\rbr{s_{i}} - f\rbr{t_{i}}- c}_+.
\]
For $S >  0$ we will shortly denote $\UTV{f}{[0,S]}{c}$ and $\DTV{f}{[0,S]}{c} $ by $\UTV{f}{S}{c}$ and $ \DTV{f}{S}{c}$ respectively. 
Naturally, $ \DTV{f}{S}{c} =\UTV{-f}{S}{c}$. We have also the following important relationship 
\be \label{add}
\TTV{f}{S}{c} =\UTV{f}{S}{c}  + \DTV{f}{S}{c}, 
\ee
see \cite{LochowskiColloquium:2013}. First we will obtain concentration results for $\UTV{W^H}{S}{c}$ and then, using \eqref{add}, analogous results for $\TTV{W^H}{S}{c}$. 

More general results on truncated variations and other related quantities (like for example strip (up-, down-) crossings and $c$-level crossings) will be presented in the second part of this paper -- Sect. \ref{sec:Concentration}. Using results obtained for truncated variation of fBm we will establish tight, exponential bounds for tails of $1/H$-variations of fBm along Lebesgue partitions. Combining them with some concentration type results for martingales \cite{Rio} and Girsanov's type theorem for variables depending on fBm's past (see \cite{Picard:2008}, \cite{Toyomu:2023}), we will obtain the a.s. convergence of these variations to deterministic limits, provided the "vertical" mesh of these partitions (corresponding to the truncation parameter $c$) converges to $0$ sufficiently fast (as a positive power of consecutive reciprocals of positive integers).

In the last, fourth section we will deal with local times of fBm. We will proceed there using similar techniques as in Sect.
\ref{sec:Concentration} to prove that normalized numbers of strip upcrossings by fBm's trajectories almost
surely converge in weak topology of $L^{1}(\R)$ to the local time of fBm (with any Hurst parameter $H\in(0,1)$). Let us mention that much stronger mode
of convergence of normalized numbers of strip upcrossings -- almost sure convergence to the local time of fBm -- was
proven in \cite{Toyomu:2023} for $H>1/2$ while a weaker mode of convergence of these quantities (almost sure, weak convergence of relevant measures) was proven in \cite{AIHPlocaltimes:2024} for any $H\in(0,1)$. 

In the Appendix we present some auxiliary results related to the estimates of $\E \UTV{W^H}{S}{c}$, weak variance of sums of increments of fBm and two important lemmas (Lemma A and Lemma B) used in Sect. \ref{sec:Concentration} and \ref{sect:loc_times}.

\noindent
{\bf Notation}. $\Z$ denotes the set of integers, $\N$ denotes the set of positive integers.

Throughout the paper $A$, $B$ etc. denote universal positive constants while $A_H, B_H$ etc. denote positive constants, which depend on $H$ only. For any positive quantities $a$, $b$
(depending on one or more parameters) we write  $a\sim b$  if there exists
a universal positive constant $C$ such that $ C^{-1} a \le b \le C a$. Similarly, we write $a \sim_{H} b$
if there exists a universal positive constant $C_H$, depending on $H$ only, such that $ C^{-1}_H a \le b \le C_H a$.

Conditional expectation of a r.v. $X$ with respect to a variable $Y$ will be denoted by $E \rbr{X|Y}$ and  with respect to a $\sigma$-field $\cal{F}$ will be denoted by $E \rbr{X|{\cal F}}$.

\section{Concentration of the truncated variation of fractional Brownian motion} 

\subsection{Main results}

\begin{theo} \label{theo-1} Let $S>0$, $H \in (0,1/2)$, $c \in \left(0,S^H\right]$.  There exist positive constants ${\bar A}_H \le 36$ and ${\bar B}_{H}$ (depending only on $H$) such that
for any $u >0$,
\[
\P \rbr{ \left| \TTV{W^H}{S}{c}-\E \TTV{W^H}{S}{c}\right| > u} \ls {\bar A}_H \exp\rbr{ - {\bar B}_H S^{-1}{c}^{\frac{1}{H} - 2} {u}^{1+2H} \min\rbr{S{c}^{1-\frac{1}{H}}, {u}}^{1-2H} }
\]
or, equivalently, introducing $v >0$ such that $u = S c^{1-\frac{1}{H}}v$,
\be \label{estHle12}
\P \rbr{ \left| \TTV{W^H}{S}{c}-\E {\TTV{W^H}{S}{c}}\right| > Sc^{1-\frac{1}{H}}v} \ls {\bar A}_H \exp\rbr{ - {\bar B}_H S{c}^{-\frac{1}{H}} \min\rbr{v^{1+2H} , v^2}}.
\ee
\end{theo}
Theorem \ref{theo-1} is an immediate consequence of the following Theorem \ref{theo-11} and \eqref{add}.
\begin{theo} \label{theo-11} Let $S>0$, $H \in (0,1/2)$, $c \in \left(0,S^H\right]$.  There exist positive constants ${\hat A}_H \le 18$ and ${\hat B}_{H}$ (depending only on $H$) such that
for any $u >0$,
\[
\P \rbr{ \left| \UTV{W^H}{S}{c}-\E \UTV{W^H}{S}{c}\right| > u} \ls {\hat A}_H \exp\rbr{ - {\hat B}_H S^{-1}{c}^{\frac{1}{H} - 2} {u}^{1+2H} \min\rbr{S{c}^{1-\frac{1}{H}}, {u}}^{1-2H} }
\]
or, equivalently, introducing $v >0$ such that $u = S c^{1-\frac{1}{H}}v$,
\be \label{estHle12}
\P \rbr{ \left| \UTV{W^H}{S}{c}-\E {\UTV{W^H}{S}{c}}\right| > Sc^{1-\frac{1}{H}}v} \ls {\hat A}_H \exp\rbr{ - {\hat B}_H S{c}^{-\frac{1}{H}} \min\rbr{v^{1+2H} , v^2}}.
\ee
The same estimates hold for $\DTV{W^H}{S}{c}$.
\end{theo}
The proof of Theorem \ref{theo-11} may be found at the end of Subsection \ref{subsect:2.2.1}. 
\begin{rema} \label{rema1} For $c \in \left(0,S^H\right]$ the quantity $S c^{1-{1}/{H}}$ is of the same order (up to a constant depending on $H$ only) as $\E {\UTV{W^H}{S}{c}}$, in our notation $\E {\UTV{W^H}{S}{c}} \sim_H S c^{1-{1}/{H}}$, see \eqref{eutv} and \eqref{eutvlower}. By symmetry, $\E {\DTV{W^H}{S}{c}} = \E {\UTV{W^H}{S}{c}}$ and by this and \eqref{add}, $\E {\TTV{W^H}{S}{c}} = 2 \cdot \E {\UTV{W^H}{S}{c}}$.
\end{rema}
\noindent
\begin{rema} \label{rema2} For $v \gs 1$ the estimate \eqref{estHle12} is optimal in this sense that there exists universal positive constant $T_H$ such that for $v \gs 1$,
\[
\P \rbr{ \left| \UTV{W^H}{S}{c}-\E {\UTV{W^H}{S}{c}}\right| > Sc^{1-\frac{1}{H}}v} \gs \frac{1}{2}\exp\rbr{ - {T}_H S{c}^{-\frac{1}{H}} {v^{1+2H}}},
\]
see Remark \ref{optimality}.
\end{rema}
\begin{theo} \label{theo-2} Let $S>0$, $H \in [1/2,1)$, $c \in \left(0,S^H\right]$.  There exist universal positive constants ${\bar A} \le 4$ and ${\bar B} \ge 2/\pi^2$  such that
for any $u >0$,
\[
\P \rbr{ \left| \TTV{W^H}{S}{c}-\E \TTV{W^H}{S}{c}\right| > u} \ls {\bar A} \exp\rbr{ - {\bar B} S^{-2H} u^{2} }
\]
or, equivalently, introducing $v >0$ such that $u = S c^{1-\frac{1}{H}}v$,
\be
\P \rbr{ \left| \TTV{W^H}{S}{c}-\E {\TTV{W^H}{S}{c}}\right| > Sc^{1-\frac{1}{H}}v} \ls {\bar A} \exp\rbr{ - {\bar B} \rbr{S{c}^{-\frac{1}{H}}}^{2-2H}  v^2}. \nonumber
\ee
\end{theo}
Theorem \ref{theo-2} is an immediate consequence of the following Theorem \ref{theo-21} and \eqref{add}.
\begin{theo} \label{theo-21} Let $S>0$, $H \in [1/2,1)$, $c \in \left(0,S^H\right]$.  There exist universal positive constants ${\hat A} \le 2$ and ${\hat B} \ge 2/\pi^2$  such that
for any $u >0$,
\be \label{estHge12}
\P \rbr{ \left| \UTV{W^H}{S}{c}-\E \UTV{W^H}{S}{c}\right| > u} \ls {\hat A} \exp\rbr{ - {\hat B} S^{-2H} u^{2} }
\ee
or, equivalently, introducing $v >0$ such that $u = S c^{1-\frac{1}{H}}v$,
\be \label{estHge12}
\P \rbr{ \left| \UTV{W^H}{S}{c}-\E \UTV{W^H}{S}{c}\right| > S c^{1-\frac{1}{H}}v} \ls {\hat A} \exp\rbr{ - {\hat B}  \rbr{S{c}^{-\frac{1}{H}}}^{2-2H}  v^2}.
\ee
The same estimates hold for $\DTV{W^H}{S}{c}$.
\end{theo}
The proof of Theorem \ref{theo-21} may be found at the end of Subsection \ref{subsect:2.2.2}. 
\begin{rema} \label{rema3} For $u \gs \max\rbr{\E {\UTV{W^H}{S}{c}}+c,S^H} \sim \E {\UTV{W^H}{S}{c}}$ the estimate \eqref{estHge12} is optimal in this sense that there exists universal positive constant $T$ such that for $u \gs  \max\rbr{\E {\UTV{W^H}{S}{c}}+c,S^H}$,
\[
\P \rbr{ \left| \UTV{W^H}{S}{c}-\E {\UTV{W^H}{S}{c}}\right| > u} \gs \frac{1}{2}\exp\rbr{ - TS^{-2H}u^2}.
\]
To see this it is sufficient to notice that $ \UTV{W^H}{S}{c} \gs  W^H_S - c$, $S^{-H}{W^H_S} \sim \ccN(0,1)$ and 
\begin{align*}
 \P \rbr{ \left| \UTV{W^H}{S}{c}-\E {\UTV{W^H}{S}{c}}\right| > u} & \gs \P \rbr{ W^H_S - c -\E {\UTV{W^H}{S}{c}} > u} \\
& =  \P \rbr{ S^{-H}{W^H_S} > S^{-H}\rbr{u  + \E {\UTV{W^H}{S}{c}}}+ c} \\
& \ge  \P \rbr{ S^{-H}{W^H_S} > 2S^{-H}u} \gs \frac{1}{2}\exp\rbr{ - TS^{-2H}u^2}.
\end{align*}
\end{rema}

By Theorem \ref{theo-1} and stationarity of the increments of fBm, for each $H\in\rbr{0,1/2}$
there exist some positive and finite constants $\bar{A}_{H},\bar{B}_{H}$,
depending on $H$ only, such that for any $s$, $t$, $c$ and $v$
such that $0\le s<t$, $0<c\le(t-s)^{H}$, $v\ge1$, 
\begin{align*}
 & \P\rbr{\left|\TTV{W^{H}}{[s,t]}c-\E\TTV{W^{H}}{[s,t]}c\right|\ge\rbr{t-s}c{}^{1-1/H}v}\\
 & \le\bar{A}_{H}\exp\rbr{-\bar{B}_{H}(t-s)c{}^{-1/H}\min\rbr{v^{1+2H},v^{2}}}\\
 & =\bar{A}_{H}\exp\rbr{-\bar{B}_{H}(t-s)c^{-1/H}v^{1+2H}}.
\end{align*}
Similarly, by Theorem \ref{theo-2} and stationarity of the increments of fBm, for each $H\in[1/2,1)$
there exist some universal positive and finite constants $\bar{A},\bar{B}>0$
such that
\begin{align*}
 & \P\rbr{\left|\TTV{W^{H}}{[s,t]}c-\E\TTV{W^{H}}{[s,t]}c\right|\ge\rbr{t-s}c{}^{1-1/H}v}\\
 & \le\bar{A}\exp\rbr{-\bar{B}\rbr{(t-s)c^{-1/H}}^{2-2H}v^{2}}.
\end{align*}
Both cases yield that for each $H\in\rbr{0,1}$, there exist$\bar{A}_{H},\bar{B}_{H}\in(0,+\ns)$
such that for $s$, $t$, $c$ and $v$ satisfying $0\le s<t$, $0<c\le(t-s)^{H}$,
$v\ge1$ one has 
\begin{align}
 & \P\rbr{\left|\TTV{W^{H}}{[s,t]}c-\E\TTV{W^{H}}{[s,t]}c\right|\ge\rbr{t-s}c{}^{1-1/H}v} \nonumber \\
 & \le\bar{A}_{H}\exp\rbr{-\bar{B}_{H}\rbr{(t-s)c^{-1/H}}^{2-\max\rbr{2H,1}}v^{1+\min(2H,1)}} \label{eq:Th12-1}.
\end{align}
An almost immediate consequence of inequality (\ref{eq:Th12-1}) is
the a.s. existence of the limits 
\begin{equation}
\lim_{c\ra0_{+}}c^{1/H-1}\TTV{W^{H}}{[s,t]}c=\mathfrak{c}_{H}(t-s),\label{eq:TTVlimit}
\end{equation}
\begin{equation}
\lim_{c\ra0_{+}}c^{1/H-1}\UTV{W^{H}}{[s,t]}c=\lim_{c\ra0_{+}}c^{1/H-1}\DTV{W^{H}}{[s,t]}c=\frac{1}{2}\mathfrak{c}_{H}(t-s).\label{eq:UTVlimit}
\end{equation}
\begin{coro} \label{cor_TV_limit}. For any $0\le s<t$,
the equalities (\ref{eq:TTVlimit}) and (\ref{eq:UTVlimit}) hold
a.s. with $\mathfrak{c}_{H}$ being the finite positive constant satisfying
for any $n\in\N$,
\[
\frac{\E\TTV{W^{H}}{[0,n]}1}{n}\le\mathfrak{c}_{H}\le\frac{\E\TTV{W^{H}}{[0,n]}1}{n}+\frac{1}{n}.
\]
\end{coro}
{\bf Proof.}
From (\ref{eq:sup_2})-(\ref{eq:sub_2}) we infer that the
sequence $a_{n}:=\E\TTV{W^{H}}{[0,n]}1$, $n\in\N$, is superadditive,
while the sequence $b_{n}:=\E\TTV{W^{H}}{[0,n]}1+1$, $n\in\N$, is
subadditive. As a result, by Fekete's lemma, there exists the limit
\[
\mathfrak{c}_{H}:=\lim_{n\ra0_{+}}\frac{\E\TTV{W^{H}}{[0,n]}1}{n}=\lim_{n\ra0_{+}}\frac{\E\TTV{W^{H}}{[0,n]}1+1}{n},
\]
moreover,
\[
\mathfrak{c}_{H}=\sup_{n\in\N}\frac{\E\TTV{W^{H}}{[0,n]}1}{n}=\inf_{n\in\N}\frac{\E\TTV{W^{H}}{[0,n]}1+1}{n}.
\]
The positivity of $\mathfrak{c}_{H}$ follows from the bound $\mathfrak{c}_{H}\ge\E\TTV{W^{H}}{[0,1]}1>0$
while its finiteness follows from the bound $\mathfrak{c}_{H}\le\E\TTV{W^{H}}{[0,1]}1+1$
and finiteness of $\E\TTV{W^{H}}{[0,1]}1$, which is an immediate
consequence of \cite[Corollary 2]{BednorzLochowski:2013}. For more
elementary proof of finiteness of $\E\TTV{W^{H}}{[0,1]}1$, based
on H\"older's continuity and Kolmogorov's continuity theorem see  \cite[Lemma 2.8]{Toyomu:2023}.

By scaling properties of fBm, for any $t>0$ we have
\begin{align*}
c^{1/H-1}\E\TTV{W^{H}}{[0,t]}c & =c^{1/H-1}\E\TTV{cW_{c^{-1/H}\cdot}^{H}}{[0,t]}c\\
 & =c^{1/H-1}c\E\TTV{W^{H}}{[0,c^{-1/H}t]}1\\
 & =\frac{\E\TTV{W^{H}}{[0,c^{-1/H}t]}1}{c^{-1/H}t}t
\end{align*}
thus 
\begin{equation}
\lim_{c\ra0_{+}}c^{1/H-1}\E\TTV{W^{H}}{[0,t]}c=\lim_{c\ra0_{+}}\frac{\E\TTV{W^{H}}{[0,c^{-1/H}t]}1}{c^{-1/H}t}t=\mathfrak{c}_{H}t.\label{eq:ETV_limit}
\end{equation}
Now, by (\ref{eq:Th12-1}) and the Borel-Cantelli lemma, the limit
\[
\lim_{n\ra+\ns}(1/n)^{1/H-1}\TTV{W^{H}}{[0,t]}{1/n}
\]
exists a.s. and is equal $\mathfrak{c}_{H}t$. Finally, for any sequence
of positive reals $\rbr{c_{n}}$ tending to $0$, for those $c_{n}$
which are smaller than $1$, we have $1/\left\lfloor 1/c_{n}\right\rfloor \ge c_{n}\ge1/\left\lceil 1/c_{n}\right\rceil $
thus (by the monotonicity of $c\mapsto\TTV{W^{H}}{[0,t]}c$ \textendash{}
recall the definition of $\TTV f{[0,t]}c$)
\[
\TTV{W^{H}}{[0,t]}{1/\left\lfloor 1/c_{n}\right\rfloor }\le\TTV{W^{H}}{[0,t]}{c_{n}}\le\TTV{W^{H}}{[0,t]}{1/\left\lceil 1/c_{n}\right\rceil }
\]
and
\begin{align}
\rbr{1/\left\lceil 1/c_{n}\right\rceil }^{^{1/H-1}}\TTV{W^{H}}{[0,t]}{1/\left\lfloor 1/c_{n}\right\rfloor } & \le c_{n}^{1/H-1}\TTV{W^{H}}{[0,t]}{c_{n}}\nonumber \\
 & \le\rbr{1/\left\lfloor 1/c_{n}\right\rfloor }^{^{1/H-1}}\TTV{W^{H}}{[0,t]}{1/\left\lceil 1/c_{n}\right\rceil }.\label{eq:3sequences}
\end{align}
Since $c_{n}\ra0_{+}$, 
\[
\lim_{n\ra+\ns}\rbr{\frac{1/\left\lceil 1/c_{n}\right\rceil }{1/\left\lfloor 1/c_{n}\right\rfloor }}^{1/H-1}=1
\]
and we have 
\begin{align*}
 & \lim_{n\ra+\ns}\rbr{1/\left\lfloor 1/c_{n}\right\rfloor }^{^{1/H-1}}\TTV{W^{H}}{[0,t]}{1/\left\lceil 1/c_{n}\right\rceil }\\
 & =\lim_{n\ra+\ns}\rbr{\frac{1/\left\lceil 1/c_{n}\right\rceil }{1/\left\lfloor 1/c_{n}\right\rfloor }}^{1/H-1}\rbr{1/\left\lfloor 1/c_{n}\right\rfloor }^{^{1/H-1}}\TTV{W^{H}}{[0,t]}{1/\left\lceil 1/c_{n}\right\rceil }\\
 & =\lim_{n\ra+\ns}\rbr{1/\left\lceil 1/c_{n}\right\rceil }^{^{1/H-1}}\TTV{W^{H}}{[0,t]}{1/\left\lceil 1/c_{n}\right\rceil }=\mathfrak{c}_{H}t\text{ a.s.}
\end{align*}
and similarly 
\[
\lim_{n\ra+\ns}\rbr{1/\left\lceil 1/c_{n}\right\rceil }^{^{1/H-1}}\TTV{W^{H}}{[0,t]}{1/\left\lfloor 1/c_{n}\right\rfloor }=\mathfrak{c}_{H}t\text{ a.s.},
\]
thus the middle term in (\ref{eq:3sequences}) tends a.s. to $\mathfrak{c}_{H}t$
as well. Finally, the convergence of $c^{1/H-1}\TTV{W^{H}}{[s,t]}c$
follows from the convergences of $c^{1/H-1}\TTV{W^{H}}{[0,t]}c$,
$c^{1/H-1}\TTV{W^{H}}{[0,s]}c$ and the estimates (see (\ref{eq:sup_2})
and (\ref{eq:sub_2})):
\[
\TTV{W^{H}}{[0,t]}c-\TTV{W^{H}}{[0,s]}c-c\le\TTV{W^{H}}{[s,t]}c\le\TTV{W^{H}}{[0,t]}c-\TTV{W^{H}}{[0,s]}c.
\]

The proof of existence of the limits (\ref{eq:UTVlimit})
is analogous. The existence of the limit 
\[
\lim_{n\ra0_{+}}\frac{\E\UTV{W^{H}}{[0,n]}1}{n}
\]
follows from Fekete's lemma, while the fact that its value equals
$\mathfrak{c}_{H}/2$ follows from symmetry of $W^{H}$ which yields
\[
\E\UTV{W^{H}}{[0,n]}1=\E\DTV{W^{H}}{[0,n]}1,\ n\in\N,
\]
and the equality 
\[
\E\TTV{W^{H}}{[0,n]}1=\E\UTV{W^{H}}{[0,n]}1+\E\DTV{W^{H}}{[0,n]}1.
\]
Next, instead of (\ref{eq:Th12-1}) one uses Theorems \ref{theo-11} and \ref{theo-21}.
\hfill $\blacksquare$

Taking into account that for $c\le(t-s)^{H}$,
\[
\E\TTV{W^{H}}{[s,t]}c\sim_{H}\rbr{t-s}c{}^{1-1/H}
\]
and that for $c\le(t-s)^{H}$, $(t-s)c^{-1/H}\ge1$, we get the following
simpler (but weaker) form of (\ref{eq:Th12-1}): for each $H\in\rbr{0,1}$,
there exist$\bar{A}_{H},\bar{B}_{H},\bar{D}_{H}\in(0,+\ns)$ such
that for $c\le(t-s)^{H}$, $v\ge1$
\begin{align}
 & \P\rbr{\TTV{W^{H}}{[s,t]}c\ge\bar{D}_{H}\rbr{t-s}c{}^{1-1/H}v}  \le\bar{A}_{H}\exp\rbr{-\bar{B}_{H}v^{1+\min(2H,1)}}. \label{eq:Th12-1-2}
\end{align}

\subsection{Proofs of Theorems \ref{theo-1} - \ref{theo-21}}
Our approach will be based on the following concentration result of Daniel J. Fresen, \cite{Fresen03072023}. Let $\psi: \R^M \ra \R$ be a continuously differentiable function, $\varphi: \R \ra \R$ be a convex function and $\rbr{g_i} = \rbr{g_i}_{i=1}^M$, $\rbr{g_i'} = \rbr{g_i'}_{i=1}^M$ be two independent, standard normal vectors in $\R^M$.
{ We have 
\be \label{Fresen}
\E \varphi\rbr{\psi\rbr{\rbr{g_i}} - \psi \rbr{\rbr{g_i'}}} \ls \E \varphi\rbr{\frac{\pi}{2}|\nabla \psi((g_i))| g' },
\ee
where $g'$ is a $1$-dimensional, $\ccN(0,1)$ random variable, idependent from $\rbr{g_i}$.}

{We will set $S=1$ and assume that $c \ls 1$ (the case of general $S>0$, $c \ls S^H$, will be obtained by scaling).}
Fix integer $N$ large enough.
For all $k\in \cbr{0,1,2,\ldots,2N}$ we assume  the following representation
\[
W^H_{k/(2N)}=\sum_i M_{ki}g_i,\;\; \mbox{where}\;\; g_i-\mbox{i.i.d.},\;\; \ccN(0,1),
\]
such that for any $k,l\in \cbr{0,1,\ldots,2N}$
\[
|(k-l)/(2N)|^{2H}=\E\rbr{W^H_{k/(2N)}-W^H_{l/(2N)}}^2=\sum_{i} \rbr{M_{ki}-M_{li}}^2. 
\]
Now, we define 
\[
F_{n,N} (x):=
\begin{cases} 0 & \text{ if } n=0, \\
\max_{0\ls l_1< k_1 < l_2<k_2 < \ldots < l_{n}<k_n\ls 2N}\sum^n_{j=1} \sum_i \rbr{M_{k_j i}-M_{l_j i}}x_i , & \text{ if }  n=1,2,\ldots, N;
\end{cases}
\]
and 
\[
F_N^c (x):=\max_{0\ls n\ls N } \rbr{F_{n,N} (x) - nc}. 
\]
For $n=0,1,2,\ldots, N$ we define
\begin{align*}
UTV_{n,N}(W^H) & := F_{n,N} ((g_i)) \\
& = \begin{cases} 0 & \text{ if } n=0, \\
  \max_{0\ls l_1< k_1 < l_2<k_2 < \ldots < l_{n}<k_n\ls 2N}\sum^n_{j=1}\rbr{W^H_{k_j /(2N)}-W^H_{l_j  /(2N) }} & \text{ if } n = 1,2,\ldots,N;
\end{cases}
\end{align*}
and
\[
UTV_N^c(W^H) : =F_N^c((g_i))=\max_{0\ls n\ls N } \rbr{UTV_{n,N}(W^H) - nc}.
\]
{By \eqref{Fresen} We have 
\[
\E\exp\rbr{\lambda (UTV_N^c-{UTV_N^c}') } \ls \E \exp\rbr{\frac{\pi}{2}\lambda |\nabla F_N^c((g_i))| g' },
\]
where $g'$ is a $1$-dimensional, $\ccN(0,1)$ random variable, idependent from $\rbr{g_i}$.}

Let 
\[
\tau_N^c(x) := \min \cbr{n = 0,1,\ldots,N: F_N^c(x) =F_{n,N} (x) - nc}.
\]
For any $n =0,1,2,\ldots $, the boundary of the set $\cbr{x: \tau_N^c(x) = n}$ has the Lebesgue measure $0$ (notice that the interior of $\cbr{x: \tau_N^c(x) = n}$ is obtained as finite number of  intersections and unions of regions where some affine functions of components of $x$ are strictly positive), therefore $(g_i)$ is a.s. inside the interior of one of the sets $\cbr{x: \tau_N^c(x) = n}$, $n=0,1,\ldots,N$. 
Obviously,  inside the interior of the set $\cbr{x: \tau_N^c(x)=n}$ we have
\begin{align}
|\nabla F_N^c(x)|^2  = |\nabla F_{n,N}(x)|^2 & \ls \max_{0\ls l_1< k_1 < l_2<k_2 < \ldots < l_{n}<k_n\ls 2N} \sum_i \rbr{\sum^n_{j=1} \rbr{ M_{k_j i}-M_{l_j i}}}^2 \nonumber \\
& = \max_{0\ls l_1< k_1 < l_2<k_2 < \ldots < l_{n}<k_n\ls 2N} \E \rbr{ \sum_i g_i  \sum^n_{j=1} \rbr{ M_{k_j i}-M_{l_j i}}}^2 \nonumber \\
& = \max_{0\ls l_1< k_1 < l_2<k_2 < \ldots < l_{n}<k_n\ls 2N} \E\rbr{\sum^n_{j=1} \rbr{W^H_{k_j/(2N)}-W^H_{l_j/(2N)}}}^2. \label{grad}
\end{align}
We consider two cases.  
\subsubsection{The case $H \in (0,1/2)$} \label{subsect:2.2.1}
If $H<1/2$ and $S=1$,  and then by \eqref{grad}, \eqref{weakvarhle12},
 \[
 |\nabla F_{n,N}((g_i))|^2 \ls n\rbr{\frac{S}{n}}^{2H} = n^{1-2H} \text{ a.e. }
 \]
Therefore, 
\begin{equation} \label{expansion}
\E\exp\rbr{\lambda (UTV_N^c-{UTV_N^c}') } \ls \E \exp\rbr{\frac{\pi}{2}\lambda |\nabla F_N^c \rbr{(g_i)}| g' }\ls \sum^N_{n = 0} \exp\rbr{\frac{\pi^2}{8}\lambda^2 n^{1-2H} }\P\rbr{\tau_N^c((g_i))=n}.
\end{equation}

Let ${UTV_N^c}'$ be an independent copy of ${UTV_N^c}$. We will prove the following
\begin{theo} \label{theorem1} Let $c \in (0,1]$, $H \in (0,1/2)$. There exist some positive constants $A_{1,H} \ls 9, A_{2,H} \ls 9$, $B_{1,H}, B_{2,H}$, $K_H$ (depending only on $H$) such that
\begin{itemize}
\item for $\lambda \gs K_H/c$
\[
\E\exp\rbr{\lambda (UTV_N^c-{UTV_N^c}') } \ls A_{1,H} \exp\rbr{B_{1,H} \lambda^{1+\frac{1}{2H}} c^{1-\frac{1}{2H}}};
\]
\item for $\lambda < K_H/c$
\[
\E\exp\rbr{\lambda (UTV_N^c-{UTV_N^c}') } \ls A_{2,H} \exp\rbr{B_{2,H}  \lambda^{2} c^{2-\frac{1}{H}}}.
\]
\end{itemize}
\end{theo}

\begin{coro} \label{corollary1} Let $c \in (0,1]$, $H \in (0,1/2)$.   There exist some positive constants  $A_H \ls 9$ and $B_{H}$ (depending only on $H$) such that
for any $u >0$,
\begin{align}
\E\exp\rbr{\lambda (UTV_N^c-{UTV_N^c}') } & \ls A_H \exp \rbr{B_H \max\cbr{ \lambda^{2} c^{2-\frac{1}{H}}, \lambda^{1+\frac{1}{2H}} c^{1-\frac{1}{2H}}}} \nonumber \\
& =  A_H \exp \rbr{B_H   \lambda^{2} c^{1-\frac{1}{2H}} \max\cbr{ \lambda, c^{-1}}^{\frac{1}{2H}-1}}.  \label{cor1}
\end{align}
\end{coro}
{\bf Proof of Corollary \ref{corollary1}.}
Let $A_{1,H}$, $A_{2,H}$, $B_{1,H}$, $B_{2,H}$ be as in Theorem \ref{theorem1}. Now \eqref{cor1} follows when we take $A_H = \max\cbr{A_{1,H}, A_{2,H}} \ls 9$, $B_H = \max \cbr{B_{1,H}, B_{2,H}}$.
\hfill $\blacksquare$

\begin{coro} \label{corollary2}
 Let $c \in (0,1]$, $H \in (0,1/2)$.  There exist positive constants ${\hat A}_H \ls 18$ and ${\hat B}_{H}$ (depending only on $H$) such that
for any $u >0$,
\[
\P \rbr{ \left| UTV_N^c-\E {UTV_N^c}\right| > u} \ls {\hat A}_H \exp\rbr{ - {\hat B}_H c^{\frac{1}{H} - 2} u^{1+2H} \min\rbr{c^{1-\frac{1}{H}}, u}^{1-2H} }
\]
or, equivalently, introducing $v >0$ such that $u = c^{1-\frac{1}{H}}v$,
\[
\P \rbr{ \left| UTV_N^c-\E {UTV_N^c}\right| > c^{1-\frac{1}{H}}v} \ls {\hat A}_H \exp\rbr{ - {\hat B}_H c^{-\frac{1}{H}} \min\rbr{v^2, v^{1+2H}} }.
\]
\end{coro}
{\bf Proof of Corollary \ref{corollary2}.}
Jensen's inequality for any real $\lambda$ yields:
\[
\exp\rbr{\lambda (UTV_N^c-\E {UTV_N^c}') } \ls E \rbr{ \exp\rbr{\lambda (UTV_N^c-{UTV_N^c}') } |  UTV_N^c},
\]
\[
 \E \exp\rbr{\lambda (UTV_N^c-\E {UTV_N^c}') } \ls \E  \exp\rbr{\lambda (UTV_N^c-{UTV_N^c}') }.
\]
Then, by Markov's inequality, \eqref{cor1} and optimization by $\lambda$ we have
\begin{align}
\P\rbr{UTV_N^c-\E {UTV_N^c} > u} & \ls \inf_{\lambda > 0} \rbr{ \E \exp\rbr{\lambda (UTV_N^c-\E {UTV_N^c}') } e^{-\lambda \cdot u}} \nonumber \\
& \ls  \inf_{\lambda > 0} \rbr{ \E \exp\rbr{\lambda (UTV_N^c-{UTV_N^c}') } e^{-\lambda \cdot u}} \nonumber \\
& \ls \inf_{\lambda > 0} \rbr{ A_H  \exp \rbr{B_H   \lambda^{2} c^{1-\frac{1}{2H}} \max\cbr{ \lambda, c^{-1}}^{\frac{1}{2H}-1}} e^{-\lambda \cdot u} }\nonumber \\
& =  A_H   \exp \rbr{ \inf_{\lambda > 0} \rbr{B_H   \lambda^{2} c^{1-\frac{1}{2H}} \max\cbr{ \lambda, c^{-1}}^{\frac{1}{2H}-1}-\lambda u}} \nonumber \\
& = 
\begin{cases}
A_H\exp\rbr{-\frac{1}{4B_H} c^{\frac{1}{H}-2} u^2} \text{ if } u \ls 2B_H c^{1-\frac{1}{H}}, \\
A_H \exp\rbr{B_H c^{-\frac{1}{H}} - {c}^{-1}u}  \text{ if } u \in \rbr{2B_H c^{1-\frac{1}{H}}, \rbr{1+\frac{1}{2H}}B_H c^{1-\frac{1}{H}}}, \\
A_H \exp\rbr{-\frac{(2H)^{2H}}{(2H+1)^{2H+1}B_H} c^{1-2H} u^{1+2H}} \text{ if } u \gs \rbr{1+\frac{1}{2H}}B_H c^{1-\frac{1}{H}} \\
\end{cases} \nonumber \\
& \ls A_H \exp \rbr{ - \min \rbr{\frac{1}{4B_H}, \frac{(2H)^{2H}}{(2H+1)^{2H+1}B_H}} \min\rbr{c^{\frac{1}{H}-2} u^2, c^{1-2H} u^{1+2H}} } \nonumber \\
& = A_H \exp\rbr{ - {\hat B}_H c^{\frac{1}{H} - 2} u^{1+2H} \min\rbr{c^{1-\frac{1}{H}}, u}^{1-2H} } \label{eq12}
\end{align}
with ${\hat B}_H =  \min \rbr{\frac{1}{4B_H}, \frac{(2H)^{2H}}{(2H+1)^{2H+1}B_H}}$.

Similarly,
\begin{align}
\P\rbr{ {UTV_N^c} - \E UTV_N^c < - u}  &= \P\rbr{\E UTV_N^c- {UTV_N^c} > u}  \nonumber \\
& \ls \inf_{\lambda > 0}  \rbr{\E \exp\rbr{\lambda (\E UTV_N^c- {UTV_N^c}') } e^{-\lambda \cdot u} }\nonumber \\
& \ls  \inf_{\lambda > 0} \rbr{ \E \exp\rbr{\lambda (UTV_N^c- {UTV_N^c}') } e^{-\lambda \cdot u} }\nonumber \\
& \ls A_H \exp\rbr{ - {\hat B}_H c^{\frac{1}{H} - 2} u^{1+2H} \min\rbr{c^{1-\frac{1}{H}}, u}^{1-2H} }. \label{eq13}
\end{align}
Now, by \eqref{eq12} and \eqref{eq13} the result follows with $\hat{A}_H = 2A_H$.
\hfill $\blacksquare$

{ 
\begin{rema} \label{optimality} The bound
\[
\P \rbr{ UTV^c_N-\E {UTV^c_N} > u} \ls {\hat A}_H \exp\rbr{ - {\hat B}_H c^{\frac{1}{H} - 2} u^{1+2H} \min\rbr{c^{1-\frac{1}{H}}, u}^{1-2H} }
\]
for $u \gs c^{1-{1}/{H}}$ takes the form.
\[
\P \rbr{ UTV^c_N-\E {UTV^c_N} > u} \ls {\hat A}_H \exp\rbr{ - {\hat B}_H c^{1 - 2H} u^{1+2H}}.
\]
The factor $u^{1+2H}$ in this bound is optimal in the sense that there exists positive constant $T_H$ such that for $u \gs c^{1-{1}/{H}}$,
\[
\P \rbr{ UTV^c_N-\E {UTV^c_N} > u} \gs  \frac{1}{2} \exp\rbr{ - T_H c^{1-2H} u^{1+2H} }
\]
thus $u^{1+2H}$ can not be replaced by $u$ raised to the higher power than $1+2H$. 
Indeed, let $c \in (0,1)$ and $u \gs c^{1-1/H}$. Taking $n = n_{c,u} := \lfloor c^{-1}u \rfloor \gs  1$ we have 
\[
S^{c,u} := \sum_{i=1}^{n_{c,u}} W^H_{(2i)/(2n_{c,u})} - W^H_{(2i-1)/(2n_{c,u})} \sim \ccN\rbr{0, Q_{H,c,u} n_{c,u}^{1-2H}},
\]
where, by \eqref{weakvarhle12}, $1/8 \ls Q_{H,c,u} \ls 1$.
We have
\[
UTV^c_N \ge S^{c,u} - n_{c,u}c, \ \E UTV^c_N \le C_H c^{1-\frac{1}{H}} \ls C_H u, \
n_{c,u}c \ls  u
\]
and thus, for $u \gs  c^{1-1/H}$ and $Z \sim \ccN\rbr{0, 1}$, $\Phi(z) = \P \rbr{Z \ls z}$,
\begin{align*}
\P \rbr{UTV^c_N \gs \E{UTV^c_N} + u}&  \gs \P\rbr{S^{c,u} \gs n_{c,u} c + \E{UTV_N^c} + u}\\
& \gs \P\rbr{S^{c,u} \gs R_H u} \\
& = \P\rbr{\frac{S^{c,u}}{Q_{H,c,u}^{1/2} n_{c,u}^{1/2-H}} \gs \frac{R_H  u}{Q_{H,c,u}^{1/2}  n_{c,u}^{1/2-H}}} \\
&= \P\rbr{Z \gs \frac{R_H }{Q_{H,c,u}^{1/2}}  n_{c,u}^{H-1/2}u} \gs 1-\Phi \rbr{\sqrt{8}{R_H }  n_{c,u}^{H-1/2}u} \\
& \gs \frac{1}{2} \exp \rbr{-T_H  n_{c,u}^{2H-1}u^2} \gs \frac{1}{2} \exp \rbr{-T_H  c^{1-2H}u^{2H+1}}.
\end{align*}
\end{rema}
}

\noindent
{\bf Proof of Theorem \ref{theorem1}.}
Let us define
\[
\widetilde{{UTV}}_N^c(W^H) : =\max_{1\ls n\ls N } \rbr{UTV_{n,N}(W^H) - nc}.
\]

Clearly, for $n=1,2,\ldots,N$ and any real $A_n$
\begin{align} \label{esttau}
 \P\rbr{\tau_N^c((g_i))=n}\ls \P\rbr{UTV_{n,N}(W^H) -nc > A_n }+\P\rbr{\widetilde{{UTV}}_N^c(W^H) \ls A_n}.
\end{align}
It is not difficult to observe that
{
\[
UTV_{N}^c(W^H) \ls \UTV{ W^H}{1}{c}.
\]
}
Now, using \eqref{eutv} we  observe that for $n =1,2,\ldots, N$,
\begin{align*}
\E UTV_{n,N}(W^H) & \ls \inf_{c \gs 1}  \rbr{\E UTV_{N}^c(W^H) + nc}  \ls \inf_{c \gs 1}  \rbr{\E \TTV{ W^H}{1}{c} + nc} \\ 
& \ls  \inf_{c\gs1}\rbr{C_H c^{1-\frac{1}{H}} + nc} \ls {C_H \rbr{n^{-H}}^{1-\frac{1}{H}} + n\cdot n^{-H}} = D_H n^{1-H},
\end{align*}
where $D_H : = C_H +1$.
Therefore, if $D_H n^{-H}\ls c/4$ (equivalently $n \gs n_{0,c,H} := \rbr{4D_H/c}^{1/H}$) 
and $A_{n}=-nc/4$ we have
\[
nc + A_n - med \rbr{UTV_{n,N}(W^H)} \gs nc + A_n - 2 \E UTV_{n,N}(W^H) \gs nc+A_{n}-2 D_H n^{1-H} \gs nc/4,
\]
where $med \rbr{UTV_{n,N}(W^H)}$ denotes the median of $UTV_{n,N}(W^H)$. The estimate of the weak variance of $UTV_{n,N}(W^H)$ is given in \eqref{weakvarhle12}. 
Thus the first summand on the right side of \eqref{esttau}, by Gaussian concentration inequality, see for example \cite[Theorem 5.4.3]{Marcus:2006fk}, may be estimated as
\begin{align*}
& \P \rbr{UTV_{n,N}(W^H) -nc > A_n} \\
& = \P\rbr{UTV_{n,N}(W^H) - med \rbr{UTV_{n,N}(W^H)} > nc + A_n - med \rbr{UTV_{n,N}(W^H)}} \\
& \ls \P\rbr{UTV_{n,N}(W^H) - med \rbr{UTV_{n,N}(W^H)}} \\
& \ls \exp\rbr{-\frac{(nc)^{2}}{32n^{1-2H}}} = \exp\rbr{-\frac{n^{1+2H}c^2}{32}}.
\end{align*}
Next we need to look for the estimate of $\P\rbr{\widetilde{{UTV}}_N^c(W^H) \ls A_n}$. It suffices
to observe that $\widetilde{{UTV}}_N^c(W^H) \gs W_{1}^{H}-W_{0}^{H}-c$. Therefore, 
\[
\P\rbr{\widetilde{{UTV}}_N^c(W^H) \ls A_n}\ls\P\rbr{W_{1}^{H}-W_{0}^{H}-c\ls-nc/4}\ls\exp\rbr{-\frac{(n-4)^{2}c^{2}}{32}}.
\]
From last two estimates and \eqref{esttau}, for $$n\gs n_{1,c,H} := \rbr{8D_H/c}^{1/H}  \gs \max\cbr{n_{0,c,H},8}$$ (notice that $D_H = C_H + 1 \gs 1$ and recall that $c \ls 1$) we have  $n-4 \gs n/2$ and 
\be \label{estptaun}
 \P\rbr{\tau_N^c((g_i))=n} \ls 2 \exp\rbr{-\frac{n^{1+2H}c^2}{128}}.
\ee

Now we will estimate the sum on the right hand side of \eqref{expansion}.

For $n\ls n_{1,c,H}$ we simply estimate $\exp\rbr{\rbr{\pi^2/8}\lambda^{2}n^{1-2H}} \ls \exp\rbr{\rbr{\pi^2/8}\lambda^{2}n_{1,c,H}^{1-2H}}$ and 
\[
\sum_{n=0}^{\min\cbr{n_{1,c,H}, N}}\exp\rbr{\frac{\pi^2}{8}\lambda^{2}n^{1-2H}}\P\rbr{\tau_N^c((g_i))=n}\ls\exp\rbr{\frac{\pi^2}{8}\lambda^{2}n_{1,c,H}^{1-2H}}.
\]
Clearly, 
\[
n_{1,c,H}^{1-2H}=L_{H}c^{2-1/H},
\]
where $L_H = \rbr{8D_H}^{1/H-2}$, thus 
\be \label{pierwsza_suma}
\sum_{n=0}^{\min\cbr{n_{1,c,H}, N}}\exp\rbr{\frac{\pi^2}{8}\lambda^{2}n^{1-2H}}\P\rbr{\tau=n}\ls\exp\rbr{ \frac{\pi^2}{8}L_H \lambda^{2}c^{2-1/H}}.
\ee

Next, by \eqref{estptaun}, for $n_2 \ge {n_{1,c,H}}$
\be \label{druga_suma}
 \sum_{n > n_2}^{N}\exp\rbr{\frac{\pi^2}{8}\lambda^{2}n^{1-2H}}\P\rbr{\tau=n} \ls  2 \sum_{n >{n_{2}}}^{\ns}\exp\rbr{\frac{\pi^2}{8}\lambda^{2}n^{1-2H}-\frac{1}{128}n^{1+2H}c^2}.
\ee
Let us introduce two auxiliary numbers. Let $n_{\max, \lambda, c, H}$ be the unique number at which the function
\[
E(n) := \frac{\pi^2}{8}\lambda^{2}n^{1-2H}-\frac{1}{128}n^{1+2H}c^2, \ n \in  [0,+\ns),
\]
attains its global maximum. We have
\[
n_{\max, \lambda, c, H} = \rbr{16\pi^2 \frac{1-2H}{1+2H}}^{\frac{1}{4H}} \rbr{\frac{\lambda}{c}}^{\frac{1}{2H}} = G_H \rbr{\frac{\lambda}{c}}^{\frac{1}{2H}}.
\]
Moreover,
\[
E\rbr{n_{\max, \lambda, c, H}} = E_H \lambda^{1+\frac{1}{2H}} c^{1-\frac{1}{2H}}
\]
for some positive constant $E_H$ depending on $H$ only. Next, let $n_{2, \lambda, c, H}$ be the number such that for $n \gs n_{2, \lambda, c, H}$
\[
E(n) \ls -\frac{1}{256}n^{1+2H}c^2.
\]
We have
\[
n_{2, \lambda, c, H} = \rbr{32\pi^2}^{\frac{1}{4H}}  \rbr{\frac{\lambda}{c}}^{\frac{1}{2H}} = J_H  \rbr{\frac{\lambda}{c}}^{\frac{1}{2H}} \gs n_{\max, \lambda, c, H}.
\]

We consider two cases. 

{\bf First case:} $n_{\max, \lambda, c, H} \gs n_{1, c, H}$. This case is equivalent with
\[
G_H \rbr{\frac{\lambda}{c}}^{\frac{1}{2H}} \gs \rbr{\frac{8D_H}{c}}^{\frac{1}{H}}
\]
or
\[
\lambda \gs \frac{K_H}{c}. 
\]
In this case we estimate
\[
\lambda^{2}c^{2-1/H} = \lambda^{1+\frac{1}{2H}}  \lambda^{1- \frac{1}{2H}} c^{2-\frac{1}{H}} \ls \lambda^{1+\frac{1}{2H}}  K_H^{1- \frac{1}{2H}} c^{\frac{1}{2H}-1} c^{2-\frac{1}{H}} = M_H  \lambda^{1+\frac{1}{2H}}  c^{1-\frac{1}{2H}}
\]
and, using \eqref{druga_suma} and \eqref{expansion}, we estimate
\begin{align}
& \sum^{N}_{n = 0} \exp\rbr{\frac{\pi^2}{8}\lambda^2 n^{1-2H} }\P\rbr{\tau_N^c((g_i))=n} \nonumber \\
&  \ls \sum^{\min\cbr{n_{2,\lambda,c,H},N}}_{n = 0} \exp\rbr{\frac{\pi^2}{8}\lambda^2 n^{1-2H} }\P\rbr{\tau_N^c((g_i))=n} +  2 \sum_{n >{n_{2,\lambda,c,H}}}^{\ns}\exp\rbr{\frac{\pi^2}{8}\lambda^{2}n^{1-2H}-\frac{1}{128}n^{1+2H}c^2} \nonumber \\
&  \ls \exp\rbr{\frac{\pi^2}{8} \lambda^2 n_{2,\lambda,c,H}^{1-2H} } +  2  \sum_{n >{n_{2, \lambda,c,H}}}^{\ns} \exp\rbr{-\frac{1}{256}n^{1+2H}c^2} \nonumber \\
&  = \exp\rbr{\frac{\pi^2}{8} \lambda^2 J_H^{1-2H} \lambda^{\frac{1}{2H}-1} c^{\frac{1}{2H}-1}} +  2  \sum_{n >{n_{2, \lambda,c,H}}}^{\ns} \exp\rbr{-\frac{1}{256}n^{1+2H}c^2}.  \label{firstcase} \end{align}
 
To estimate the sum $ \sum_{n >{n_{2, \lambda,c,H}}}^{\ns} \ldots$ we notice that in the first case 
$$n_{2,\lambda,c,H} \gs n_{\max,\lambda,c,H} \gs n_{1,c,H} \gs  \rbr{\frac{8}{c}}^{1/H} $$
hence for $n > n_{2, \lambda,c,H}$
\[
n^{1+2H} c^2 \gs n \rbr{\frac{8}{c}}^{{2H}/{H}} c^2= n \frac{8^2}{c^2} c^2 \gs 64 n
\]
and
\be \label{geom}
 \sum_{n >{n_{2, \lambda,c,H}}}^{\ns} \exp\rbr{-\frac{1}{256}n^{1+2H}c^2} \ls \sum_{n=1}^{\ns} e^{-n/4} <4 .
\ee 
Finally, from \eqref{firstcase} and \eqref{geom} we have 
\[
 \sum^N_{n = 0} \exp\rbr{\frac{\pi^2}{8} \lambda^2 n^{1-2H} }\P\rbr{\tau_N^c((g_i))=n} \ls A_{1,H} \exp\rbr{B_{1,H}  \lambda^{1+\frac{1}{2H}} c^{1-\frac{1}{2H}}} 
\]
for $A_{1,H} \gs 1+ 8$, $B_{1,H} \gs \rbr{\pi^2/8} J_H^{1-2H}$.

{\bf Second case:} $n_{\max, \lambda, c, H} < n_{1, c, H}$.

In this case
\[
\frac{J_H}{G_H} n_{1, c, H} > \frac{J_H}{G_H} n_{\max, \lambda, c, H} =  \frac{J_H}{G_H}G_H \rbr{\frac{\lambda}{c}}^{\frac{1}{2H}}  = J_H \rbr{\frac{\lambda}{c}}^{\frac{1}{2H}} = n_{2, \lambda, c, H}.
\]
We estimate \eqref{expansion}
\begin{align}
& \sum^N_{n = 0} \exp\rbr{\frac{\pi^2}{8}\lambda^2 n^{1-2H} }\P\rbr{\tau_N^c((g_i))=n} \nonumber \\
& \ls \sum_{n=0}^{\min\cbr{(J_H/G_H) n_{1,c,H}, N}}\exp\rbr{\frac{\pi^2}{8}\lambda^{2}n^{1-2H}}\P\rbr{\tau_N^c((g_i))=n}  + 2  \sum_{n >(J_H/G_H) n_{1,c,H}}^{\ns} \exp\rbr{-\frac{1}{256}n^{1+2H}c^2}, \nonumber
\end{align}
where in the second sum we used $n > \frac{J_H}{G_H} n_{1, c, H} > n_{2, \lambda, c, H}$ which allows to estimate $ \exp\rbr{\frac{1}{2}\lambda^2 n^{1-2H} }\P\rbr{\tau_N^c=n} $ by $\exp\rbr{-\frac{1}{256}n^{1+2H}c^2}$. 
The first sum is estimated as in \eqref{pierwsza_suma},
\begin{align*}
 \sum_{n=0}^{\min\cbr{(J_H/G_H) n_{1,c,H}, N}} & \exp\rbr{\frac{1}{2}\lambda^{2}n^{1-2H}}\P\rbr{\tau_N^c((g_i))=n} \\
& \ls \exp\rbr{\frac{\pi^2}{8}\lambda^{2}\frac{J_H^{1-2H}}{G_H^{1-2H}}n_{1,c,H}^{1-2H}} \ls  \exp\rbr{\frac{\pi^2}{8}\frac{J_H^{1-2H}}{G_H^{1-2H}}L_H \lambda^{2}c^{2-1/H}}.
\end{align*}
To estimate the second sum we use $n > (J_H/G_H) n_{1,c,H} \gs n_{1,c,H} \gs (8/c)^{1/H}$ and proceed as in \eqref{geom}  to obtain
\[
2  \sum_{n >(J_H/G_H) n_{1,c,H}}^{\ns} \exp\rbr{-\frac{1}{256}n^{1+2H}c^2} < 8.
\]
Finally, in the second case,
\[
 \sum^N_{n = 0} \exp\rbr{\frac{\pi^2}{8}\lambda^2 n^{1-2H} }\P\rbr{\tau_N^c((g_i))=n} \ls A_{2,H} \exp\rbr{B_{2,H}  \lambda^{2} c^{2-{1}/{H}}} 
\]
for $A_{2,H} \gs 1+ 8$, $B_{2,H} \gs \rbr{\pi^2/8}\rbr{J_H^{1-2H}/G_H^{1-2H}}L_H$.
\hfill $\blacksquare$

\noindent
{\bf Proof of Theorem \ref{theo-11}.} By a.s. continuity of trajectories of $W^H$, $UTV_N^c$ tends a.s. to $\UTV{W^H} {1}{c}$ as $N$ tends to $+\ns$, hence, by Corollary \ref{corollary2}, for $c \in (0,1]$
\be \label{cor2utvc}
\P \rbr{ \left| \UTV{W^H} {1}{c}-\E {\UTV{W^H} {1}{c}}\right| > u} \ls {\hat A}_H \exp\rbr{ - {\hat B}_H c^{\frac{1}{H} - 2} u^{1+2H} \min\rbr{c^{1-\frac{1}{H}}, u}^{1-2H} }.
\ee
Now, take any $S >0$ and $c \le S^H$. By self-similarity of $W^H$ we have
\[
\UTV{W^H}{S}{c} \sim \UTV{S^H \cdot W^H}{\rbr{S^H}^{-1/H}S}{c} = \UTV{S^H \cdot W^H}{1}{c} = S^H \cdot \UTV{W^H}{1}{S^{-H} c}.
\] 
This together with \eqref{cor2utvc} (substituting $S^{-H}c$ in place of $c$ and $S^{-H}u$ in place of $u$) yields
\begin{align*}
\P & \rbr{ \left| \UTV{W^H} {S}{c}-\E {\UTV{W^H} {S}{c}}\right| > u} \\
& = \P \rbr{ \left| S^H \cdot \UTV{W^H}{1}{S^{-H} c} -S^H \cdot \E \UTV{W^H}{1}{S^{-H} c}\right| > u} \\
& = \P \rbr{ \left| \UTV{W^H}{1}{S^{-H} c} - \E \UTV{W^H}{1}{S^{-H} c}\right| > S^{-H}u} \\
& \ls {\hat A}_H \exp\rbr{ - {\hat B}_H \rbr{S^{-H}c}^{\frac{1}{H} - 2} \rbr{S^{-H}u}^{1+2H} \min\rbr{\rbr{S^{-H}c}^{1-\frac{1}{H}}, {S^{-H}u}}^{1-2H} } \\
& = {\hat A}_H \exp\rbr{ - {\hat B}_H S^{-1}{c}^{\frac{1}{H} - 2} {u}^{1+2H} \min\rbr{S{c}^{1-\frac{1}{H}}, {u}}^{1-2H} }.
\end{align*}
Next, substituting $u = Sc^{1-\frac{1}{H}}v$ we get 
\begin{align*}
\P & \rbr{ \left| \UTV{W^H} {S}{c}-\E {\UTV{W^H} {S}{c}}\right| > Sc^{1-\frac{1}{H}}v} \\
& \ls {\hat A}_H \exp\rbr{ - {\hat B}_H S^{-1}{c}^{\frac{1}{H} - 2} \rbr{Sc^{1-\frac{1}{H}}v}^{1+2H} \min\rbr{S{c}^{1-\frac{1}{H}}, {Sc^{1-\frac{1}{H}}v}}^{1-2H} } \\
& = {\hat A}_H \exp\rbr{ - {\hat B}_H S{c}^{-\frac{1}{H}} \min\rbr{v^{1+2H} , v^2}}.
\end{align*}
\hfill $\blacksquare$
\subsubsection{The case $H \in [1/2,1)$} \label{subsect:2.2.2}
If $H \in [1/2, 1)$ and $S=1$,  and then by \eqref{grad}, \eqref{weakvarhge12},
 \[
 |\nabla F_{n,N}((g_i))|^2 \ls S^{2H} = 1 \text{ a.e.}
 \]
Therefore, 
\begin{equation} \label{expansion2}
\E\exp\rbr{\lambda (UTV_N^c-{UTV_N^c}') } \ls \E \exp\rbr{\frac{\pi}{2}\lambda |\nabla F_N^c| g' }\ls  \exp\rbr{\frac{\pi^2}{8}\lambda^2 }.
\end{equation}

\begin{theo} \label{theorem2} Let $c \in (0,1]$, $H \in [1/2)$.  There exist universal positive constants ${\hat A} \le 2$ and ${\hat B} \ge 2/\pi^2$  such that
for any $u >0$,
\[
\P \rbr{ \left| UTV_N^c-\E {UTV_N^c}\right| > u} \ls {\hat A} \exp\rbr{ - {\hat B} u^{2} }.
\]
\end{theo}
{\bf Proof of Theorem  \ref{theorem2}.} Similarly as in the proof of Corollary \ref{corollary2}, if ${UTV_N^c}'$ is independent copy of ${UTV_N^c}$, then for any real $\lambda$:
\[
 \E \exp\rbr{\lambda (UTV_N^c-\E {UTV_N^c}') } \ls \E  \exp\rbr{\lambda (UTV_N^c-{UTV_N^c}') }.
\]
Then, by Markov's inequality, \eqref{expansion2} and optimization by $\lambda$ we have
\begin{align}
\P\rbr{UTV_N^c-\E {UTV_N^c} > u} & \ls \inf_{\lambda > 0} \rbr{ \E \exp\rbr{\lambda (UTV_N^c-\E {UTV_N^c}') } e^{-\lambda \cdot u}} \nonumber \\
& \ls  \inf_{\lambda > 0} \rbr{ \exp\rbr{\frac{\pi^2}{8}\lambda^2 }e^{-\lambda \cdot u}}  = \exp\rbr{-\frac{2}{\pi^2}u^2} = \exp\rbr{-{\hat B}u^2}.
\end{align}
with ${\hat B}=\rbr{2/\pi^2} $.

Similarly, 
\[
 \E \exp\rbr{\lambda (\E UTV_N^c- {UTV_N^c}') } \ls \E  \exp\rbr{\lambda (UTV_N^c-{UTV_N^c}') }.
\]
and 
\begin{align}
\P\rbr{ {UTV_N^c} - \E UTV_N^c < - u}  &= \P\rbr{\E UTV_N^c- {UTV_N^c} > u}  \nonumber \\
& \ls \inf_{\lambda > 0}  \rbr{\E \exp\rbr{\lambda (\E UTV_N^c- {UTV_N^c}') } e^{-\lambda \cdot u} }\nonumber \\
& \ls  \inf_{\lambda > 0} \rbr{ \exp\rbr{\frac{\pi^2}{8}\lambda^2 } e^{-\lambda \cdot u} }\nonumber = \exp\rbr{-\frac{2}{\pi^2}u^2} = \exp\rbr{-{\hat B}u^2}.
\end{align}
\hfill $\blacksquare$

\noindent
{\bf Proof of Theorem \ref{theo-21}.} Similarly as in the proof of Theorem \ref{theo-11}, we use a.s. continuity of trajectories of fBm and Theorem \ref{theorem2} to obtain
\be
\P \rbr{ \left| \UTV{W^H} {1}{c}-\E {\UTV{W^H} {1}{c}}\right| > u} \ls {\hat A} \exp\rbr{ - {\hat B} u^{2} }.
\ee
Next, using self-similarity of fBm, we have
\begin{align*}
\P & \rbr{ \left| \UTV{W^H} {S}{c}-\E {\UTV{W^H} {S}{c}}\right| > u} \\
& = \P \rbr{ \left| S^H \cdot \UTV{W^H}{1}{S^{-H} c} -S^H \cdot \E \UTV{W^H}{1}{S^{-H} c}\right| > u} \\
& = \P \rbr{ \left| \UTV{W^H}{1}{S^{-H} c} - \E \UTV{W^H}{1}{S^{-H} c}\right| > S^{-H}u} \\
& \ls {\hat A} \exp\rbr{ - {\hat B} \rbr{S^{-H}u}^2} = {\hat A} \exp\rbr{ - {\hat B} S^{-2H}u^2}. 
\end{align*}
\hfill $\blacksquare$

\section{Concentration of $1/H$-variation of fBm along the Lebesgue partitions
and its a.s. convergence \label{sec:Concentration}}

\subsection{Numbers of strip and $c$-level crossings by a continuous function}

Now, for the continuous path $f:[0,+\ns)\ra\R$, $0\le s<t$ and $c>0$
we define the number $U_{s,t}\rbr{c,f}$ of times (the graph of) $f$
\emph{upcrosses} the strip $\cbr{(x,y)\in\R^{2}:y\in[0,c]}$. Formally,
we define it as the authors of \cite{Toyomu:2023}, as follows 
\[
U_{s,t}\rbr{c,f}:=\#\cbr{(u,v):s\le u<v\le t,f(u)=0,\ f(v)=c,\forall r\in(u,v)\ f(r)\in(0,c)},
\]
where $\#A$ denotes the cardinality of the set $A$. Similarly, one
may define the number $D_{s,t}\rbr{c,f}$ of times (the graph of)
$f$ \emph{downcrosses} the strip $\cbr{(x,y)\in\R^{2}:y\in[0,c]}$
and the number $N_{s,t}\rbr{c,f}$ of times (the graph of) $f$ \emph{crosses}
the strip $\cbr{(x,y)\in\R^{2}:y\in[0,c]}$:
\[
N_{s,t}\rbr{c,f}=U_{s,t}\rbr{c,f}+D_{s,t}\rbr{c,f}.
\]
The just defined numbers are superadditive functions of the interval
$[s,t]$, that is, for any $u\in(s,t)$,
\begin{equation}
U_{s,t}(c,f)\ge U_{s,u}(c,f)+U_{u,t}(c,f)\label{eq:sup-1}
\end{equation}
and analogous relationships hold for $D_{\cdot,\cdot}\rbr{c,f}$ and
$N_{\cdot,\cdot}\rbr{c,f}$\@. On the other hand, $U_{\cdot,\cdot}(c,f)+1$
is a subadditive function of the interval $[s,t]$, that is, for any
$u\in(s,t)$,
\begin{equation}
U_{s,t}(c,f)+1\le U_{s,u}(c,f)+1+U_{u,t}(c,f)+1\label{eq:sub-1}
\end{equation}
and analogous relationships hold for $D_{\cdot,\cdot}\rbr{c,f}$ and
$N_{\cdot,\cdot}\rbr{c,f}$\@. 

Let us introduce three more related quantities called \emph{the total
numbers of $c$-level (up-, down-) crossings by $f$}. The total numbers
of $c$-level upcrossings and downcrossings by $f$ are defined respectively
as 
\[
KU_{s,t}(c,f):=\sum_{p\in\Z}U_{s,t}\rbr{c,f-p\cdot c},\quad KD_{s,t}(c,f):=\sum_{p\in\Z}D_{s,t}\rbr{c,f-p\cdot c}
\]
while the total numbers of $c$-level crossings by $f$ is respectively
defined as 
\begin{equation}
K_{s,t}(c,f):=\sum_{p\in\Z}N_{s,t}\rbr{c,f-p\cdot c}.\label{eq:Leb_var}
\end{equation}
Immediately from the definition it follows that for any $q\in\Z$,
$K_{s,t}(c,f+q\cdot c)=K_{s,t}(c,f)$. It is not difficult to see
that $K_{s,t}(c,f)$ is a superadditive function of the interval $[s,t]$,
that is, for any $u\in(s,t)$,
\begin{equation}
K_{s,t}(c,f)\ge K_{s,u}(c,f)+K_{u,t}(c,f).\label{eq:sup}
\end{equation}
On the other hand, $K_{s,t}(c,f)+1$ is a subadditive function of
the interval $[s,t]$, that is, for any $u\in(s,t)$,
\begin{equation}
K_{s,t}(c,f)+1\le K_{s,u}(c,f)+1+K_{u,t}(c,f)+1.\label{eq:sub}
\end{equation}
 It is also not difficult to see that $c^{1/H}K_{s,t}(c,f)$ is the
$1/H$- variation of $f$ along the Lebesgue partition of $[s,t]$
for the grid $c\Z$. Similarly, for $\rho\in\R$ $c^{1/H}K_{s,t}(c,f-\rho)$
is the $1/H$- variation of $f-\rho$ along the Lebesgue partition
for the grid $c\Z$ or $1/H$ - variation of $f$ along the Lebesgue
partition for the grid $c\Z+\rho$.

The numbers $U_{s,t}(c,f)$, $D_{s,t}(c,f)$, $N_{s,t}(c,f)$ are
very much related to the truncated variations via relationships
\begin{equation}
\UTV f{[s,t]}c=\int_{\R}U_{s,t}\rbr{c,f-a}\dd a,\quad\DTV f{[s,t]}c=\int_{\R}D_{s,t}\rbr{c,f-a}\dd a\label{eq:Banach_ind_gen_0}
\end{equation}
and
\begin{equation}
\TTV f{[s,t]}c=\int_{\R}N_{s,t}\rbr{c,f-a}\dd a\label{eq:Banach_ind_gen}
\end{equation}
($f-a$ denotes here the function $g:[0,+\ns)\ra\R$ such that $g(u)=f(u)-a$
for $u\in[0,+\ns)$), see \cite{LochowskiColloquium:2017}. From (\ref{eq:Banach_ind_gen})
and the definition of $K_{s,t}$ we also have
\begin{align}
\TTV f{[s,t]}c & =\int_{\R}N_{s,t}\rbr{c,f-a}\dd a=\sum_{p\in\Z}\int_{pc-c/2}^{pc+c/2}N_{s,t}\rbr{c,f-a}\dd a\nonumber \\
 & =\sum_{p\in\Z}\int_{-c/2}^{c/2}N_{s,t}\rbr{c,f-\rho-p\cdot c}\dd{\rho}=\int_{-c/2}^{c/2}\sum_{p\in\Z}N_{s,t}\rbr{c,f-\rho-p\cdot c}\dd{\rho}\nonumber \\
 & =\int_{-c/2}^{c/2}K_{s,t}(c,f-\rho)\dd{\rho}\label{eq:Banach_ind}
\end{align}
and similar relationships hold between pairs $\UTV f{[s,t]}c$ and
$KU_{s,t}\rbr{c,f-a}$ and $\DTV f{[s,t]}c$ and $KD_{s,t}\rbr{c,f-a}$. The
relationship (\ref{eq:Banach_ind}) is the reason for introduction
of one more quantity - $\bar{K}_{s,t}(c,f)$ - defined in the following
way
\begin{equation}
\bar{K}_{s,t}(c,f):=\frac{1}{c}\TTV f{[s,t]}c=\frac{1}{c}\int_{-c/2}^{c/2}K_{s,t}(c,f-\rho)\dd{\rho}.\label{eq:Banach_ind_1}
\end{equation}
It is worth to mention here that similar relationships to (\ref{eq:sup})-(\ref{eq:sub})
hold for the truncated variations, and they have the form: for any
$u\in(s,t)$,
\begin{equation}
\TTV f{[s,t]}c\ge\TTV f{[s,u]}c+\TTV f{[u,t]}c,\label{eq:sup_2}
\end{equation}
\begin{equation}
\TTV f{[s,t]}c\le\TTV f{[s,u]}c+\TTV f{[u,t]}c+c;\label{eq:sub_2}
\end{equation}
and analogous relationships hold for $\UTV f{[\cdot,\cdot]}c$ and
$\DTV f{[\cdot,\cdot]}c$. They follow from (\ref{eq:sup})-(\ref{eq:sub})
and (\ref{eq:Banach_ind}).

In the previous section we established strong concentration of the
truncated variation of fBm around its mean. In this section, using
the relationship between $K_{s,t}\rbr{c,W^{H}-\rho}$ and $\TTV{W^{H}}{[s,t]}c$
we will investigate the concentration of the $1/H$ - variation of
$W^{H}$ along the Lebesgue partition for the grid $c\Z+\rho$.

\subsection{Tails of total numbers of level crossings by fBm}

We start with estimates of $\E\sup_{\rho\in\R}K_{s,t}\rbr{c,W^{H}-\rho}$
and tails of $\sup_{\rho\in\R}K_{s,t}\rbr{c,W^{H}-\rho}$. These estimates
will follow from the estimates of tails of $\TTV{W^{H}}{[s,t]}c$. 

Let us notice that if $K_{s,t}\rbr{c,f-\rho}=k>0$ then there exist
$k+1$ times $s\le t_{1}<t_{2}<\ldots<t_{k+1}\le t$ such that $\left|f\rbr{t_{i+1}}-f\rbr{t_{i}}\right|=c$
for $i=1,2,\ldots,k$, hence 
\[
\TTV f{[s,t]}{c/2}\ge\sum_{i=1}^{k}\rbr{\left|f\rbr{t_{i+1}}-f\rbr{t_{i}}\right|-\frac{c}{2}}_{+}\ge k\frac{c}{2}=\frac{c}{2}K_{s,t}\rbr{c,f-\rho}.
\]
As a result we have the estimate
\begin{equation}
\sup_{\rho\in\R}K_{s,t}\rbr{c,f-\rho}=\sup_{\rho\in(-c/2,c/2]}K_{s,t}\rbr{c,f-\rho}\le\frac{2}{c}\TTV f{[s,t]}{c/2}.\label{eq:oszac_Kst}
\end{equation}
On the other side, for any $\rho\in\R$, 
\begin{equation}
K_{s,t}\rbr{c,f-\rho}\ge\frac{1}{c}\TTV f{[s,t]}{2c}.\label{eq:oszac_Kst_below}
\end{equation}
To prove this estimate let us take $\varepsilon>0$ and times $s\le t_{1}<t_{2}<\ldots<t_{k+1}\le t$
such that 
\[
\TTV f{[s,t]}{2c}\le\sum_{i=1}^{k}\rbr{\left|f\rbr{t_{i+1}}-f\rbr{t_{i}}\right|-2c}_{+}+\varepsilon.
\]
Let $I$ denote the set of indices $i\in\cbr{1,2,\ldots,k}$ such
that$\rbr{\left|f\rbr{t_{i+1}}-f\rbr{t_{i}}\right|-2c}_{+}\neq0$.
If $i\in I$ then $\Delta_{i}:=\left|f\rbr{t_{i+1}}-f\rbr{t_{i}}\right|>2c$
and there are at least $\lceil\Delta_{i}/c\rceil-1\ge2$ numbers from
the grid $c\Z+\rho$ in the interval $\rbr{f\rbr{t_{i}}\wedge f\rbr{t_{i+1}},f\rbr{t_{i}}\vee f\rbr{t_{i+1}}}$
corresponding to at least $\lceil\Delta_{i}/c\rceil-2\ge1$ upcrossing(s)
or downcrossing(s) by $f$ the strips of the form $\cbr{(x,y)\in\R^{2}:y\in\sbr{p\cdot c+\rho,(p+1)c+\rho}}$,
$p\in\Z$, during the time interval $\rbr{t_{i},t_{i+1}}$. As a result
we obtain
\[
K_{s,t}\rbr{c,f-\rho}\ge\sum_{i\in I}\rbr{\left\lceil \frac{\Delta_{i}}{c}\right\rceil -2}_{+}\ge\frac{1}{c}\sum_{i\in I}\rbr{\Delta_{i}-2c}_{+}\ge\frac{1}{c}\rbr{\TTV f{[s,t]}{2c}-\varepsilon}.
\]
Since $\varepsilon$ may be arbitrary close to $0$, we get (\ref{eq:oszac_Kst_below}).

Now, using (\ref{eq:oszac_Kst}), (\ref{eq:oszac_Kst_below}) and
the fact that for $c\in\left(0,(t-s)^{H}\right]$ (see Appendix)
\begin{equation}
\E\TTV{W^{H}}{[s,t]}c\sim_{H}(t-s)c^{1-1/H}.\label{eq:oszac_ETV}
\end{equation}
we infer that for $c>0$ such that $2c\le(t-s)^{H}$
\begin{equation}
\E\sup_{\rho\in\R}K_{s,t}(c,W^{H}-\rho)\sim_{H}(t-s)c^{-1/H}.\label{eq:oszac_EKst}
\end{equation}
Finally, using (\ref{eq:oszac_Kst}) and Theorems \ref{theo-1}
and \ref{theo-2}, we will estimate the tails of $\sup_{\rho\in\R}K_{s,t}\rbr{c,f-\rho}$.
In order to obtain these estimates, let us notice that by 
Theorem \ref{theo-1}, for each $H\in\rbr{0,1/2}$, for $c\le(t-s)^{H}/2$,
$v\ge1$ one has $(t-s)(c/2)^{-1/H}\ge2^{1/H}\ge1$ and for some $\bar{C}_{H}\in(1,+\ns)$
such that 
\[
\E\TTV{W^{H}}{[s,t]}{c/2}+\rbr{t-s}(c/2)^{-1/H}\le\bar{C}_{H}\rbr{t-s}c^{-1/H}
\]
 we have

\begin{align*}
 & \P\rbr{\sup_{\rho\in\R}K_{s,t}\rbr{c,W^{H}-\rho}\ge\bar{C}_{H}\rbr{t-s}c^{-1/H}v}\\
 & \le\P\rbr{\sup_{\rho\in\R}K_{s,t}\rbr{c,W^{H}-\rho}\ge\frac{2}{c}\E\TTV{W^{H}}{[s,t]}{c/2}+\rbr{t-s}(c/2)^{-1/H}v}\\
 & \le\P\rbr{\frac{2}{c}\TTV{W^{H}}{[s,t]}{c/2}\ge\frac{2}{c}\E\TTV{W^{H}}{[s,t]}{c/2}+\rbr{t-s}(c/2)^{-1/H}v}\\
 & =\P\rbr{\TTV{W^{H}}{[s,t]}{c/2}-\E\TTV{W^{H}}{[s,t]}{c/2}\ge\rbr{t-s}(c/2)^{1-1/H}v}\\
 & \le\bar{A}_{H}\exp\rbr{-\bar{B}_{H}(t-s)(c/2)^{-1/H}\min\rbr{v^{1+2H},v^{2}}}\\
 & \le\bar{A}_{H}\exp\rbr{-\bar{B}_{H}\min\rbr{v^{1+2H},v^{2}}}=\bar{A}_{H}\exp\rbr{-\bar{B}_{H}v^{1+2H}}.
\end{align*}
Similarly, by Theorem \ref{theo-2}, for each $H\in[1/2,1)$,
for some $\bar{C}\in(1,+\ns)$, $c\le(t-s)^{H}/2$, $v\ge1$ one has
\begin{align*}
 & \P\rbr{\sup_{\rho\in\R}K_{s,t}\rbr{c,W^{H}-\rho}\ge\bar{C}\rbr{t-s}c^{-1/H}v}\\
 & \le\P\rbr{\sup_{\rho\in\R}K_{s,t}\rbr{c,W^{H}-\rho}\ge\frac{2}{c}\E\TTV{W^{H}}{[s,t]}{c/2}+\rbr{t-s}(c/2)^{-1/H}v}\\
 & \le\P\rbr{\TTV{W^{H}}{[s,t]}{c/2}-\E\TTV{W^{H}}{[s,t]}{c/2}\ge\rbr{t-s}(c/2)^{-1/H}v}\\
 & \le\bar{A}\exp\rbr{-\bar{B}\rbr{(t-s)(c/2)^{-1/H}}^{2-2H}v^{2}}\le\bar{A}\exp\rbr{-\bar{B}v^{2}}.
\end{align*}
Both cases yield that for each $H\in\rbr{0,1}$, there exist$\bar{A}_{H},\bar{B}_{H},\bar{C}_{H}\in(0,+\ns)$
such that for $c\le(t-s)^{H}/2$, $v\ge1$ one has 
\begin{equation}
\P\rbr{\sup_{\rho\in\R}K_{s,t}\rbr{c,W^{H}-\rho}\ge\bar{C}_{H}\rbr{t-s}c^{-1/H}v}\le\bar{A}_{H}\exp\rbr{-\bar{B}_{H}v^{1+\min(2H,1)}}.\label{eq:Th12}
\end{equation}

\subsection{Decomposition of numbers of level crossings\label{subsec:Decomposition}}

A classical result by Mandelbrot and van Ness \cite{MandelbrotVanNess} gives the
following representation of $W^{H}$:
\begin{equation}
W_{u}^{H}=\frac{1}{C(H)}\int_{-\ns}^{0}\cbr{(u-r)^{H-1/2}-(-r)^{H-1/2}}\dd{W_{r}}+\frac{1}{C(H)}\int_{0}^{u}(u-r)^{H-1/2}\dd{W_{r}},\label{eq:vN_M}
\end{equation}
where $W_{r}$, $-\ns<r<+\ns$, is a symmetric standard Brownian motion
and
\begin{align*}
C(H) & =\cbr{\int_{-\ns}^{0}\rbr{(1-r)^{H-1/2}-(-r)^{H-1/2}}^{2}\dd r+\int_{0}^{1}(1-r)^{2H-1}\dd r}^{1/2}\\
 & =\cbr{\int_{0}^{+\ns}\rbr{(1+r)^{H-1/2}-r{}^{H-1/2}}^{2}\dd r+\frac{1}{2H}}^{1/2}.
\end{align*}
Let us fix $\varepsilon>0$. To obtain quantitative bounds on the
probability 
\[
\P\rbr{\left|c^{1/H}K_{0,1}\rbr{c,W^{H}-\rho}-c^{1/H}\E\bar{K}_{0,1}\rbr{c,W^{H}}\right|>\varepsilon}
\]
we will decompose $K_{0,1}\rbr{c,W^{H}-\rho}$ in a similar way as
it is done in stochastic sewing lemma setting, but additionally we
will take advantage of (\ref{eq:Th12}) and some concentration results
for martingales. Let $L$ be an integer greater or equal $2$ and
consider the decomposition of the segment $[0,1]$:
\[
[0,1]=\sbr{0,\frac{1}{L}}\cup\sbr{\frac{1}{L},\frac{2}{L}}\cup\ldots\cup\sbr{\frac{L-1}{L},1}.
\]
Next, each of the segments $\sbr{(l-1)/L,l/L}$, $l=1,2,\ldots,L$,
is further decomposed as 
\[
\sbr{\frac{l-1}{L},\frac{l}{L}}=\sbr{\frac{l-1}{L},\frac{l-1}{L}+\frac{1}{L^{2}}}\cup\sbr{\frac{l-1}{L}+\frac{1}{L^{2}},\frac{l-1}{L}+\frac{2}{L^{2}}}\cup\ldots\cup\sbr{\frac{l-1}{L}+\frac{L-1}{L^{2}},\frac{l}{L}}.
\]
For fixed $m=1,2,\ldots,L$ we consider the segments
\[
\sbr{\frac{m-1}{L^{2}},\frac{m}{L^{2}}},\sbr{\frac{1}{L}+\frac{m-1}{L^{2}},\frac{1}{L}+\frac{m}{L^{2}}},\ldots,\sbr{\frac{L-1}{L}+\frac{m-1}{L^{2}},\frac{L-1}{L}+\frac{m}{L^{2}}},
\]
define 
\begin{equation}
[s(l,m;L),t(l,m;L)]:=\sbr{\frac{l-1}{L}+\frac{m-1}{L^{2}},\frac{l-1}{L}+\frac{m}{L^{2}}},\ l=1,2,\ldots,L,\label{eq:stlm_def}
\end{equation}
and then define the martingale ${\cal M}^{m}=\rbr{{\cal M}_{1}^{m},\ {\cal M}_{2}^{m},\ldots,\ {\cal M}_{L}^{m}}$
as
\begin{align*}
{\cal M}_{n}^{m} & :=c^{1/H}\sum_{l=2}^{n}\cbr{K_{s(l,m;L),t(l,m;L)}\rbr{c,W^{H}-\rho}-E\sbr{K_{s(l,m;L),t(l,m;L)}\rbr{c,W^{H}-\rho}|{\cal F}_{t(l-1,m;L)}}},
\end{align*}
where ${\cal F}_{t}=\sigma\rbr{W_{r},r\le t}$ is the natural filtration
of the symmetric standard Brownian motion $W$ appearing in (\ref{eq:vN_M}).
(By the convention $\sum_{l=2}^{1}\ldots=0$, we have ${\cal M}_{1}^{m}=0$.)

We have the following estimates stemming from the superadditivity
(\ref{eq:sup}) of $K_{s,t}\rbr{\cdot,\cdot}$ and subadditivity (\ref{eq:sub})
of $K_{s,t}\rbr{\cdot,\cdot}+1$ as a function of the interval $[s,t]$:
\begin{equation}
K_{0,1}\rbr{c,W^{H}-\rho}\ge K_{0,1/L}\rbr{c,W^{H}-\rho}+\sum_{l=2}^{L}\sum_{m=1}^{L}K_{s(l,m;L),t(l,m;L)}\rbr{c,W^{H}-\rho},\label{eq:supp}
\end{equation}
\begin{equation}
K_{0,1}\rbr{c,W^{H}-\rho}\le K_{0,1/L}\rbr{c,W^{H}-\rho}+\sum_{l=2}^{L}\sum_{m=1}^{L}K_{s(l,m;L),t(l,m;L)}\rbr{c,W^{H}-\rho}+L^{2}.\label{eq:subb}
\end{equation}
Next, we decompose the sum $\sum_{l=2}^{L}\sum_{m=1}^{L}K_{s(l,m;L),t(l,m;L)}\rbr{c,W^{H}-\rho}$:
\begin{align}
 & c^{1/H}\sum_{l=2}^{L}\sum_{m=1}^{L}c^{1/H}K_{s(l,m;L),t(l,m;L)}\rbr{c,W^{H}-\rho} \nonumber \\
 & =\sum_{m=1}^{L}{\cal M}_{L}^{m}+\sum_{m=1}^{L}\sum_{l=2}^{n}c^{1/H}E\sbr{K_{s(l,m;L),t(l,m;L)}\rbr{c,W^{H}-\rho}|{\cal F}_{t(l-1,m;L)}}. \label{eq:dekompozycja}
\end{align}

\subsection{Estimates of tails of ${\cal M}_{L}^{m}$\label{subsec:Martingale}}

Let us recall (\ref{eq:stlm_def}) and notice that $t(l,m;L)-s(l,m;L)=L^{-2}$.
This observation together with (\ref{eq:Th12}) gives for $c\le L^{-2H}/2$
and $v\ge1$
\[
\P\rbr{c^{1/H}K_{s(l,m;L),t(l,m;L)}\rbr{c,W^{H}-\rho}\ge\bar{C}_{H}L^{-2}v}\le\bar{A}_{H}\exp\rbr{-\bar{B}_{H}v^{1+\min(2H,1)}}.
\]
Now, by Lemma B (see Appendix), for any $c\le L^{-2H}/2$, $x\ge0$, $v>0$,
\[
\P\rbr{\left|{\cal M}_{L}^{m}\right|\ge x}\le2\exp\rbr{-\frac{L^{3}x^{2}}{10\rbr{\bar{C}_{H}v}^{2}}}+6L\tilde{A}_{H}\exp\rbr{-\frac{\bar{B}_{H}}{2}v^{1+\min(2H,1)}}.
\]
Finally, taking $\varepsilon=x\cdot L$ we get 
\begin{align*}
 & \P\rbr{\left|\sum_{m=1}^{L}{\cal M}_{L}^{m}\right|\ge\varepsilon}\le\P\rbr{\sum_{m=1}^{L}\left|{\cal M}_{L}^{m}\right|\ge\varepsilon}\le\sum_{m=1}^{L}\P\rbr{\left|{\cal M}_{L}^{m}\right|\ge\varepsilon/L}\\
 & \le2L\exp\rbr{-\frac{L\varepsilon^{2}}{10\rbr{\bar{C}_{H}v}^{2}}}+6L^{2}\tilde{A}_{H}\exp\rbr{-\frac{\bar{B}_{H}}{2}v^{1+\min(2H,1)}}.
\end{align*}
In particular, taking $\varepsilon=L^{-1/4}$, $v=L^{1/6}$ we get
for any $c\le L^{-2H}/2$, 
\begin{align}
 & \P\rbr{\left|\sum_{m=1}^{L}{\cal M}_{L}^{m}\right|\ge L^{-1/4}}\nonumber \\
 & \le2L\exp\rbr{-\frac{L^{1/6}}{10\bar{C}_{H}^{2}}}+6L^{2}\tilde{A}_{H}\exp\rbr{-\frac{\bar{B}_{H}}{2}L^{(1+\min(2H,1))/6}}\nonumber \\
 & \le\rbr{1+6\tilde{A}_{H}}L^{2}\exp\rbr{-\min\rbr{\frac{1}{10\bar{C}_{H}^{2}},\frac{\bar{B}_{H}}{2}}L^{1/6}}.\label{eq:estim_mart_tails}
\end{align}

\subsection{Estimates of tails of $c^{1/H}E\protect\sbr{K_{s(l,m;L),t(l,m;L)}\protect\rbr{c,W^{H}-\rho}|{\cal F}_{t(l-1,m;L)}}$}

To estimate the tails of
\[
c^{1/H}E\sbr{K_{s(l,m;L),t(l,m;L)}\rbr{c,W^{H}-\rho}|{\cal F}_{t(l-1,m;L)}}
\]
 we decompose 
\begin{align}
 & c^{1/H}E\sbr{K_{s(l,m;L),t(l,m;L)}\rbr{c,W^{H}-\rho}|{\cal F}_{t(l-1,m;L)}}\nonumber \\
 & =c^{1/H}E\sbr{K_{s(l,m;L),t(l,m;L)}\rbr{c,W^{H}-\rho}-\bar{K}_{s(l,m;L),t(l,m;L)}\rbr{c,W^{H}}|{\cal F}_{t(l-1,m;L)}}\label{eq:skl1}\\
 & \quad+c^{1/H}E\sbr{\bar{K}_{s(l,m;L),t(l,m;L)}\rbr{c,W^{H}}|{\cal F}_{t(l-1,m;L)}}.\label{eq:skl2}
\end{align}
The tails of the term (\ref{eq:skl1}) will be estimated in the spirit
of \cite[Theorem A.1]{Picard:2008}, to be more precise, we will use \cite[Lemma 2.17]{Toyomu:2023}.

Let us fix $l=2,3,\ldots,L$, $m=1,2,\ldots,L$. Let us recall (\ref{eq:stlm_def})
and denote
\[
v=t(l-1,m;L)=\frac{l-2}{L}+\frac{m}{L^{2}},\quad s=s(l,m;L)=\frac{l-1}{L}+\frac{m-1}{L^{2}},\quad t=t(l,m;L)=\frac{l-1}{L}+\frac{m}{L^{2}}.
\]
We have 
\[
(s-v)^{-1}(t-v)^{1-H}=\rbr{\frac{1}{L}-\frac{1}{L^{2}}}^{-1}\rbr{\frac{1}{L}}^{1-H}\le2L^{H}.
\]
Now, using \cite[Lemma 2.17]{Toyomu:2023} for just defined $v$, $s$ and $t$,
we obtain for any $\rho'\in[-c/2,c/2]$ (with probability $1$):
\begin{align*}
 & \left|E\sbr{K_{s,t}\rbr{c,W^{H}-\rho}-K_{s,t}\rbr{c,W^{H}-\rho-\rho'}|{\cal F}_{v}}\right|\\
 & \le e^{b_{H}4L^{2H}|\rho'|^{2}}E\sbr{K_{s,t}\rbr{c,W^{H}-\rho}^{2}|{\cal F}_{v}}^{1/2}b_{H}^{1/2}2L^{H}|\rho'|\\
 & \le e^{b_{H}L^{2H}c^{2}}E\sbr{\sup_{\sigma}K_{s,t}\rbr{c,W^{H}-\sigma}^{2}|{\cal F}_{v}}^{1/2}b_{H}^{1/2}L^{H}c,
\end{align*}
where $b_{H}$ is a positive constant depending on $H$ only. Hence
(also with probability $1$)
\begin{align}
 & \left|E\sbr{K_{s,t}\rbr{c,W^{H}-\rho}-\bar{K}_{s,t}\rbr{c,W^{H}}|{\cal F}_{v}}\right|=\left|E\sbr{K_{s,t}\rbr{c,W^{H}-\rho}-\bar{K}_{s,t}\rbr{c,W^{H}-\rho}|{\cal F}_{v}}\right|\nonumber \\
 & =\left|\frac{1}{c}\int_{-c/2}^{c/2}E\sbr{K_{s,t}\rbr{c,W^{H}-\rho}-K_{s,t}\rbr{c,W^{H}-\rho-\rho'}|{\cal F}_{v}}\dd{\rho'}\right|\nonumber \\
 & \le e^{b_{H}L^{2H}c^{2}}E\sbr{\sup_{\sigma}K_{s,t}\rbr{c,W^{H}-\sigma}^{2}|{\cal F}_{v}}^{1/2}b_{H}^{1/2}L^{H}c.\label{eq:Lemma317}
\end{align}
For $c\le L^{-2H}/2$ one has $L^{2H}c^{2}\le L^{-2H}/4<1/4$ hence
$e^{b_{H}L^{2H}c^{2}}\le e^{b_{H}/4}$. Denoting $g_{H}=e^{b_{H}/4}b_{H}^{1/2}$,
for $c\le L^{-2H}/2$ and any $\tilde{x}\ge0$ we have the estimate
\begin{align}
 & \P\rbr{c^{1/H}\left|E\sbr{K_{s,t}\rbr{c,W^{H}-\rho}-\bar{K}_{s,t}\rbr{c,W^{H}}|{\cal F}_{v}}\right|\ge\tilde{x}}\nonumber \\
 & \le\P\rbr{e^{b_{H}L^{2H}c^{2}}E\sbr{\sup_{\sigma}K_{s,t}\rbr{c,W^{H}-\sigma}^{2}|{\cal F}_{v}}^{1/2}b_{H}^{1/2}L^{H}c^{1+1/H}\ge\tilde{x}}\nonumber \\
 & \le\P\rbr{e^{b_{H}/4}E\sbr{\sup_{\sigma}K_{s,t}\rbr{c,W^{H}-\sigma}^{2}|{\cal F}_{v}}^{1/2}b_{H}^{1/2}L^{H}c^{1+1/H}\ge\tilde{x}} \nonumber \\
 & =\P\rbr{E\sbr{\sup_{\sigma}K_{s,t}\rbr{c,W^{H}-\sigma}^{2}|{\cal F}_{v}}\ge\rbr{\frac{\tilde{x}}{g_{H}L^{H}c^{1+1/H}}}^{2}},\label{eq:lastlast}
\end{align}
Now, using (\ref{eq:Th12}), equality $t-s=L^{-2}$ and Lemma A (see Appendix), for
$c\le L^{-2H}/2$, $v\ge0$ and $\tilde{x}$ such that
\[
\frac{\tilde{x}}{g_{H}L^{H}c^{1+1/H}}=\bar{C}_{H}\rbr{t-s}c^{-1/H}v=\bar{C}_{H}L^{-2}c^{-1/H}v
\]
 we have 
\begin{align*}
 & \P\rbr{E\sbr{\sup_{\sigma}K_{s,t}\rbr{c,W^{H}-\sigma}^{2}|{\cal F}_{v}}\ge\rbr{\frac{\tilde{x}}{g_{H}L^{H}c^{1+1/H}}}^{2}}\\
 & =\P\rbr{E\sbr{\sup_{\sigma}K_{s,t}\rbr{c,W^{H}-\sigma}^{2}|{\cal F}_{v}}\ge\rbr{\bar{C}_{H}L^{-2}c^{-1/H}v}^{2}}\\
 & \le\tilde{A}_{H}\exp\rbr{-\frac{1}{2}\bar{B}_{H}v^{1+\min(2H,1)}}\\
 & =\tilde{A}_{H}\exp\rbr{-\frac{1}{2}\frac{\bar{B}_{H}}{\rbr{g_{H}\bar{C}_{H}}^{1+\min(2H,1)}}\rbr{\frac{L^{2-H}}{c}\tilde{x}}^{\min(1+2H,2)}}.
\end{align*}
Finally, by the just obtained estimate, defining 
\[
\bar{G}_{H}:=\frac{1}{2}\frac{\bar{B}_{H}}{\rbr{g_{H}\bar{C}_{H}}^{1+\min(2H,1)}}
\]
 and taking $\tilde{\varepsilon}=\tilde{x}\cdot L^{2}$, for $c\le L^{-2H}/2$,
we get 
\begin{align*}
 & \P\rbr{\left|\sum_{m=1}^{L}\sum_{n=2}^{L}c^{1/H}E\sbr{K_{s(l,m;L),t(l,m;L)}\rbr{c,W^{H}-\rho}-\bar{K}_{s(l,m;L),t(l,m;L)}\rbr{c,W^{H}}|{\cal F}_{t(l-1,m;L)}}\right|\ge\tilde{\varepsilon}}\\
 & \le\sum_{m=1}^{L}\sum_{l=2}^{L}\P\rbr{c^{1/H}\left|E\sbr{K_{s(l,m;L),t(l,m;L)}\rbr{c,W^{H}-\rho}-\bar{K}_{s(l,m;L),t(l,m;L)}\rbr{c,W^{H}}|{\cal F}_{t(l-1,m;L)}}\right|\ge\tilde{\varepsilon}/L^{2}}\\
 & \le\tilde{A}_{H}L^{2}\exp\rbr{-\bar{G}_{H}\rbr{\frac{\tilde{\varepsilon}}{L^{H}c}}^{1+\min(2H,1)}}.
\end{align*}
In particular, taking $L=2,3,\ldots$, $\tilde{\varepsilon}=L^{-1/4}$,
$c\le\min\cbr{L^{-H-(1/4)-1/(6+6\min(2H,1))},L^{-2H}/2}$ we get 
\begin{align}
 & \P\rbr{\left|\sum_{m=1}^{L}\sum_{n=2}^{L}c^{1/H}E\sbr{K_{s(l,m;L),t(l,m;L)}\rbr{c,W^{H}-\rho}-\bar{K}_{s(l,m;L),t(l,m;L)}\rbr{c,W^{H}}|{\cal F}_{t(l-1,m;L)}}\right|\ge L^{-1/4}}\nonumber \\
 & \le\tilde{A}_{H}L^{2}\exp\rbr{-\bar{G}_{H}L^{(1+\min(2H,1))/(6+6\min(2H,1))}}=\tilde{A}_{H}L^{2}\exp\rbr{-\bar{G}_{H}L^{1/6}}.\label{eq:estim_cond_exp_tails}
\end{align}

\subsection{Estimates of tails of $\left|K_{0,1}\protect\rbr{c,W^{H}-\rho}-\protect\E\bar{K}_{0,1}\protect\rbr{c,W^{H}}\right|$
and a.s. convergence of $c^{1/H}K_{0,1}\protect\rbr{c,W^{H}-\rho}$}

Now we will summarize obtained results to estimate tails of 
\[
c^{1/H}\left|K_{0,1}\rbr{c,W^{H}-\rho}-\E\bar{K}_{0,1}\rbr{c,W^{H}}\right|.
\]
First, using (\ref{eq:supp}), (\ref{eq:subb}) and (\ref{eq:dekompozycja}),
we get 
\begin{align}
 & c^{1/H}\left|K_{0,1}\rbr{c,W^{H}-\rho}-\E\bar{K}_{0,1}\rbr{c,W^{H}}\right|\nonumber \\
 & \le c^{1/H}K_{0,1/L}\rbr{c,W^{H}-\rho}+\left|\sum_{m=1}^{L}{\cal M}_{L}^{m}\right|+c^{1/H}L^{2}\nonumber \\
 & \quad+\left|\sum_{m=1}^{L}\sum_{l=2}^{n}c^{1/H}E\sbr{K_{s(l,m;L),t(l,m;L)}\rbr{c,W^{H}-\rho}-\bar{K}_{s(l,m;L),t(l,m;L)}\rbr{c,W^{H}}|{\cal F}_{t(l-1,m;L)}}\right|\nonumber \\
 & \quad+\left|\sum_{m=1}^{L}\sum_{l=2}^{n}c^{1/H}E\sbr{\bar{K}_{s(l,m;L),t(l,m;L)}\rbr{c,W^{H}}-\E\bar{K}_{s(l,m;L),t(l,m;L)}\rbr{c,W^{H}}|{\cal F}_{t(l-1,m;L)}}\right|\nonumber \\
 & \quad+c^{1/H}\left|\sum_{m=1}^{L}\sum_{l=2}^{n}\E\bar{K}_{s(l,m;L),t(l,m;L)}\rbr{c,W^{H}}-\E\bar{K}_{0,1}\rbr{c,W^{H}}\right|.\label{eq:variation_decom}
\end{align}

The tails of $K_{0,1/L}\rbr{c,W^{H}-\rho}$ may be estimated directly
from (\ref{eq:Th12}):
\begin{align}
 & \P\rbr{c^{1/H}K_{0,1/L}\rbr{c,W^{H}-\rho}\ge L^{-1/4}}\nonumber \\
 & =\P\rbr{K_{0,1/L}\rbr{c,W^{H}-\rho}\ge\bar{C}_{H}(1/L)c^{-1/H}\rbr{L^{3/4}/\bar{C}_{H}}}\nonumber \\
 & \le\bar{A}_{H}\exp\rbr{-\bar{B}_{H}\rbr{L^{3/4}/\bar{C}_{H}}^{1+\min(2H,1)}}.\label{eq:tails1}
\end{align}

The tails of 
\[
c^{1/H}E\sbr{\bar{K}_{s(l,m;L),t(l,m;L)}\rbr{c,W^{H}}-\E\bar{K}_{s(l,m;L),t(l,m;L)}\rbr{c,W^{H}}|{\cal F}_{t(l-1,m;L)}}
\]
may be estimated directly from the (\ref{eq:Th12-1}) and Lemma A.
By (\ref{eq:Th12-1}), using the equality
\[
\bar{K}_{s(l,m;L),t(l,m;L)}\rbr{c,W^{H}}=\frac{1}{c}\TTV{W^{H}}{\sbr{s(l,m;L),t(l,m;L)}}c
\]
and taking into account that $t(l,m;L)-s(l,m;L)=L^{-2}$, for $c\le L^{-2H}$,
$v\ge1$, we have
\begin{align}
 & \P\rbr{c^{1/H}\left|\bar{K}_{s(l,m;L),t(l,m;L)}\rbr{c,W^{H}}-\E\bar{K}_{s(l,m;L),t(l,m;L)}\rbr{c,W^{H}}\right|\ge L^{-2}v}\nonumber \\
 & =\P\rbr{c^{1/H-1}\left|\TTV{W^{H}}{\sbr{s(l,m;L),t(l,m;L)}}c-\E\TTV{W^{H}}{\sbr{s(l,m;L),t(l,m;L)}}c\right|\ge L^{-2}v}\nonumber \\
 & =\P\rbr{\left|\TTV{W^{H}}{\sbr{s(l,m;L),t(l,m;L)}}c-\E\TTV{W^{H}}{\sbr{s(l,m;L),t(l,m;L)}}c\right|\ge L^{-2}c^{1-1/H}v}\nonumber \\
 & \le\bar{A}_{H}\exp\rbr{-\bar{B}_{H}\rbr{L^{-2}c^{-1/H}}^{2-\max(2H,1)}v^{1+\min(2H,1)}},\label{eq:ittttt}
\end{align}
which, by Lemma A yields for any $v\ge0$
\begin{align}
 & \P\rbr{E\sbr{c^{1/H}\left|\bar{K}_{s(l,m;L),t(l,m;L)}\rbr{c,W^{H}}-\E\bar{K}_{s(l,m;L),t(l,m;L)}\rbr{c,W^{H}}\right||{\cal F}_{t(l-1,m;L)}}\ge L^{-2}v}\nonumber \\
 & \le\tilde{A}_{H}\exp\rbr{-\rbr{\bar{B}_{H}/2}\rbr{L^{-2}c^{-1/H}}^{2-\max(2H,1)}v^{1+\min(2H,1)}}.\label{eq:itt}
\end{align}
Let us notice that for sufficiently small $c$, one has
\begin{equation}
\rbr{L^{-2}c^{-1/H}}^{2-\max(2H,1)}\ge L^{1/6}\rbr{L^{1/4}}^{1+\min(2H,1)}.\label{eq:small_c}
\end{equation}
Easy algebraic manipulations give that (\ref{eq:small_c}) holds for
$c\le L^{-\beta}$ where
\[
\beta=H\rbr{2+\frac{\frac{1}{6}+\frac{1}{4}+\frac{\min(2H,1)}{4}}{2-\max(2H,1)}}.
\]
(\ref{eq:itt}) applied for $v=L^{-1/4}$ together with (\ref{eq:small_c})
give
\begin{align*}
 & \P\rbr{E\sbr{c^{1/H}\left|\bar{K}_{s(l,m;L),t(l,m;L)}\rbr{c,W^{H}}-\E\bar{K}_{s(l,m;L),t(l,m;L)}\rbr{c,W^{H}}\right||{\cal F}_{t(l-1,m;L)}}\ge L^{-1/4}/L^{2}}\\
 & \le\tilde{A}_{H}\exp\rbr{-\rbr{\bar{B}_{H}/2}\rbr{L^{-2}c^{-1/H}}^{2-\max(2H,1)}\rbr{L^{-1/4}}^{1+\min(2H,1)}}\\
 & \le\tilde{A}_{H}\exp\rbr{-\rbr{\bar{B}_{H}/2}L^{1/6}\rbr{L^{1/4}}^{1+\min(2H,1)}\rbr{L^{-1/4}}^{1+\min(2H,1)}}\\
 & =\tilde{A}_{H}\exp\rbr{-\rbr{\bar{B}_{H}/2}L^{1/6}}
\end{align*}
 and as a result we get 
\begin{align}
 & \P\rbr{\left|\sum_{m=1}^{L}\sum_{l=2}^{n}E\sbr{c^{1/H}\bar{K}_{s(l,m;L),t(l,m;L)}\rbr{c,W^{H}}-\E\bar{K}_{s(l,m;L),t(l,m;L)}\rbr{c,W^{H}}|{\cal F}_{t(l-1,m;L)}}\right|\ge L^{-1/4}}\nonumber \\
 & \le\sum_{m=1}^{L}\sum_{l=2}^{n}\P\rbr{\left|E\sbr{c^{1/H}\bar{K}_{s(l,m;L),t(l,m;L)}\rbr{c,W^{H}}-\E\bar{K}_{s(l,m;L),t(l,m;L)}\rbr{c,W^{H}}|{\cal F}_{t(l-1,m;L)}}\right|\ge L^{-1/4}/L^{2}}\nonumber \\
 & \le\sum_{m=1}^{L}\sum_{l=2}^{n}\tilde{A}_{H}\exp\rbr{-\rbr{\bar{B}_{H}/2}L^{1/6}}\le\tilde{A}_{H}L^{2}\exp\rbr{-\rbr{\bar{B}_{H}/2}L^{1/6}}.\label{eq:tails2}
\end{align}

Finally, let us notice that for $c\le L^{-H}$
\begin{align}
 & c^{1/H}\left|\sum_{m=1}^{L}\sum_{l=2}^{L}\E\bar{K}_{s(l,m;L),t(l,m;L)}\rbr{c,W^{H}}-\E\bar{K}_{0,1}\rbr{c,W^{H}}\right|\nonumber \\
 & \le c^{1/H}\rbr{\E\bar{K}_{0,1/L}\rbr{c,W^{H}}+L^{2}} \le c^{1/H}\rbr{2C_Hc^{-1/H}L^{-1}+L^{2}}\nonumber \\
 & =2C_HL^{-1}+c^{1/H}L^{2},\label{eq:laaast}
\end{align}
where we used (\ref{eq:sub_2}) and (\ref{eutv}).
Summarizing the estimates of tails of all random quantities appearing
in (\ref{eq:variation_decom}), that is (\ref{eq:tails1}), (\ref{eq:estim_mart_tails}),
(\ref{eq:estim_cond_exp_tails}) and (\ref{eq:tails2}) we get the
following concentration result for $c^{1/H}K_{0,1}\rbr{c,W^{H}-\rho}$.

\begin{fact}\label{fact_Kst_concentration} There exist universal positive
constants $Q_{H}$ and $R_{H}$, depending on $H$ only, such that
for $L=2,3,\ldots$, and 
\[
c\le\min\cbr{L^{-2H}/2,L^{-\alpha},L^{-\beta}},
\]
where $\alpha=H+(1/4)+1/(6+6\min(2H,1))$, $\beta=H\rbr{2+\frac{\frac{1}{6}+\frac{1}{4}+\frac{\min(2H,1)}{4}}{2-\max(2H,1)}}$,
one has 
\begin{align}
 & \P\rbr{c^{1/H}\left|K_{0,1}\rbr{c,W^{H}-\rho}-\E\bar{K}_{0,1}\rbr{c,W^{H}}\right|\ge4L^{-1/4}+2C_HL^{-1}+2c^{1/H}L^{2}}\nonumber \\
 & \le R_{H}L^{2}\exp\rbr{-Q_{H}L^{1/6}}.\label{eq:Fact_Kst_estim}
\end{align}
\end{fact}
An easy consequence of Fact \ref{fact_Kst_concentration}
is the following corollary (cf. \cite[Theorem 1.1]{Toyomu:2023}).

\begin{coro}\label{cor_Kst_convergence} Let $\rbr{c_{n}}$
be a sequence such that $c_{n}=O\rbr{n^{-\eta}}$ for some $\eta>0$.
Then for any real number $\rho$, 
\[
\lim_{n\ra+\ns}c_{n}^{1/H}K_{0,1}\rbr{c,W^{H}-\rho}=\mathfrak{c}_{H}\text{ a.s.},
\]
where $\mathfrak{c}_{H}$ is defined in Corollary \ref{cor_TV_limit}.
As a result, by scaling, for any $s,t$ such that
$0\le s<t$,
\[
\lim_{n\ra+\ns}c_{n}^{1/H}K_{s,t}\rbr{c,W^{H}-\rho}=\mathfrak{c}_{H}(t-s)\text{ a.s.}.
\]
\end{coro}
{\bf Proof.} By the assumption, there exists some positive
constant $\kappa>0$ such that for $n\in\N$, $c_{n}\le\kappa n^{-\eta}$.
For $n\in\N$ we define
\[
L_{n}:=\left\lfloor \rbr{2c_{n}}^{-1/\max\cbr{3H,\alpha,\beta}}\right\rfloor \ge\left\lfloor \rbr{2\kappa}^{-1/\max\cbr{3H,\alpha,\beta}}n^{\eta/\max\cbr{3H,\alpha,\beta}}\right\rfloor ,
\]
where $\alpha$ and $\beta$ are the same as in the statement of Fact
\ref{fact_Kst_concentration}. For sufficiently large $n$ we have
$L_{n}\ge2$ and
\begin{align*}
L_{n}^{\max\cbr{3H,\alpha,\beta}} & \le\rbr{2c_{n}}^{-1}\text{ hence }L_{n}^{-\max\cbr{3H,\alpha,\beta}}\ge2c_{n}.
\end{align*}
Thus, for sufficiently large $n$, the estimate (\ref{eq:Fact_Kst_estim})
holds for $c=c_{n}$ with $L=L_{n}$. On the other hand, 
\[
L_{n}\le\rbr{2c_{n}}^{-1/\max\cbr{3H,\alpha,\beta}},
\]
and if $2c_{n}\le1$ then
\[
2c_{n}^{1/H}L_{n}^{2}\le2^{1/H}c_{n}^{1/H}\rbr{2c_{n}}^{-2/\max\cbr{3H,\alpha,\beta}}\le\rbr{2c_{n}}^{1/(3H)}.
\]
Thus, for any $\varepsilon>0$, for sufficiently large $n\in\N$,
$2c_{n}^{1/H}L_{n}^{2}\le\varepsilon$. Now, since the sum $\sum_{n\in\N}R_{H}L_{n}^{2}\exp\rbr{-Q_{H}L_{n}^{1/6}}$
is finite, by the Borel-Cantelli lemma the difference $c_{n}^{1/H}\left|K_{0,1}\rbr{c_{n},W^{H}-\rho}-\E\bar{K}_{0,1}\rbr{c_{n},W^{H}}\right|$
tends a.s. to $0$. By definition (\ref{eq:Banach_ind_1}) of $K_{s,t}(c,f)$
and (\ref{eq:ETV_limit}), $c_{n}^{1/H}\E\bar{K}_{0,1}\rbr{c_{n},W^{H}}$
tends to $\mathfrak{c}_{H}$ a.s., thus $c_{n}^{1/H}K_{0,1}\rbr{c_{n},W^{H}-\rho}$
tends to $\mathfrak{c}_{H}$ a.s. as well.

To prove the a.s. convergence of $c_{n}^{1/H}K_{s,t}\rbr{c,W^{H}-\rho}$
let us first consider $c_{n}^{1/H}K_{0,t}\rbr{c,W^{H}-\rho}$. By
scaling, 
\[
c_{n}^{1/H}K_{0,t}\rbr{c_{n},W^{H}-\rho}=^{d}c_{n}^{1/H}K_{0,t}\rbr{c_{n},t^{H}W_{t^{-1}\cdot}^{H}-\rho}=c_{n}^{1/H}K_{0,t}\rbr{t^{-H}c_{n},W_{t^{-1}\cdot}^{H}-t^{-H}\rho}
\]
($=^{d}$ denotes equality in distribution) and
\begin{align*}
\lim_{n\ra+\ns}c_{n}^{1/H}K_{0,t}\rbr{c_{n},W^{H}-\rho} & =^{d}\lim_{n\ra+\ns}c_{n}^{1/H}K_{0,t}\rbr{t^{-H}c_{n},W_{t^{-1}\cdot}^{H}-t^{-H}\rho}\\
 & =\lim_{n\ra+\ns}t\rbr{t^{-H}c_{n}}^{1/H}K_{0,1}\rbr{t^{-H}c_{n},W^{H}-t^{-H}\rho}\\
 & =\mathfrak{c}_{H}t\text{ a.s.}
\end{align*}
Finally, the convergence of $c_{n}^{1/H}K_{s,t}\rbr{c_{n},W^{H}-\rho}$
follows from (\ref{eq:sup}) and (\ref{eq:sub}). 
\hfill$\blacksquare$

\begin{rema} \label{rem_Kst_not_monotonic} It is not clear if the a.s.
convergence of $c_{n}^{1/H}K_{0,1}\rbr{c_{n},W^{H}-\rho}$ holds for
any sequence $\rbr{c_{n}}$ such that $c_{n}\ra0_{+}$. For given
continuous $f:[0,+\ns)\ra\R$, the quatity $K_{0,1}\rbr{c,f}$ may
behave in a much less stable way than $\bar{K}_{0,1}\rbr{c,f}$. For
example, the mapping $(0,1)\ni c\mapsto K_{0,1}\rbr{c,f}$ may be
non-monotonic. 
\end{rema}
An easy consequence of Corollary \ref{cor_Kst_convergence} are the
convergences of $$c_{n}^{1/H}KU_{0,1}\rbr{c_{n},W^{H}-\rho}  \text{ and } c_{n}^{1/H}KD_{0,1}\rbr{c_{n},W^{H}-\rho}.$$

\begin{coro} Let $\rbr{c_{n}}$ be a sequence such
that $c_{n}=O\rbr{n^{-\eta}}$ for some $\eta>0$. Then for any real
number $\rho$ and $s$, $t$ such that $0\le s<t$,
\[
\lim_{n\ra+\ns}c_{n}^{1/H}KU_{s,t}\rbr{c,W^{H}-\rho}=\lim_{n\ra+\ns}c_{n}^{1/H}KU_{s,t}\rbr{c,W^{H}-\rho}=\frac{1}{2}\mathfrak{c}_{H}(t-s)\text{ a.s.},
\]
where $\mathfrak{c}_{H}$ is defined in Corollary \ref{cor_TV_limit}. 
\end{coro}
{\bf Proof.} We have 
\begin{equation}
\left|KU_{s,t}\rbr{c,W^{H}-\rho}-KD_{s,t}\rbr{c,W^{H}-\rho}\right|\le2\max_{u\in[s,t]}\left|W_{u}^{H}\right|/c+3\label{eq:ku-kd_estim}
\end{equation}
This follows from the following observations: for $0\le s<t$,
\[
\left|U_{s,t}\rbr{c,W^{H}-p\cdot c-\rho}-D_{s,t}\rbr{c,W^{H}-p\cdot c-\rho}\right|\le1
\]
and for $p\in\Z$ such that $pc>\max_{u\in[s,t]}\left|W_{u}^{H}\right|+c$
we have $U_{s,t}\rbr{c,W^{H}-p\cdot c-\rho}=D_{s,t}\rbr{c,W^{H}-p\cdot c-\rho}=0$.
As a result we have 
\begin{align*}
 & \left|KU_{s,t}\rbr{c,W^{H}-\rho}-KD_{s,t}\rbr{c,W^{H}-\rho}\right|\\
 & =\left|\sum_{p\in\Z}\cbr{U_{s,t}\rbr{c,W^{H}-p\cdot c-\rho}-D_{s,t}\rbr{c,W^{H}-p\cdot c-\rho}}\right|\\
 & \le\sum_{p\in\Z}\left|U_{s,t}\rbr{c,W^{H}-p\cdot c-\rho}-D_{s,t}\rbr{c,W^{H}-p\cdot c-\rho}\right|\\
 & \le\sum_{p\in\Z,p\le\max_{u\in[s,t]}\left|W_{u}^{H}\right|/c+1}\left|U_{s,t}\rbr{c,W^{H}-p\cdot c-\rho}-D_{s,t}\rbr{c,W^{H}-p\cdot c-\rho}\right|\\
 & \le2\max_{u\in[s,t]}\left|W_{u}^{H}\right|/c+3.
\end{align*}
Now, since $K_{s,t}\rbr{c,W^{H}-\rho}=KU_{s,t}\rbr{c,W^{H}-\rho}+KD_{s,t}\rbr{c,W^{H}-\rho}$,
we have the estimate
\begin{align*}
 & 2c^{1/H}KU_{s,t}\rbr{c,W^{H}-\rho}\\
 & =c^{1/H}K{}_{s,t}\rbr{c,W^{H}-\rho}+c^{1/H}KU_{s,t}\rbr{c,W^{H}-\rho}-KD_{s,t}\rbr{c,W^{H}-\rho}\\
 & \le c^{1/H}K{}_{s,t}\rbr{c,W^{H}-\rho}+c^{1/H}\left|KU_{s,t}\rbr{c,W^{H}-\rho}-KD_{s,t}\rbr{c,W^{H}-\rho}\right|\\
 & \le c^{1/H}K{}_{s,t}\rbr{c,W^{H}-\rho}+2c^{1/H-1}\max_{u\in[s,t]}\left|W_{u}^{H}\right|+3c^{1/H}
\end{align*}
and 
\begin{align*}
 & 2c^{1/H}KU_{s,t}\rbr{c,W^{H}-\rho}\\
 & \le c^{1/H}K{}_{s,t}\rbr{c,W^{H}-\rho}-2c^{1/H-1}\max_{u\in[s,t]}\left|W_{u}^{H}\right|-3c^{1/H}.
\end{align*}
Since 
\[
\lim_{c\ra0_{+}}c^{1/H-1}\max_{u\in[s,t]}\left|W_{u}^{H}\right|=0\text{ a.s.}
\]
we finally have 
\[
\lim_{n\ra+\ns}2c^{1/H}KU_{s,t}\rbr{c,W^{H}-\rho}=\lim_{n\ra+\ns}c_{n}^{1/H}K_{s,t}\rbr{c,W^{H}-\rho}=\mathfrak{c}_{H}(t-s)\text{ a.s}.
\]
The convergence of $c_{n}^{1/H}KD_{s,t}\rbr{c,W^{H}-\rho}$ is proven
analogously.
\hfill $\blacksquare$

\section{Local times} \label{sect:loc_times}

We will now proceed -- using similar techniques as in Sect. \ref{sec:Concentration}
-- to prove that normalized numbers of upcrossings almost
surely converge in $L^{1}(\R)$ to the local time of fBm with any
Hurst parameter $H\in(0,1)$. 

To prove the claimed convergence we need to prove for any $g\in L^{\infty}(\R)$
the convergence 
\[
\int_{\R}c^{1/H-1}U_{0,1}\rbr{c,W^{H}-a}g(a)\dd a\ra\frac{\mathfrak{c}_{H}}{2}\int_{\R}{\cal L}_{1}^{a}\rbr{W^{H}}g(a)\dd a\ \text{a.s.,}\ \text{as}\ c\ra0_{+},
\]
where ${\cal L}_{1}\rbr{W^{H}}$ denotes the local time of $W^{H}$
at time $1$, defined as the denisty (with resp. to the Lebesgue measure)
of the occupation measure 
\[
\mu\rbr A:=\int_{0}^{1}{\bf 1}_{A}\rbr{W_{t}^{H}}\dd t,\quad A\in{\cal B}\rbr{\R},
\]
and $\mathfrak{c}_{H}$ is the same positive constant as in Corollary
\ref{cor_TV_limit}.

\emph{Outline}. Let us present the outline of our reasoning. Let $\psi_{1}:\R\ra[0,+\ns)$
be a mollifier \textendash{} a smooth ($C^{\ns}$) function such that
$\int_{-1/2}^{1/2}$ $\psi_{1}(\rho)\dd{\rho}=1$ and $\psi_{1}(\rho)=0$
whenever $\left|\rho\right|>1/2$. Next, for $r>0$ we define $\psi_{r}:\R\ra[0,+\ns)$
by
\[
\psi_{r}(\rho):=\frac{1}{r}\psi_{1}(\rho/r).
\]
For any $r>0$, $\psi_{r}$ is again a smooth ($C^{\ns}$), non-negative
function, and such that $\int_{-r/2}^{r/2}$ $\psi_{r}(\rho)\dd{\rho}=1$
and $\psi_{r}(\rho)=0$ whenever $\left|\rho\right|>r/2$. Next, let
us consider the differences
\begin{align}
 & \int_{\R}c^{1/H-1}U_{0,1}\rbr{c,W^{H}-a}g(a)\dd a-\int_{-r/2}^{r/2}\cbr{\int_{\R}\psi_{r}(\rho)c^{1/H-1}U_{0,1}\rbr{c,W^{H}-a-\rho}g(a)\dd a}\dd{\rho} \nonumber \\
 & =\int_{-r/2}^{r/2}\psi_{r}(\rho)\int_{\R}c^{1/H-1}\cbr{U_{0,1}\rbr{c,W^{H}-a}-U_{0,1}\rbr{c,W^{H}-a-\rho}}g(a)\dd a\dd{\rho}\label{eq:dife1}
\end{align}
and
\begin{equation}
\int_{\R}{\cal L}_{1}^{a}\rbr{W^{H}}g(a)\dd a-\int_{-r/2}^{r/2}\psi(\rho)\left\{ \int_{\R}{\cal L}_{1}^{a-\rho}\rbr{W^{H}}g(a)\dd a\right\} \dd{\rho}.\label{eq:dife2}
\end{equation}
The a.s. convergence of the first difference (\ref{eq:dife1}) to
$0$, for some sequences of positive reals $\rbr{c_{n}}$ and $\rbr{r_{n}}$,
will be established similarly as in the previous section, while the
second difference (\ref{eq:dife2}) may be estimated from the H\"older
regularity of the local time of the fBm in the space variable as follows:
\begin{align*}
 & \left|\int_{\R}{\cal L}_{1}^{a}\rbr{W^{H}}g(a)\dd a-\int_{-r/2}^{r/2}\psi(\rho)\left\{ \int_{\R}{\cal L}_{1}^{a-\rho}\rbr{W^{H}}g(a)\dd a\right\} \dd{\rho}\right|\\
 & =\left|\int_{-r/2}^{r/2}\left\{ \int_{\R}\psi(\rho)\rbr{{\cal L}_{1}^{a}\rbr{W^{H}}-{\cal L}_{1}^{a-\rho}\rbr{W^{H}}}g(a)\dd a\right\} \dd{\rho}\right|\\
 & \le\int_{-r/2}^{r/2}\psi(\rho)\int_{\R}\left|{\cal L}_{1}^{a}\rbr{W^{H}}-{\cal L}_{1}^{a-\rho}\rbr{W^{H}}\right|\left|g(a)\right|\dd a\dd{\rho}\\
 & \le\cbr{\sup_{\left|\rho\right|\le r/2}\left|{\cal L}_{1}^{a}\rbr{W^{H}}-{\cal L}_{1}^{a-\rho}\rbr{W^{H}}\right|}\rbr{2\sup_{t\in[0,1]}W_{t}^{H}+r}\left\Vert g\right\Vert _{\infty},
\end{align*}
which converges a.s. as $r\ra0_{+}$, since 
\[
\sup_{\left|\rho\right|\le r/2}\frac{\left|{\cal L}_{1}^{a}\rbr{W^{H}}-{\cal L}_{1}^{a-\rho}\rbr{W^{H}}\right|}{\left|\rho\right|^{\gamma}}<\ns\text{ a.s.}
\]
for any $\gamma\le\min\rbr{1,(1-H)/(2H)}$, see \cite[Subsect. 3.30]{GemanHorowitz80}. 

Next, we will prove the a.s. convergence of the difference of mollified
versions, i.e.
\begin{align*}
 & \int_{-r/2}^{r/2}\cbr{\int_{\R}\psi_{r}(\rho)c^{1/H-1}U_{0,1}\rbr{c,W^{H}-a-\rho}g(a)\dd a}\dd{\rho}\\
 & -\int_{-r/2}^{r/2}\psi(\rho)\left\{ \int_{\R}{\cal L}_{1}^{a-\rho}\rbr{W^{H}}g(a)\dd a\right\} \dd{\rho}
\end{align*}
for the same sequences of positive reals $\rbr{c_{n}}$ and $\rbr{r_{n}}$
as those used in establishing the a.s. convergence of (\ref{eq:dife1}).

This way we will establish a.s. convergence of $\int_{\R}c^{1/H-1}U_{0,1}\rbr{c,W^{H}-a}g(a)\dd a$
to $\int_{\R}{\cal L}_{1}^{a}\rbr{W^{H}}g(a)\dd a$ along the mentioned
sequence $\rbr{c_{n}}$. Since $c_{n}$s  tend to $0$ sufficiently
slowly, we will be able, using also the monotonicity of the mapping
$c\mapsto U_{0,1}\rbr{c,f}$, to drop the assumption about $c_{n}$s
and establish the a.s. convergence $\int_{\R}c^{1/H-1}U_{0,1}\rbr{c,W^{H}-a}g(a)\dd a$
to $\int_{\R}{\cal L}_{1}^{a}\rbr{W^{H}}g(a)\dd a$ as $c\ra0_{+}$.

\subsection{Estimation of difference (\ref{eq:dife1})}

To estimate the first difference (\ref{eq:dife1}), let us choose
a positive integer $L$ greater or equal $2$. Using superadditivity
(\ref{eq:sup-1}) and subadditivity (\ref{eq:sub-1}), for $m\le L-1$
we naturally have
\begin{align*}
 & U_{s(l,m;L),t(l,m;L)}\rbr{c,W^{H}-a}-U_{s(l,m;L),t(l,m;L)}\rbr{c,W^{H}-a-\rho}\\
 & +U_{s(l,m+1;L),t(l,m+1;L)}\rbr{c,W^{H}-a}-U_{s(l,m+1;L),t(l,m+1;L)}\rbr{c,W^{H}-a-\rho}\\
 & \le U_{s(l,m;L),t(l,m+1;L)}\rbr{c,W^{H}-a}-U_{s(l,m;L),t(l,m+1;L)}\rbr{c,W^{H}-a-\rho}+1
\end{align*}
and 
\begin{align*}
 & U_{s(l,m;L),t(l,m;L)}\rbr{c,W^{H}-a}-U_{s(l,m;L),t(l,m;L)}\rbr{c,W^{H}-a-\rho}\\
 & +U_{s(l,m+1;L),t(l,m+1;L)}\rbr{c,W^{H}-a}-U_{s(l,m+1;L),t(l,m+1;L)}\rbr{c,W^{H}-a-\rho}\\
 & \ge U_{s(l,m;L),t(l,m+1;L)}\rbr{c,W^{H}-a}-U_{s(l,m;L),t(l,m+1;L)}\rbr{c,W^{H}-a-\rho}-1.
\end{align*}
These yield for any real $a$ and $\rho$ the estimates 
\begin{align*}
 & U_{0,1}\rbr{c,W^{H}-a}-U_{0,1}\rbr{c,W^{H}-a-\rho}\\
 & \le\sum_{l=1}^{L}\sum_{m=1}^{L}\cbr{U_{s(l,m;L),t(l,m;L)}\rbr{c,W^{H}-a}-U_{s(l,m;L),t(l,m;L)}\rbr{c,W^{H}-a-\rho}}+L^{2}
\end{align*}
and 
\begin{align*}
 & U_{0,1}\rbr{c,W^{H}-a}-U_{0,1}\rbr{c,W^{H}-a-\rho}\\
 & \ge\sum_{l=1}^{L}\sum_{m=1}^{L}\cbr{U_{s(l,m;L),t(l,m;L)}\rbr{c,W^{H}-a}-U_{s(l,m;L),t(l,m;L)}\rbr{c,W^{H}-a-\rho}}-L^{2}.
\end{align*}
Taking also into account that for $a$ and $\rho$ such that $\left|\rho\right|\le r/2$
and $\left|a\right|>\max_{t\in[0,1]}\left|W_{t}^{H}\right|+r/2$,
\[
U_{s(l,m;L),t(l,m;L)}\rbr{c,W^{H}-a}=0,\quad U_{s(l,m;L),t(l,m;L)}\rbr{c,W^{H}-a-\rho}=0
\]
we get 
\begin{align}
 & \int_{-r/2}^{r/2}\psi_{r}(\rho)\int_{\R}c^{1/H-1}\cbr{U_{0,1}\rbr{c,W^{H}-a}-U_{0,1}\rbr{c,W^{H}-a-\rho}}g(a)\dd a\dd{\rho}\nonumber \\
 & \le\sum_{l=1}^{L}\sum_{m=1}^{L}\int_{-r/2}^{r/2}\psi_{r}(\rho)\int_{\R}c^{1/H-1}\Delta U_{s(l,m;L),t(l,m;L)}\rbr{c,W^{H}-a,\rho}g(a)\dd a\dd{\rho}\nonumber \\
 & \quad+\int_{-r/2}^{r/2}\psi_{r}(\rho)\int_{\cbr{a\in\R:\left|a\right|\le\max_{t\in[0,1]}\left|W_{t}^{H}\right|+r/2}}c^{1/H-1}L^{2}g(a)\dd a\dd{\rho}\nonumber \\
 & \le\sum_{l=1}^{L}\sum_{m=1}^{L}\int_{-r/2}^{r/2}\psi_{r}(\rho)\int_{\R}c^{1/H-1}\Delta U_{s(l,m;L),t(l,m;L)}\rbr{c,W^{H}-a,\rho}g(a)\dd a\dd{\rho}\nonumber \\
 & \quad+\rbr{2\max_{t\in[0,1]}\left|W_{t}^{H}\right|+r}c^{1/H-1}L^{2}\left\Vert g\right\Vert _{\ns}\label{eq:upper}
\end{align}
and 
\begin{align}
 & \int_{-r/2}^{r/2}\psi_{r}(\rho)\int_{\R}c^{1/H-1}\cbr{U_{0,1}\rbr{c,W^{H}-a}-U_{0,1}\rbr{c,W^{H}-a-\rho}}g(a)\dd a\dd{\rho}\nonumber \\
 & \ge\sum_{l=1}^{L}\sum_{m=1}^{L}\int_{-r/2}^{r/2}\psi_{r}(\rho)\int_{\R}c^{1/H-1}\Delta U_{s(l,m;L),t(l,m;L)}\rbr{c,W^{H}-a,\rho}g(a)\dd a\dd{\rho}\nonumber \\
 & \quad-\rbr{2\max_{t\in[0,1]}\left|W_{t}^{H}\right|+r}c^{1/H-1}L^{2}\left\Vert g\right\Vert _{\ns},\label{eq:lower}
\end{align}
where 
\begin{align*}
 & \Delta U_{s(l,m;L),t(l,m;L)}\rbr{c,W^{H}-a,\rho}\\
 & :=c^{1/H-1}\cbr{U_{s(l,m;L),t(l,m;L)}\rbr{c,W^{H}-a}-U_{s(l,m;L),t(l,m;L)}\rbr{c,W^{H}-a-\rho}}.
\end{align*}
Similarly as in Subsect. \ref{subsec:Decomposition}, for each $m=1,2,\ldots,L$
we introduce the martingale 
\[
{\cal U}_{n}^{m}:=\sum_{l=2}^{n}\psi_{r}(\rho)\int_{-r/2}^{r/2}\int_{\R}\Delta UE_{s(l,m;L),t(l,m;L)}\rbr{c,W^{H}-a,\rho}g(a)\dd a\dd{\rho},
\]
where $n=1,2,\ldots,L$ and
\begin{align*}
 & \Delta UE_{s(l,m;L),t(l,m;L)}\rbr{c,W^{H}-a,\rho}\\
 & :=\Delta U_{s(l,m;L),t(l,m;L)}\rbr{c,W^{H}-a,\rho}-E\sbr{\Delta U_{s(l,m;L),t(l,m;L)}\rbr{c,W^{H}-a,\rho}|{\cal F}_{t(l-1,m;L)}}.
\end{align*}
Now, by (\ref{eq:upper}) and (\ref{eq:lower}) we decompose 
\begin{align}
 & \int_{-r/2}^{r/2}\psi_{r}(\rho)\int_{\R}c^{1/H-1}\cbr{U_{0,1}\rbr{c,W^{H}-a}-U_{0,1}\rbr{c,W^{H}-a-\rho}}g(a)\dd a\dd{\rho}\nonumber \\
 & =\int_{-r/2}^{r/2}\psi_{r}(\rho)\int_{\R}c^{1/H-1}\cbr{U_{0,1/L}\rbr{c,W^{H}-a}-U_{0,1/L}\rbr{c,W^{H}-a-\rho}}g(a)\dd a\dd{\rho}\nonumber \\
 & \quad+\sum_{m=1}^{L}{\cal U}_{L}^{m}+R\rbr{c,L,g,r}\nonumber \\
 & \quad+\sum_{m=1}^{L}\sum_{l=2}^{L}\int_{-r/2}^{r/2}\psi_{r}(\rho)\int_{\R}E\sbr{\Delta U_{s(l,m;L),t(l,m;L)}\rbr{c,W^{H}-a,\rho}|{\cal F}_{t(l-1,m;L)}}g(a)\dd a\dd{\rho},\label{eq:dife1_decomp}
\end{align}
where 
\begin{equation}
\left|R\rbr{c,L,g,r}\right|\le\rbr{2\max_{t\in[0,1]}\rbr{W_{t}^{H}}+r}\left\Vert g\right\Vert _{\ns}c^{1/H-1}L^{2}.\label{eq:remainder}
\end{equation}

\subsubsection{Tails of the sum of martingales }

Now let us notice that for any $\sigma\in\R$ and $0\le s<t<+\ns$
\begin{align*}
 & \left|\int_{\R}c^{1/H-1}U_{s,t}\rbr{c,W^{H}-a-\sigma}g(a)\dd a\right|\\
 & \le c^{1/H-1}\left\Vert g\right\Vert _{\ns}\int_{\R}U_{s,t}\rbr{c,W^{H}-a-\sigma}\dd a\\
 & =c^{1/H-1}\left\Vert g\right\Vert _{\ns}\sum_{p\in\Z}\int_{pc-c/2-\sigma}^{pc+c/2-\sigma}U_{s,t}\rbr{c,W^{H}-a-\sigma}\dd a\\
 & =\cbr{a+\sigma=pc+\alpha}=c^{1/H-1}\left\Vert g\right\Vert _{\ns}\sum_{p\in\Z}\int_{-c/2}^{c/2}U_{s,t}\rbr{c,W^{H}-pc-\alpha}\dd{\alpha}\\
 & =\left\Vert g\right\Vert _{\ns}c^{1/H-1}\int_{-c/2}^{c/2}\sum_{p\in\Z}U_{s,t}\rbr{c,W^{H}-pc-\alpha}\dd{\alpha}\\
 & \le\left\Vert g\right\Vert _{\ns}c^{1/H-1}\int_{-c/2}^{c/2}K_{s,t}\rbr{c,W^{H}-\alpha}\dd{\alpha}=\left\Vert g\right\Vert _{\ns}c^{1/H-1}\TTV{W^{H}}{[s,t]}c.
\end{align*}
As a result (taking $\sigma=0$ and $\sigma=\rho$), we also have
\begin{align}
 & \left|\int_{-r/2}^{r/2}\psi_{r}(\rho)\int_{\R}c^{1/H-1}\cbr{U_{s,t}\rbr{c,W^{H}-a}-U_{s,t}\rbr{c,W^{H}-a-\rho}}g(a)\dd a\dd{\rho}\right|\nonumber \\
 & \le\int_{-r/2}^{r/2}\psi_{r}(\rho)2\left\Vert g\right\Vert _{\ns}c^{1/H-1}\TTV{W^{H}}{[s,t]}c\dd{\rho}\nonumber \\
 & =2\left\Vert g\right\Vert _{\ns}c^{1/H-1}\TTV{W^{H}}{[s,t]}c.\label{eq:oszac}
\end{align}
Now, by (\ref{eq:oszac}), the equality $t(l,m;L)-s(l,m;L)=L^{-2}$
and (\ref{eq:Th12-1-2}) we see that for any $c\le L^{-2H}$, $v\ge1$
\begin{align*}
 & \P\rbr{\left|\int_{-r/2}^{r/2}\psi_{r}(\rho)\int_{\R}\Delta U_{s(l,m;L),t(l,m;L)}\rbr{c,W^{H}-a,\rho}g(a)\dd a\dd{\rho}\right|\ge2\bar{D}_{H}\left\Vert g\right\Vert _{\ns}L^{-2}v}\\
 & \le\P\rbr{2\left\Vert g\right\Vert _{\ns}c^{1/H-1}\TTV{W^{H}}{[s(l,m;L),t(l,m;L)]}c\ge2\bar{D}_{H}\left\Vert g\right\Vert _{\ns}L^{-2}v}\\
 & =\P\rbr{\TTV{W^{H}}{[s(l,m;L),t(l,m;L)]}c\ge\bar{D}_{H}L^{-2}c^{1-1/H}}\\
 & \le\bar{A}_{H}\exp\rbr{-\bar{B}_{H}v^{\min(1+2H,2)}}.
\end{align*}
As a result, by Lemma B, for any $x\ge0$, $c\le L^{-2H}$ and $v>0$
we have the estimate 
\[
\P\rbr{\left|{\cal U}_{L}^{m}\right|\ge x}\le2\exp\rbr{-\frac{L^{3}x^{2}}{40\rbr{\bar{D}_{H}\left\Vert g\right\Vert _{\ns}v}^{2}}}+6L\tilde{A}_{H}\exp\rbr{-(\bar{B}_{H}/2)v^{\min(1+2H,2)}}.
\]
Similarly as in Subsect. \ref{subsec:Martingale}, taking $\varepsilon=x\cdot L=L^{-1/4}$,
$v=L^{1/6}$ we get 
\begin{align}
 & \P\rbr{\left|\sum_{m=1}^{L}{\cal U}_{L}^{m}\right|\ge L^{-1/4}}\le\sum_{m=1}^{L}\P\rbr{\left|U_{L}^{m}\right|\ge\varepsilon/L}\nonumber \\
 & \le2L\exp\rbr{-\frac{L\varepsilon^{2}}{40\rbr{\bar{D}_{H}\left\Vert g\right\Vert _{\ns}v}^{2}}}+6L^{2}\tilde{A}_{H}\exp\rbr{-\rbr{\bar{B}_{H}/2}v^{\min(1+2H,2)}}\nonumber \\
 & \le\rbr{1+6\tilde{A}_{H}}L^{2}\exp\rbr{-\min\rbr{\frac{1}{40\bar{D}_{H}^{2}\left\Vert g\right\Vert _{\ns}^{2}},\frac{\bar{B}_{H}}{2}}L^{1/6}}.\label{eq:mart_local}
\end{align}

\subsubsection{Tails of the sum of conditional expectations }

Now we will deal with the quantites

\[
\int_{-r/2}^{r/2}\psi_{r}(\rho)\int_{\R}E\sbr{\Delta U_{s(l,m;L),t(l,m;L)}\rbr{c,W^{H}-a,\rho}|{\cal F}_{t(l-1,m;L)}}g(a)\dd a\dd{\rho}.
\]
We will proceed in a similar way as in Sect. 2.5 with $E\sbr{K_{s(l,m;L),t(l,m;L)}\rbr{c,W^{H}-a,\rho}|{\cal F}_{t(l-1,m;L)}}$.
Let us fix $l=2,3,\ldots,L$, $m=1,2,\ldots,L$. Let us recall (\ref{eq:stlm_def})
and denote
\[
v=t(l-1,m;L)=\frac{l-2}{L}+\frac{m}{L^{2}},\quad s=s(l,m;L)=\frac{l-1}{L}+\frac{m-1}{L^{2}},\quad t=t(l,m;L)=\frac{l-1}{L}+\frac{m}{L^{2}}.
\]
We have 
\[
(s-v)^{-1}(t-v)^{1-H}=\rbr{\frac{1}{L}-\frac{1}{L^{2}}}^{-1}\rbr{\frac{1}{L}}^{1-H}\le2L^{H}.
\]
Now, using a very similar reasoning as in the proof {[}AoP, Lemma
2.17{]} for just defined $v$, $s$ and $t$ we obtain (with probability
$1$) 
\begin{align*}
 & \left|E\sbr{\Delta U_{s,t}\rbr{c,W^{H}-a,\rho}|{\cal F}_{v}}\right|\\
 & \le e^{b_{H}4L^{2H}|\rho|^{2}}c^{1/H-1}E\sbr{U_{s,t}\rbr{c,W^{H}-a}^{2}|{\cal F}_{v}}^{1/2}b_{H}^{1/2}2L^{H}|\rho|
\end{align*}
where $b_{H}$ is the same as in Sect. (\ref{sec:Concentration}).
Hence (also with probability $1$)
\begin{align*}
 & \left|\int_{-r/2}^{r/2}\psi_{r}(\rho)\int_{\R}E\sbr{\Delta U_{s,t}\rbr{c,W^{H}-a,\rho}|{\cal F}_{v}}g(a)\dd a\dd{\rho}\right|\\
 & =\left|\int_{-r/2}^{r/2}\psi_{r}(\rho)E\sbr{\int_{\R}\Delta U_{s,t}\rbr{c,W^{H}-a,\rho}g(a)\dd a|{\cal F}_{v}}\dd{\rho}\right|\\
 & \le\int_{-r/2}^{r/2}\psi_{r}(\rho)e^{b_{H}4L^{2H}|\rho|^{2}}E\sbr{\rbr{\int_{\R}c^{1/H-1}U_{s,t}\rbr{c,W^{H}-a}g(a)\dd a}^{2}|{\cal F}_{v}}^{1/2}\rbr{b_{H}^{1/2}2L^{H}|\rho|}\dd{\rho}\\
 & \le\left\Vert g\right\Vert _{\ns}e^{b_{H}L^{2H}r^{2}}\rbr{b_{H}^{1/2}L^{H}r}E\sbr{\rbr{c^{1/H-1}\int_{\R}U_{s,t}\rbr{c,W^{H}-a}\dd a}^{2}|{\cal F}_{v}}^{1/2}\\
 & \le\left\Vert g\right\Vert _{\ns}e^{b_{H}L^{2H}r^{2}}\rbr{b_{H}^{1/2}L^{H}r}E\sbr{\rbr{c^{1/H-1}\TTV{W^{H}}{[s,t]}c}^{2}|{\cal F}_{v}}^{1/2}.
\end{align*}
For $r\le L^{-H}/2$ one has $L^{2H}r^{2}\le1/4$ hence $e^{b_{H}L^{2H}c^{2}}\le e^{b_{H}/4}$
and for any $\tilde{x}\ge0$ we have the estimate (with the same $g_{H}=e^{b_{H}/4}b_{H}^{1/2}$
as before):
\begin{align*}
 & \P\rbr{\left|\int_{-r/2}^{r/2}\psi_{r}(\rho)\int_{\R}E\sbr{\Delta U_{s,t}\rbr{c,W^{H}-a,\rho}|{\cal F}_{v}}g(a)\dd a\dd{\rho}\right|\ge\tilde{x}}\\
 & \le\P\rbr{E\sbr{\rbr{c^{1/H-1}\TTV{W^{H}}{[s,t]}c}^{2}|{\cal F}_{v}}^{1/2}\ge\frac{\tilde{x}}{\left\Vert g\right\Vert _{\ns}e^{b_{H}L^{2H}r^{2}}\rbr{b_{H}^{1/2}L^{H}r}}}\\
 & \le\P\rbr{E\sbr{\rbr{\TTV{W^{H}}{[s,t]}c}^{2}|{\cal F}_{v}}\ge\rbr{c^{1-1/H}\frac{\tilde{x}}{g_{H}\left\Vert g\right\Vert _{\ns}L^{H}r}}^{2}}.
\end{align*}
Now, using (\ref{eq:Th12-1-2}), equality $t-s=L^{-2}$ and Lemma
A, for $r>0$, $c\le L^{-2H}$, $v\ge0$ and $\tilde{x}$ such that
\[
c^{1-1/H}\frac{\tilde{x}}{g_{H}\left\Vert g\right\Vert _{\ns}L^{H}r}=\bar{D}_{H}\rbr{t-s}c^{1-1/H}v=\bar{D}_{H}L^{-2}c^{1-1/H}v
\]
 we have 
\begin{align*}
 & \P\rbr{E\sbr{\rbr{\TTV{W^{H}}{[s,t]}c}^{2}|{\cal F}_{v}}\ge\rbr{c^{1-1/H}\frac{\tilde{x}}{g_{H}\left\Vert g\right\Vert _{\ns}L^{H}r}}^{2}}\\
 & =\P\rbr{E\sbr{\rbr{\TTV{W^{H}}{[s,t]}c}^{2}|{\cal F}_{v}}\ge\rbr{\bar{D}_{H}L^{-2}c^{1-1/H}v}^{2}}\\
 & \le\tilde{A}_{H}\exp\rbr{-\frac{1}{2}\bar{B}_{H}v^{1+\min(2H,1)}}\\
 & =\tilde{A}_{H}\exp\rbr{-\frac{1}{2}\frac{\bar{B}_{H}}{\rbr{g_{H}\bar{D}_{H}}^{1+\min(2H,1)}}\rbr{\frac{L^{2-H}}{\left\Vert g\right\Vert _{\ns}r}\tilde{x}}^{\min(1+2H,2)}}.
\end{align*}
Denoting 
\[
\bar{T}_{H}:=\frac{1}{2}\frac{\bar{B}_{H}}{\rbr{g_{H}\bar{D}_{H}}^{1+\min(2H,1)}},
\]
 and taking $\tilde{\varepsilon}=\tilde{x}\cdot L^{2}$, for $r\le L^{-H}/2$,
$c\le L^{-2H}$, we get 
\begin{align*}
 & \P\rbr{\left|\sum_{m=1}^{L}\sum_{n=2}^{L}\int_{-r/2}^{r/2}\psi_{r}(\rho)\int_{\R}E\sbr{\Delta U_{s(l,m;L),t(l,m;L)}\rbr{c,W^{H}-a,\rho}|{\cal F}_{t(l-1,m;L)}}g(a)\dd a\dd{\rho}\right|\ge\tilde{\varepsilon}}\\
 & \le\sum_{m=1}^{L}\sum_{l=2}^{L}\P\rbr{c^{1/H}\left|E\sbr{K_{s(l,m;L),t(l,m;L)}\rbr{c,W^{H}-\rho}-\bar{K}_{s(l,m;L),t(l,m;L)}\rbr{c,W^{H}}|{\cal F}_{t(l-1,m;L)}}\right|\ge\tilde{\varepsilon}/L^{2}}\\
 & \le\tilde{A}_{H}L^{2}\exp\rbr{-\bar{T}_{H}\rbr{\frac{\tilde{\varepsilon}}{\left\Vert g\right\Vert _{\ns}L^{H}r}}^{1+\min(2H,1)}}.
\end{align*}
In particular, taking $\tilde{\varepsilon}=L^{-1/4}$, $r\le\min\cbr{L^{-H-(1/4)-1/(6+6\min(2H,1))},L^{-H}/2}$
we get 
\begin{align}
 & \P\rbr{\left|\sum_{m=1}^{L}\sum_{n=2}^{L}\int_{-r/2}^{r/2}\psi_{r}(\rho)\int_{\R}E\sbr{\Delta U_{s(l,m;L),t(l,m;L)}\rbr{c,W^{H}-a,\rho}|{\cal F}_{t(l-1,m;L)}}g(a)\dd a\dd{\rho}\right|\ge L^{-1/4}}\nonumber \\
 & \le\tilde{A}_{H}L^{2}\exp\rbr{-\bar{T}_{H}\left\Vert g\right\Vert _{\ns}^{-1-\min(2H,1)}L^{1/6}}.\label{eq:cond_exp_local}
\end{align}

\subsubsection{Tails of the difference (\ref{eq:dife1}) and its convergence}

In order to obtain sufficient estimates of the tails of the difference
(\ref{eq:dife1}) and to establish its a.s. convergence to $0$ as
$c\ra0_{+}$ we also need to estimate the tails of the first summand
in (\ref{eq:dife1_decomp}), i.e.
\[
\int_{-r/2}^{r/2}\psi_{r}(\rho)\int_{\R}c^{1/H-1}\cbr{U_{0,1/L}\rbr{c,W^{H}-a}-U_{0,1/L}\rbr{c,W^{H}-a-\rho}}g(a)\dd a\dd{\rho}.
\]
Using (\ref{eq:oszac}) we have 
\begin{align*}
 & \left|\int_{-r/2}^{r/2}\psi_{r}(\rho)\int_{\R}c^{1/H-1}\cbr{U_{0,1/L}\rbr{c,W^{H}-a}-U_{0,1/L}\rbr{c,W^{H}-a-\rho}}g(a)\dd a\dd{\rho}\right|\\
 & \le2\left\Vert g\right\Vert _{\ns}c^{1/H-1}\TTV{W^{H}}{[1,1/L]}c
\end{align*}
and by (\ref{eq:Th12-1-2}) and Lemma A:
\begin{align}
 & \P\rbr{\left|\int_{-r/2}^{r/2}\psi_{r}(\rho)\int_{\R}c^{1/H-1}\cbr{U_{0,1/L}\rbr{c,W^{H}-a}-U_{0,1/L}\rbr{c,W^{H}-a-\rho}}g(a)\dd a\dd{\rho}\right|\ge L^{-1/4}}\nonumber \\
 & \le\P\rbr{2\left\Vert g\right\Vert _{\ns}c^{1/H-1}\TTV{W^{H}}{[1,1/L]}c\ge\bar{D}_{H}(1/L)L^{3/4}/\bar{D}_{H}}\nonumber \\
 & =\P\rbr{\TTV{W^{H}}{[1,1/L]}c\ge\bar{D}_{H}(1/L)c^{1-1/H}L^{3/4}/\rbr{2\bar{D}_{H}\left\Vert g\right\Vert _{\ns}}} \nonumber  \\
 & \le\tilde{A}_{H}\exp\rbr{-\rbr{\bar{B}_{H}/2}\rbr{L^{3/4}/\rbr{2\bar{D}_{H}\left\Vert g\right\Vert _{\ns}}}^{1+\min(2H,1)}}.\label{eq:first_local}
\end{align}
Now, having established (\ref{eq:first_local}), (\ref{eq:mart_local})
and (\ref{eq:cond_exp_local}) we obtain the following fact.

\begin{fact}\label{fact_local_u_molli} There exist universal positive constants
$Q_{H,\left\Vert g\right\Vert _{\ns}}$ and $R_{H,\left\Vert g\right\Vert _{\ns}}$,
depending on $H$ and $\left\Vert g\right\Vert _{\ns}$ $H$ only,
such that for $L=2,3,\ldots$, $c\le L^{-2H}$ and 
\begin{equation}
r\le\min\cbr{L^{-H}/2,L^{-H-(1/4)-1/(6+6\min(2H,1))}},\label{eq:r_cond}
\end{equation}
one has 
\begin{align}
 & \P\rbr{\left|\int_{-r/2}^{r/2}\psi_{r}(\rho)\int_{\R}c^{1/H-1}\cbr{U_{0,1}\rbr{c,W^{H}-a}-U_{0,1}\rbr{c,W^{H}-a-\rho}}g(a)\dd a\dd{\rho}\right|\ge3L^{-1/4}+\left|R\rbr{c,L,g,r}\right|}\nonumber \\
 & \le R_{H,\left\Vert g\right\Vert _{\ns}}L^{2}\exp\rbr{-Q_{H,\left\Vert g\right\Vert _{\ns}}L^{1/6}},\label{eq:fact_local_conc}
\end{align}
where $\left|R\rbr{c,L,g,\psi_{r}}\right|$ is estimated in formula
(\ref{eq:remainder}).
\end{fact}

\begin{coro}\label{cor_local_u_molli}
Let $\rbr{c_{n}}$ be a sequence of positive reals such that $c_{n}=O\rbr{n^{-\eta}}$
for some $\eta>0$ and let $r_{n}=c_{n}^{(1-H)\xi/2}$, $n\in\N$,
where $\xi\in(0,1)$. Then the difference 
\begin{equation}
\int_{\R}c_{n}^{1/H-1}U_{0,1}\rbr{c_{n},W^{H}-a}g(a)\dd a-\int_{-r_{n}/2}^{r_{n}/2}\cbr{\int_{\R}\psi_{r_{n}}(\rho)c_{n}^{1/H-1}U_{0,1}\rbr{c_{n},W^{H}-a-\rho}g(a)\dd a}\dd{\rho}\label{eq:cor_diff}
\end{equation}
 converges a.s. to $0$ as $n\ra+\ns$. 
 \end{coro}
{\bf Proof.} There exists some $\kappa>0$ such that $c_{n}\le\kappa n^{-\eta}$.
Let us define 
\[
L_{n}:=\left\lfloor c_{n}^{-(1-H)\xi/(2H+2)}\right\rfloor \ge\left\lfloor \kappa^{-(1-H)\xi/(2H+2)}n^{\eta(1-H)\xi/(2H+2)}\right\rfloor 
\]
and let us fix $\varepsilon>0$. For sufficiently large $n\in\N$
we have $L_{n}\ge2$, 
\[
L_{n}^{-2H}\ge c_{n}^{2H(1-H)\xi/(2H+2)}=c_{n}^{H(1-H)\xi/(H+1)}\ge c_{n},
\]
\[
r_{n}=c_{n}^{(1-H)\xi/2}=\rbr{c_{n}^{-(1-H)\xi/(2H+2)}}^{-H-1}\le L_{n}^{-H-1}\le L_{n}^{-H}/2,
\]
\[
3L_{n}^{-1/4}\le\varepsilon/2
\]
and
\begin{align*}
\left|R\rbr{c_{n},L_{n},g,r_{n}}\right| & \le\rbr{2\max_{t\in[0,1]}\rbr{W_{t}^{H}}+r_{n}}\left\Vert g\right\Vert _{\ns}c_{n}^{1/H-1}L_{n}^{2}\\
 & \le\rbr{2\max_{t\in[0,1]}\rbr{W_{t}^{H}}+1}\left\Vert g\right\Vert _{\ns}c_{n}^{1/H-1}c_{n}^{-2(1-H)\xi/(2H+2)}\\
 & =\rbr{2\max_{t\in[0,1]}\rbr{W_{t}^{H}}+1}\left\Vert g\right\Vert _{\ns}c_{n}^{1/H-1}c_{n}^{-2(1-H)\xi/(2H+2)}\\
 & \le\rbr{2\max_{t\in[0,1]}\rbr{W_{t}^{H}}+1}\left\Vert g\right\Vert _{\ns}c_{n}^{(1/H-1)/(H+1)}\le\varepsilon/2.
\end{align*}
Thus, for such $n$, (\ref{eq:fact_local_conc}) gives an estimate
for the probability that 
\[
\left|\int_{-r/2}^{r/2}\psi_{r}(\rho)\int_{\R}c^{1/H-1}\cbr{U_{0,1}\rbr{c,W^{H}-a}-U_{0,1}\rbr{c,W^{H}-a-\rho}}g(a)\dd a\dd{\rho}\right|\ge\varepsilon.
\]
Since the sum 
\[
\sum_{n\in\N}R_{H,\left\Vert g\right\Vert _{\ns}}L_{n}^{2}\exp\rbr{-Q_{H,\left\Vert g\right\Vert _{\ns}}L_{n}^{1/6}}
\]
is finite and $\varepsilon$ may be as close to $0$ as we please,
by the Borel-Cantelli lemma the difference (\ref{eq:cor_diff}) tends
to $0$ a.s. 
\hfill $\blacksquare$

\subsection{Convergence of the difference of mollified versions\label{subsec:conv_mollified}}

In this subsection we will establish the a.s. convergence of the difference
of mollified versions of $U$ and ${\cal L}$. More precisely, we
will prove the following fact.

\begin{fact} \label{fact_molli_diff} Assume that the terms of the sequence
$\rbr{c_{n}}$ are positive and tend to $0$ in such a way that
$c_{n}=O\rbr{n^{-\eta}}$ for some $\eta>0$ and let $r_{n}=c_{n}^{\xi}$,
$n\in\N$, where $\xi\in(0,1)$. Then the difference
\begin{align*}
 & \int_{-r_{n}/2}^{r_{n}/2}\cbr{\int_{\R}\psi_{r_{n}}(\rho)c_{n}^{1/H-1}U_{0,1}\rbr{c_{n},W^{H}-a-\rho}g(a)\dd a}\dd{\rho}\\
 & -\frac{\mathfrak{c}_{H}}{2}\int_{-r_{n}/2}^{r_{n}/2}\psi_{r_{n}}(\rho)\left\{ \int_{\R}{\cal L}_{1}^{a-\rho}\rbr{W^{H}}g(a)\dd a\right\} \rho,
\end{align*}
where $\mathfrak{c}_{H}$ is defined in Corollary
\ref{cor_TV_limit}, tends a.s. to $0$ as $n\ra+\ns$.
\end{fact}
{\bf Proof.} First, let us notice that for any $c>0$ and $r>0$ 
\begin{align*}
 & \int_{-r/2}^{r/2}\cbr{\int_{\R}\psi(\rho)c^{1/H-1}U_{0,1}\rbr{c,W^{H}-a-\rho}g(a)\dd a}\dd{\rho}\\
 & =\left\{ b=a+\rho\right\} =\int_{-r/2}^{r/2}\cbr{\int_{\R}\psi(\rho)c^{1/H-1}U_{0,1}\rbr{c,W^{H}-b}g(b-\rho)\dd b}\dd{\rho}\\
 & =\int_{\R}c^{1/H-1}U_{0,1}\rbr{c,W^{H}-b}\cbr{\int_{-r/2}^{r/2}\psi(\rho)g(b-\rho)\dd{\rho}}\dd b\\
 & =\int_{\R}c^{1/H-1}U_{0,1}\rbr{c,W^{H}-a}\bar{g_{r}}(a)\dd a
\end{align*}
and similarly
\begin{align*}
\int_{-r/2}^{r/2}\psi(\rho)\left\{ \int_{\R}{\cal L}_{1}^{a-\rho}\rbr{W^{H}}g(a)\dd a\right\} \dd r & =\int_{\R}{\cal L}_{1}^{a}\rbr{W^{H}}\bar{g_{r}}(a)\dd a,
\end{align*}
where 
\[
\bar{g}_{r}(a)=\int_{-r/2}^{r/2}\psi_{r}(\rho)g(a-\rho)\dd{\rho}
\]
is as smooth as $\psi_{1}$. 

Now, in order to estimate 
\begin{align*}
 & \int_{\R}c_{n}^{1/H-1}U_{0,1}\rbr{c_{n},W^{H}-a}\bar{g}_{r_{n}}(a)\dd a-\frac{\mathfrak{c}_{H}}{2}\int_{\R}{\cal L}_{1}^{a}\rbr{W^{H}}\bar{g}_{r_{n}}(a)\dd a\\
 & =\int_{\R}\bar{g}_{r_{n}}(a)c_{n}^{1/H-1}U_{0,1}\rbr{c_{n},W^{H}-a}\dd a-\int_{0}^{1}\bar{g}_{r_{n}}\rbr{W_{t}^{H}}\dd t
\end{align*}
we decompose the time interval $[0,1]$ into sufficiently large number
of subintervals where $W^{H}$ does not change much. Let $\zeta\in(\xi,1)$.
For $n\in\N$ such that $r_{n}<1$ we set 
\[
\delta_{n}:=\frac{1}{\left\lfloor 1/r_{n}^{1/\rbr{H\cdot\zeta}}\right\rfloor }\begin{cases}
\ge r_{n}^{1/\rbr{H\cdot\zeta}},\\
\le2r_{n}^{1/\rbr{H\cdot\zeta}}.
\end{cases}
\]
Now, let $\lambda\in(H\cdot\zeta,H)$. What will be important for
us is that 
\[
\delta_{n}^{\lambda}=O\rbr{r_{n}^{\lambda/\rbr{H\cdot\zeta}}}=o\rbr{r_{n}}\text{ as }n\ra+\ns
\]
as well as the fact that 
\[
c_{n}^{1/H}=r_{n}^{1/\rbr{H\cdot\xi}}=\rbr{r_{n}^{1/\rbr{H\cdot\zeta}}}^{\zeta/\xi}=O\rbr{\delta_{n}^{\zeta/\xi}}=o\rbr{\delta_{n}}\text{ as }n\ra+\ns.
\]

Now we define times: 
\[
\tau_{k}^{n}=k\delta_{n},\ k=0,1,\ldots,N(n),
\]
where
\[
N(n)=\delta_{n}^{-1}=\left\lfloor 1/r_{n}^{1/\rbr{H\cdot\zeta}}\right\rfloor .
\]
and decompose
\begin{align}
 & \int_{\R}\bar{g}_{r_{n}}(a)c_{n}^{1/H-1}U_{0,1}\rbr{c_{n},W^{H}-a}\dd a\nonumber \\
 & =\int_{\R}\sum_{k=1}^{N(n)}\bar{g}_{r_{n}}(a)c_{n}^{1/H-1}U_{\tau_{k-1}^{n},\tau_{k}^{n}}\rbr{c_{n},W^{H}-a}\dd a\nonumber \\
 & \qquad+\int_{\R}\bar{g}_{r_{n}}(a)c_{n}^{1/H-1}\cbr{U_{0,1}\rbr{c_{n},W^{H}-a}-\sum_{k=1}^{N(n)}U_{\tau_{k-1}^{n},\tau_{k}^{n}}\rbr{c_{n},W^{H}-a}}\dd a\nonumber \\
 & =\sum_{k=1}^{N(n)}\bar{g}_{r_{n}}\rbr{W_{\tau_{k-1}^{n}}^{H}}\int_{\R}c_{n}^{1/H-1}U_{\tau_{k-1}^{n},\tau_{k}^{n}}\rbr{c_{n},W^{H}-a}\dd a\nonumber \\
 & \qquad+\int_{\R}\sum_{k=1}^{N(n)}\rbr{\bar{g}_{r_{n}}(a)-\bar{g}_{r_{n}}\rbr{W_{\tau_{k-1}^{n}}^{H}}}c_{n}^{1/H-1}U_{\tau_{k-1}^{n},\tau_{k}^{n}}\rbr{c_{n},W^{H}-a}\dd a\nonumber \\
 & \qquad+\int_{\R}\bar{g}_{r_{n}}(a)c_{n}^{1/H-1}\cbr{U_{0,1}\rbr{c_{n},W^{H}-a}-\sum_{k=1}^{N(n)}U_{\tau_{k-1}^{n},\tau_{k}^{n}}\rbr{c_{n},W^{H}-a}}\dd a.\label{eq:decomp_loc_moll}
\end{align}
We will prove that the last two obtained terms tend to $0$ a.s. 

To estimate 
\[
\left|\int_{\R}\sum_{k=1}^{N(n)}\rbr{\bar{g}_{r_{n}}(a)-\bar{g}_{r_{n}}\rbr{W_{\tau_{k-1}^{n}}^{H}}}c_{n}^{1/H-1}U_{\tau_{k-1}^{n},\tau_{k}^{n}}\rbr{c_{n},W^{H}-a}\dd a\right|
\]
we use H\"older continuity of $W^{H}$ and write 
\begin{align}
 & \left|\bar{g}_{r_{n}}(a)-\bar{g}_{r_{n}}\rbr{W_{\tau_{k-1}^{n}}^{H}}\right|\nonumber \\
 & \le\sup_{b\in\R}\left|\bar{g}_{r_{n}}^{'}(b)\right|\sup_{\left|u-v\right|\le\delta_{n}}\left|W_{v}^{H}-W_{u}^{H}\right| \nonumber \\
 & =\sup_{b\in\R}\left|\bar{g}_{r_{n}}^{'}(b)\right|\delta_{n}^{\lambda}\sup_{\left|u-v\right|\le\delta_{n}}\frac{\left|W_{v}^{H}-W_{u}^{H}\right|}{\delta_{n}^{\lambda}}, \label{eq:estim_123}
\end{align}
where $\sup_{\left|u-v\right|\le\delta_{n}}\left|W_{v}^{H}-W_{u}^{H}\right|\delta_{n}^{-\lambda}$
is a.s. finite. Let us notice that for any $r>0$ and $b\in\R$,
\begin{align*}
\left|\bar{g}_{r}^{'}(b)\right| & =\left|\int_{-r/2}^{r/2}\psi_{r}^{'}(\rho)g(b-\rho)\dd{\rho}\right|=\left|\int_{-r/2}^{r/2}\frac{1}{r^{2}}\psi_{1}^{'}(\rho/r)g(b-\rho)\dd{\rho}\right|\\
 & =\left\{ u=\rho/r,\dd u=(1/r)\dd{\rho}\right\} =\left|\int_{-1/2}^{1/2}\frac{1}{r^{2}}\psi_{1}^{'}(u)g(b-u\cdot r)r\dd u\right|\\
 & \le\frac{1}{r}\int_{-1/2}^{1/2}\left|\psi_{1}^{'}(u)g(b-u\cdot r)\right|\dd u\le\frac{1}{r}\left\Vert g\right\Vert _{\ns}\int_{-1/2}^{1/2}\left|\psi_{1}^{'}(u)\right|\dd u,
\end{align*}
which yields
\begin{equation}
\sup_{b\in\R}\left|\bar{g}_{r}^{'}(b)\right|\le\frac{1}{r}\left\Vert g\right\Vert _{\ns}\int_{-1/2}^{1/2}\left|\psi_{1}^{'}(u)\right|\dd u.\label{eq:psi_prim_est}
\end{equation}
By (\ref{eq:estim_123}), (\ref{eq:psi_prim_est}) and the analog
of (\ref{eq:sup_2}) for $\UTV{W^{H}}{[\cdot,\cdot]}{c_{n}}$, 
\begin{align*}
 & \left|\int_{\R}\sum_{k=1}^{N(n)}\rbr{\bar{g}_{r_{n}}(a)-\bar{g}_{r_{n}}\rbr{W_{\tau_{k-1}^{n}}^{H}}}c_{n}^{1/H-1}U_{\tau_{k-1}^{n},\tau_{k}^{n}}\rbr{c_{n},W^{H}-a}\dd a\right|\\
 & \le\frac{\delta_{n}^{\lambda}}{r_{n}}\rbr{\sup_{\left|u-v\right|\le\delta_{n}}\frac{\left|W_{v}^{H}-W_{u}^{H}\right|}{\delta_{n}^{\lambda}}}\left\Vert g\right\Vert _{\ns}\rbr{\int_{-1/2}^{1/2}\left|\psi_{1}^{'}(u)\right|\dd u}\int_{\R}c_{n}^{1/H-1}\sum_{k=1}^{N(n)}U_{\tau_{k-1}^{n},\tau_{k}^{n}}\rbr{c_{n},W^{H}-a}\dd a\\
 & \le\frac{\delta_{n}^{\lambda}}{r_{n}}\rbr{\sup_{\left|u-v\right|\le\delta_{n}}\frac{\left|W_{v}^{H}-W_{u}^{H}\right|}{\delta_{n}^{\lambda}}}\left\Vert g\right\Vert _{\ns}\rbr{\int_{-1/2}^{1/2}\left|\psi_{1}^{'}(u)\right|\dd u}c_{n}^{1/H-1}\UTV{W^{H}}{[0,1]}{c_{n}}.
\end{align*}
Since $\sup_{\left|u-v\right|\le\delta_{n}}\left|W_{v}^{H}-W_{u}^{H}\right|\delta_{n}^{-\lambda}$
is a.s. finite,
\[
c_{n}^{1/H-1}\UTV{W^{H}}{[0,1]}{c_{n}}
\]
tends a.s. to a finite number $\mathfrak{c}_{H}/2$ and $\delta_{n}^{\lambda}/r_{n}\ra0$
as $n\ra+\ns$, we get that 
\begin{equation}
\left|\int_{\R}\sum_{k=1}^{N(n)}\rbr{\bar{g}_{r_{n}}(a)-\bar{g}_{r_{n}}\rbr{W_{\tau_{k-1}^{n}}^{H}}}c_{n}^{1/H-1}U_{\tau_{k-1}^{n},\tau_{k}^{n}}\rbr{c_{n},W^{H}-a}\dd a\right|\ra0\text{ a.s.}\label{eq:conv_loc1}
\end{equation}

Now we will estimate 
\[
\int_{\R}\bar{g}_{r_{n}}(a)c_{n}^{1/H-1}\cbr{U_{0,1}\rbr{c_{n},W^{H}-a}-\sum_{k=1}^{N(n)}U_{\tau_{k-1}^{n},\tau_{k}^{n}}\rbr{c_{n},W^{H}-a}}\dd a.
\]
Using (\ref{eq:sup-1}) we get that 
\[
U_{0,1}\rbr{c_{n},W^{H}-a}-\sum_{k=1}^{N(n)}U_{\tau_{k-1}^{n},\tau_{k}^{n}}\rbr{c_{n},W^{H}-a}\ge0.
\]
This observation together with (\ref{eq:Banach_ind_gen_0}) and the
analog of (\ref{eq:sub_2}) for $\UTV{W^{H}}{[\cdot,\cdot]}{c_{n}}$
yield
\begin{align*}
 & \left|\int_{\R}\bar{g}_{r_{n}}(a)c_{n}^{1/H-1}\cbr{U_{0,1}\rbr{c_{n},W^{H}-a}-\sum_{k=1}^{N(n)}U_{\tau_{k-1}^{n},\tau_{k}^{n}}\rbr{c_{n},W^{H}-a}}\dd a\right|\\
 & \le\int_{\R}\left|\bar{g}_{r_{n}}(a)\right|c_{n}^{1/H-1}\cbr{U_{0,1}\rbr{c_{n},W^{H}-a}-\sum_{k=1}^{N(n)}U_{\tau_{k-1}^{n},\tau_{k}^{n}}\rbr{c_{n},W^{H}-a}}\dd a\\
 & \le\left\Vert g\right\Vert _{\ns}c_{n}^{1/H-1}\int_{\R}\cbr{U_{0,1}\rbr{c_{n},W^{H}-a}-\sum_{k=1}^{N(n)}U_{\tau_{k-1}^{n},\tau_{k}^{n}}\rbr{c_{n},W^{H}-a}}\dd a\\
 & =\left\Vert g\right\Vert _{\ns}c_{n}^{1/H-1}\cbr{\UTV{W^{H}}{[0,1]}{c_{n}}-\sum_{k=1}^{N(n)}\UTV{W^{H}}{\sbr{\tau_{k-1}^{n},\tau_{k}^{n}}}{c_{n}}}\\
 & \le\left\Vert g\right\Vert _{\ns}c_{n}^{1/H-1}N(n)c_{n}=\left\Vert g\right\Vert _{\ns}c_{n}^{1/H}N(n)=\left\Vert g\right\Vert _{\ns}c_{n}^{1/H}\delta_{n}^{-1}.
\end{align*}
Since $c_{n}^{1/H}\delta_{n}^{-1}\ra0$ as $n\ra+\ns$, hence
\begin{equation}
\left|\int_{\R}\bar{g}_{r_{n}}(a)c_{n}^{1/H-1}\cbr{U_{0,1}\rbr{c_{n},W^{H}-a}-\sum_{k=1}^{N(n)}U_{\tau_{k-1}^{n},\tau_{k}^{n}}\rbr{c_{n},W^{H}-a}}\dd a\right|\ra0\text{ a.s.}\label{eq:conv_loc2}
\end{equation}
as $n\ra+\ns$.

Finally, let us notice that the first term in the last eq. in (\ref{eq:decomp_loc_moll}),
i.e. 
\[
\sum_{k=1}^{N(n)}\bar{g}_{r_{n}}\rbr{W_{\tau_{k-1}^{n}}^{H}}\int_{\R}c_{n}^{1/H-1}U_{\tau_{k-1}^{n},\tau_{k}^{n}}\rbr{c_{n},W^{H}-a}\dd a,
\]
tends a.s. to $\rbr{\mathfrak{c}_{H}/2}\int_{0}^{1}\bar{g}_{r}\rbr{W_{t}^{H}}\dd t$
as $c\ra0_{+}$. To see this we write 
\begin{align*}
 & \sum_{k=1}^{N(n)}\bar{g}_{r_{n}}\rbr{W_{\tau_{k-1}^{n}}^{H}}\int_{\R}c_{n}^{1/H-1}U_{\tau_{k-1}^{n},\tau_{k}^{n}}\rbr{c_{n},W^{H}-a}\dd a-\frac{\mathfrak{c}_{H}}{2}\int_{0}^{1}\bar{g}_{r}\rbr{W_{t}^{H}}\dd t\\
 & =\sum_{k=1}^{N(n)}\bar{g}_{r_{n}}\rbr{W_{\tau_{k-1}^{n}}^{H}}\int_{\R}c_{n}^{1/H-1}\cbr{U_{\tau_{k-1}^{n},\tau_{k}^{n}}\rbr{c_{n},W^{H}-a}-\E U_{\tau_{k-1}^{n},\tau_{k}^{n}}\rbr{c_{n},W^{H}-a}}\dd a\\
 & \qquad+\sum_{k=1}^{N(n)}\bar{g}_{r_{n}}\rbr{W_{\tau_{k-1}^{n}}^{H}}\cbr{c_{n}^{1/H-1}\int_{\R}\E U_{\tau_{k-1}^{n},\tau_{k}^{n}}\rbr{c_{n},W^{H}-a}\dd a-\frac{\mathfrak{c}_{H}}{2}\rbr{\tau_{k}-\tau_{k-1}}}\\
 & \qquad+\frac{\mathfrak{c}_{H}}{2}\sum_{k=1}^{N(n)}\bar{g}_{r_{n}}\rbr{W_{\tau_{k-1}^{n}}^{H}}\rbr{\tau_{k}^{n}-\tau_{k-1}^{n}}-\frac{\mathfrak{c}_{H}}{2}\int_{0}^{1}\bar{g}_{r_{n}}\rbr{W_{t}^{H}}\dd t.
\end{align*}

The first term
\[
\sum_{k=1}^{N(n)}\bar{g}_{r_{n}}\rbr{W_{\tau_{k-1}^{n}}^{H}}\int_{\R}c_{n}^{1/H-1}\cbr{U_{\tau_{k-1}^{n},\tau_{k}^{n}}\rbr{c_{n},W^{H}-a}-\E U_{\tau_{k-1}^{n},\tau_{k}^{n}}\rbr{c_{n},W^{H}-a}}\dd a
\]
may be estimated as
\begin{align*}
 & \left|\sum_{k=1}^{N(n)}\bar{g}_{r_{n}}\rbr{W_{\tau_{k-1}^{n}}^{H}}\int_{\R}c_{n}^{1/H-1}\cbr{U_{\tau_{k-1}^{n},\tau_{k}^{n}}\rbr{c_{n},W^{H}-a}-\E U_{\tau_{k-1}^{n},\tau_{k}^{n}}\rbr{c_{n},W^{H}-a}}\dd a\right|\\
 & \le\left\Vert g\right\Vert _{\ns}\sum_{k=1}^{N(n)}c_{n}^{1/H-1}\left|\UTV{W^{H}}{\sbr{\tau_{k-1}^{n},\tau_{k}^{n}}}{c_{n}}-\UTV{W^{H}}{\sbr{\tau_{k-1}^{n},\tau_{k}^{n}}}{c_{n}}\right|.
\end{align*}
Now, using (\ref{eq:Th12-1}) and Lemma A, since $c_{n}\le\delta_{n}^{H}$,
for any $\varepsilon>0$ we estimate
\begin{align*}
 & \P\rbr{\sum_{k=1}^{N(n)}c_{n}^{1/H-1}\left|\UTV{W^{H}}{\sbr{\tau_{k-1}^{n},\tau_{k}^{n}}}{c_{n}}-\UTV{W^{H}}{\sbr{\tau_{k-1}^{n},\tau_{k}^{n}}}{c_{n}}\right|\ge\varepsilon}\\
 & \le N(n)\P\rbr{c_{n}^{1/H-1}\left|\UTV{W^{H}}{\sbr{\tau_{k-1}^{n},\tau_{k}^{n}}}{c_{n}}-\UTV{W^{H}}{\sbr{\tau_{k-1}^{n},\tau_{k}^{n}}}{c_{n}}\right|\ge\varepsilon/N(n)}\\
 & =N(n)\P\rbr{c_{n}^{1/H-1}\left|\UTV{W^{H}}{\sbr{\tau_{k-1}^{n},\tau_{k}^{n}}}{c_{n}}-\UTV{W^{H}}{\sbr{\tau_{k-1}^{n},\tau_{k}^{n}}}{c_{n}}\right|\ge\delta_{n}\varepsilon}\\
 & \le N(n)\tilde{A}_{H}\exp\rbr{-\rbr{\bar{B}_{H}/2}\rbr{\delta_{n}c_{n}^{-1/H}}^{2-\max\rbr{2H,1}}\varepsilon^{1+\min\rbr{2H,1}}}\\
 & =\rbr{\tilde{A}_{H}/\delta_{n}}\exp\rbr{-\rbr{\bar{B}_{H}/2}\rbr{\delta_{n}c_{n}^{-1/H}}^{2-\max\rbr{2H,1}}\varepsilon^{1+\min\rbr{2H,1}}}
\end{align*}
from which we infer that 
\[
\sum_{k=1}^{N(n)}c_{n}^{1/H-1}\left|\UTV{W^{H}}{\sbr{\tau_{k-1}^{n},\tau_{k}^{n}}}{c_{n}}-\UTV{W^{H}}{\sbr{\tau_{k-1}^{n},\tau_{k}^{n}}}{c_{n}}\right|
\]
tends a.s. to $0$ as $c\ra0_{+}$. 

The second term may be estimated as follows
\begin{align*}
 & \left|\sum_{k=1}^{N(n)}\bar{g}_{r_{n}}\rbr{W_{\tau_{k-1}^{n}}^{H}}\cbr{c_{n}^{1/H-1}\int_{\R}\E U_{\tau_{k-1}^{n},\tau_{k}^{n}}\rbr{c_{n},W^{H}-a}\dd a-\frac{\mathfrak{c}_{H}}{2}\rbr{\tau_{k}^{n}-\tau_{k-1}^{n}}}\right|\\
 & \le\left\Vert g\right\Vert _{\ns}\sum_{k=1}^{N(n)}\left|c_{n}^{1/H-1}\E\UTV{W^{H}}{\sbr{\tau_{k-1}^{n},\tau_{k}^{n}}}{c_{n}}-\frac{\mathfrak{c}_{H}}{2}\rbr{\tau_{k}^{n}-\tau_{k-1}^{n}}\right|\\
 & =\left\Vert g\right\Vert _{\ns}\sum_{k=1}^{N(n)}\left|c_{n}^{1/H-1}\E\UTV{W^{H}}{\sbr{\tau_{k-1}^{n},\tau_{k}^{n}}}{c_{n}}-\frac{\mathfrak{c}_{H}}{2}\delta_{n}\right|.
\end{align*}
But we know, by scaling (recall the proof of Corollary \ref{cor_TV_limit}),
that
\[
c_{n}^{1/H-1}\E\UTV{W^{H}}{\sbr{\tau_{k-1}^{n},\tau_{k}^{n}}}{c_{n}}=\frac{\E\UTV{W^{H}}{\sbr{0,c_{n}^{-1/H}\delta_{n}}}1}{c_{n}^{-1/H}\delta_{n}}\delta_{n}
\]
moreover, by superadditivity of $d\mapsto\E\UTV{W^{H}}{\sbr{0,c_{n}^{-1/H}\delta_{n}d}}1$,
$d\in\N$, and subadditivity of $d\mapsto\E\UTV{W^{H}}{\sbr{0,c_{n}^{-1/H}\delta_{n}d}}1+1$,
$d\in\N$, we have, (Fekete's lemma), that
\[
\mathfrak{d}_{H}:=\lim_{d\ra+\ns}\frac{\E\UTV{W^{H}}{\sbr{0,c_{n}^{-1/H}\delta_{n}d}}1}{d}
\]
exists and satisfies 
\[
\sup_{d\in\N}\frac{\E\UTV{W^{H}}{\sbr{0,c_{n}^{-1/H}\delta_{n}d}}1}{d}\le\mathfrak{d}_{H}\le\inf_{d\in\N}\frac{\E\UTV{W^{H}}{\sbr{0,c_{n}^{-1/H}\delta_{n}d}}1+1}{d}.
\]
In particular,
\[
\E\UTV{W^{H}}{\sbr{0,c_{n}^{-1/H}\delta_{n}}}1\le\mathfrak{d}_{H}\le\E\UTV{W^{H}}{\sbr{0,c_{n}^{-1/H}\delta_{n}}}1+1.
\]
Notice, that 
\[
\mathfrak{d}_{H}=c_{n}^{-1/H}\delta_{n}\frac{\mathfrak{c}_{H}}{2}
\]
thus
\[
\frac{\E\UTV{W^{H}}{\sbr{0,c_{n}^{-1/H}\delta_{n}}}1}{c_{n}^{-1/H}}\le\frac{\mathfrak{c}_{H}}{2}\delta_{n}\le\frac{\E\UTV{W^{H}}{\sbr{0,c_{n}^{-1/H}\delta_{n}}}1}{c_{n}^{-1/H}}+\frac{1}{c_{n}^{-1/H}},
\]
hence 
\begin{align*}
\left|c_{n}^{1/H-1}\E\UTV{W^{H}}{\sbr{\tau_{k-1}^{n},\tau_{k}^{n}}}{c_{n}}-\frac{\mathfrak{c}_{H}}{2}\delta_{n}\right| & =\left|\frac{\E\UTV{W^{H}}{\sbr{0,c_{n}^{-1/H}\delta_{n}}}1}{c_{n}^{-1/H}\delta_{n}}\delta_{n}-\frac{\mathfrak{c}_{H}}{2}\delta_{n}\right|\le c_{n}^{1/H}
\end{align*}
and we have the estimate
\begin{align*}
 & \left|\sum_{k=1}^{N(n)}\bar{g}_{r_{n}}\rbr{W_{\tau_{k-1}^{n}}^{H}}\cbr{c_{n}^{1/H-1}\int_{\R}\E U_{\tau_{k-1}^{n},\tau_{k}^{n}}\rbr{c_{n},W^{H}-a}\dd a-\frac{\mathfrak{c}_{H}}{2}\rbr{\tau_{k}^{n}-\tau_{k-1}^{n}}}\right|\\
 & \le\left\Vert g\right\Vert _{\ns}\sum_{k=1}^{N(n)}\left|c_{n}^{1/H-1}\E\UTV{W^{H}}{\sbr{\tau_{k-1}^{n},\tau_{k}^{n}}}{c_{n}}-\frac{\mathfrak{c}_{H}}{2}\delta_{n}\right|\\
 & \le\left\Vert g\right\Vert _{\ns}\sum_{k=1}^{N(n)}c_{n}^{1/H}=\left\Vert g\right\Vert _{\ns}c_{n}^{1/H}N(n)=\left\Vert g\right\Vert _{\ns}c_{n}^{1/H}\delta_{n}^{-1}\ra0\text{ as }n\ra+\ns.
\end{align*}

The third term may be estimated, using (\ref{eq:psi_prim_est}) in
the following way: 
\begin{align*}
 & \left|\frac{\mathfrak{c}_{H}}{2}\sum_{k=1}^{N(n)}\bar{g}_{r_{n}}\rbr{W_{\tau_{k-1}^{n}}^{H}}\rbr{\tau_{k}-\tau_{k-1}}-\frac{\mathfrak{c}_{H}}{2}\int_{0}^{1}\bar{g}_{r_{n}}\rbr{W_{t}^{H}}\dd t\right|\\
 & =\frac{\mathfrak{c}_{H}}{2}\left|\int_{0}^{1}\hat{g}_{r_{n}}\rbr t-\bar{g}_{r_{n}}\rbr{W_{t}^{H}}\dd t\right|\dd t\\
 & \le\frac{\mathfrak{c}_{H}}{2}\int_{0}^{1}\sup_{b\in\R}\left|\bar{g}_{r_{n}}^{'}(b)\right|\sup_{\left|u-v\right|\le\delta_{n}}\left|W_{v}^{H}-W_{u}^{H}\right|\dd t\\
 & =\frac{\mathfrak{c}_{H}}{2}\left|\bar{g}_{r}^{'}(b)\right|\sup_{b\in\R}\left|\bar{g}_{r_{n}}^{'}(b)\right|\delta_{n}^{\lambda}\sup_{\left|u-v\right|\le\delta_{n}}\frac{\left|W_{v}^{H}-W_{u}^{H}\right|}{\delta_{n}^{\lambda}}\\
 & \le\frac{\mathfrak{c}_{H}}{2}\sup_{\left|u-v\right|\le\delta_{n}}\frac{\left|W_{v}^{H}-W_{u}^{H}\right|}{\delta_{n}^{\lambda}}\left\Vert g\right\Vert _{\ns}\rbr{\int_{-1/2}^{1/2}\left|\psi_{1}^{'}(u)\right|\dd u}\frac{\delta_{n}^{\lambda}}{r_{n}},
\end{align*}
where 
\[
\hat{g}_{r_{n}}\rbr t=\sum_{k=1}^{N(n)}\bar{g}_{r_{n}}\rbr{W_{\tau_{k-1}^{n}}^{H}}{\bf 1}_{\left[\tau_{k-1}^{n},\tau_{k}^{n}\right)}(t)+\bar{g}_{r_{n}}\rbr{W_{\tau_{N(n)-1}}^{H}}{\bf 1}_{\cbr 1}(t)
\]
is a stepwise approximation of $\bar{g}_{r_{n}}$. 
\hfill $\blacksquare$

\subsection{Almost sure weak convergence (in $L^{1}\protect\rbr{\protect\R}$)
of normalized strip crossings to the local time of fBm}

Having established (along any sequence of positive reals $\rbr{c_{n}}$
such that $c_{n}=O\rbr{n^{-\eta}}$, $n\in\N$, for some $\eta>0$
and proper choice of $r_{n}s$):
\begin{itemize}
\item a.s. convergence of the difference (\ref{eq:dife1}) in Corollary
\ref{cor_local_u_molli},
\item a.s. convergence of the difference (\ref{eq:dife2}) at the beginning
of this section,
\item a.s. convergence of the difference of mollified versions of $U$ and
${\cal L}$ in Fact \ref{fact_molli_diff}, 
\end{itemize}
we have the following result.

\begin{fact} \label{fact_local_final} Assume that the terms of the sequence
$\rbr{c_{n}}$ are positive and tend to $0$ in such a way that
$c_{n}=O\rbr{n^{-\eta}}$ for some $\eta>0$, then for any $g\in L^{\ns}\rbr{\R},$
\[
\int_{\R}c_{n}^{1/H-1}U_{0,1}\rbr{c_{n},W^{H}-a}g(a)\dd a
\]
converges a.s. to
\[
\frac{\mathfrak{c}_{H}}{2}\int_{\R}{\cal L}_{1}^{a}g(a)\dd a=\frac{\mathfrak{c}_{H}}{2}\int_{0}^{1}g\rbr{W_{t}^{H}}\dd t
\]
as $n\ra+\ns$, where$\mathfrak{c}_{H}$ is defined in Corollary \ref{cor_TV_limit}.
\end{fact}

Taking advantage of the fact that for any $a\in\R$ the mapping $(0,\ns)\ni c\mapsto U_{0,1}\rbr{c,W^{H}-a}$
is non-increasing (this is not the case of $K_{0,1}\rbr{c,W^{H}-a}$,
see Remark \ref{rem_Kst_not_monotonic}), we may strengthen Fact \ref{fact_local_final}
and obtain the following theorem.

\begin{theo} \label{th_final_local}  For any $g\in L^{\ns}\rbr{\R}$
and any sequence $\rbr{c_{n}}$ of positive reals tending to $0$
\[
\int_{\R}c_{n}^{1/H-1}U_{0,1}\rbr{c_{n},W^{H}-a}g(a)\dd a
\]
converges a.s. to
\[
\frac{\mathfrak{c}_{H}}{2}\int_{\R}{\cal L}_{1}^{a}g(a)\dd a=\frac{\mathfrak{c}_{H}}{2}\int_{0}^{1}g\rbr{W_{t}^{H}}\dd t.
\]
\end{theo}
{\bf Proof.} Let us decompose 
\[
g=\max\rbr{g,0}-\max\rbr{-g,0}.
\]
By Fact \ref{fact_local_final}, 
\begin{equation}
\int_{\R}(1/n)^{1/H-1}U_{0,1}\rbr{1/n,W^{H}-a}\max\rbr{g(a),0}\dd a\ra\frac{\mathfrak{c}_{H}}{2}\int_{0}^{1}\max\rbr{g\rbr{W_{t}^{H}},0}\dd t\label{eq:limit}
\end{equation}
a.s. as $n\ra+\ns$. Now, for any sequence of positive
reals $\rbr{c_{n}}$ tending to $0$, for those $c_{n}$ which are
smaller than $1$, we have $1/\left\lfloor 1/c_{n}\right\rfloor \ge c_{n}\ge1/\left\lceil 1/c_{n}\right\rceil $
and similarly as in the proof of Corollary \ref{cor_TV_limit}, using
the monotonicity of $(0,\ns)\ni c\mapsto U_{0,1}\rbr{c,W^{H}-a}$
and positivity of $\max\rbr{g(a),0}$ we have 
\begin{align*}
 & \rbr{1/\left\lfloor 1/c_{n}\right\rfloor }^{1/H-1}\int_{\R}U_{0,1}\rbr{1/\left\lceil 1/c_{n}\right\rceil ,W^{H}-a}\max\rbr{g(a),0}\dd a\\
 & \ge c_{n}^{1/H-1}\int_{\R}U_{0,1}\rbr{c_{n},W^{H}-a}\max\rbr{g(a),0}\dd a\\
 & \ge\rbr{1/\left\lceil 1/c_{n}\right\rceil }^{1/H-1}\int_{\R}U_{0,1}\rbr{1/\left\lfloor 1/c_{n}\right\rfloor ,W^{H}-a}\max\rbr{g(a),0}\dd a.
\end{align*}
However, since 
\[
\lim_{n\ra+\ns}\frac{\rbr{1/\left\lfloor 1/c_{n}\right\rfloor }^{1/H-1}}{\rbr{1/\left\lceil 1/c_{n}\right\rceil }^{1/H-1}}=1,
\]
all terms converge to the same limit as 
\[
\rbr{1/\left\lfloor 1/c_{n}\right\rfloor }^{1/H-1}\int_{\R}U_{0,1}\rbr{1/\left\lfloor 1/c_{n}\right\rfloor ,W^{H}-a}\max\rbr{g(a),0}\dd a
\]
which coincides with the limit in (\ref{eq:limit}), i.e. $\rbr{\mathfrak{c}_{H}/{2}}\int_{0}^{1}\max\rbr{g\rbr{W_{t}^{H}},0}\dd t$. 

Similar reasoning applies naturally to $\max\rbr{-g,0}$ and taking
the difference $\max\rbr{g,0}-\max\rbr{-g,0}$ we get the claimed
convergence of $\int_{\R}c_{n}^{1/H-1}U_{0,1}\rbr{c_{n},W^{H}-a}g(a)\dd a$
for any $g\in L^{\ns}\rbr{\R}.$

\hfill $\blacksquare$

\section*{Appendix}

\subsection*{Estimates of expectation of $ \UTV{W^H}{S}{c}$}
 From self-similarity $W^H_s \sim cW^H_{c^{-1/H}s}$, $s \ge 0$, and subadditivity of $\UTV{f}{[\cdot,\cdot]}{c}+c$ as a function of the interval:
$$\UTV{f}{[a,b]}{c} + c +\UTV{f}{[b,d]}{c} + c \gs \UTV{f}{[a,d]}{c} +c,$$ we get  
\begin{align*}
\E \UTV{W^H}{S}{c}& =  \E \UTV{c \cdot W^H}{c^{-1/H}S}{c} =  c \cdot \E \UTV{W^H}{c^{-1/H}S}{1} \\
&\ls  c \sum_{k=1}^{\lfloor c^{-1/H}S \rfloor +1}\cbr{ \E \UTV{W^H}{\sbr{k-1,k}}{1} +1}\\
&= c \rbr{\lfloor c^{-1/H}S \rfloor +1} \cbr{\E \UTV{W^H}{1}{1}+1} \\
& \ls 2 S c^{1-1/H} \cbr{\E \UTV{W^H}{1}{1}+1}.
\end{align*}
for $c$ such that $c^{-1/H}S \gs 1 $, which is equivalent with $c \ls S^H$.
Thus
\begin{equation} \label{eutv}
\E \UTV{W^H} {S}{c}  \ls C_H S c^{1-1/H}, \ \text{ for } c\ls S^H.
\end{equation}
On the other hand, from superadditivity of $\UTV{f}{[\cdot,\cdot]}{c}$ as a function of the interval:
$$\UTV{f}{[a,b]}{c}  +\UTV{f}{[b,d]}{c} \ls \UTV{f}{[a,d]}{c}, $$ we get  
\begin{align*}
\E \UTV{W^H}{S}{c}& =  \E \UTV{c \cdot W^H}{c^{-1/H}S}{c} =  c \cdot \E \UTV{W^H}{c^{-1/H}S}{1} \\
&\gs  c \sum_{k=1}^{\lfloor c^{-1/H}S \rfloor } \E \UTV{W^H}{\sbr{k-1,k}}{1} \\
&= c \rbr{\lfloor c^{-1/H}S \rfloor } \E \UTV{W^H}{1}{1} \\
& \gs \frac{1}{2} S c^{1-1/H} \E \UTV{W^H}{1}{1}
\end{align*}
for $c$ such that $c^{-1/H}S \gs 1 $. As a result, for some positive $\tilde{C}_H$,
\begin{equation} \label{eutvlower}
\E \UTV{W^H} {S}{c}  \gs \tilde{C}_H Sc^{1-1/H}, \ \text{ for } c\ls S^H.
\end{equation}

\subsection*{Weak variance of sums of increments of fractional Brownian motion}
In this subsection we will estimate the weak second moment of the sums of $n$ increments ($n=1,2,\ldots$) of fractional Brownian motion, that is
\begin{align*}
V_n\rbr{W^H,S} & : = \sup_{0 \ls s_1 < t_1 < s_2 < t_2 < \dots < s_n <t_n \ls S} \E \rbr{ \sum_{i=1}^n \rbr{ W_{t_{i}}^{H}-W_{s_{i}}^{H}}}^2 \\
& = \sup_{0 \ls s_1 < t_1 < s_2 < t_2 < \dots < s_n <t_n \ls S}  \sum_{i=1}^n  \sum_{j=1}^n   \E \rbr{ W_{t_{i}}^{H}-W_{s_{i}}^{H}}\rbr{ W_{t_{j}}^{H}-W_{s_{j}}^{H}}.
\end{align*}
From the formula $
\E \rbr{W^H_s W^H_t} = \frac{1}{2}\rbr{|s|^{2H}+ |t|^{2H} - |t-s|^{2H}}
$
we infer for $1\ls i < j \ls n$
\be \label{covH}
\E \rbr{ W_{t_{i}}^{H}-W_{s_{i}}^{H}}\rbr{ W_{t_{j}}^{H}-W_{s_{j}}^{H}} = \frac{1}{2}\rbr{|t_j - s_i|^{2H}+ |t_i - s_j|^{2H} - |s_j-s_i|^{2H} - |t_j-t_i|^{2H} }.
\ee
\subsubsection*{The case $0 < H<1/2$}
{\bf Lower bound}. To estimate $V_n\rbr{W^H,S}$ from below we consider the variable
\[
S_{n} := \sum_{i=1}^n \rbr{ W_{2iS/(2n)}^{H}-W_{(2i-1)S/(2n)}^{H}}
\]
and, using \eqref{covH},  calculate
\begin{align}
\E S_n^2 & = \sum_{i=1}^n  \sum_{j=1}^n   \E \rbr{ W_{2iS/(2n)}^{H}-W_{(2i-1)S/(2n)}^{H}}\rbr{ W_{(2j)S/(2n)}^{H}-W_{(2j-1)S/(2n)}^{H}}  \nonumber \\
& = \sum_{i=1}^n \E \rbr{ W_{2iS/(2n)}^{H}-W_{(2i-1)S/(2n)}^{H}}^2 + 2 \sum_{1\ls i < j \ls n}  \E \rbr{ W_{2iS/(2n)}^{H}-W_{(2i-1)S/(2n)}^{H}}\rbr{ W_{(2j)S/(2n)}^{H}-W_{(2j-1)S/(2n)}^{H}}  \nonumber \\
& = \sum_{i=1}^n \rbr{\frac{S}{2n}}^{2H} + 2 \sum_{m=1}^{n-1}  \sum_{1\ls i < j \ls n, j-i = m}  \E \rbr{ W_{2iS/(2n)}^{H}-W_{(2i-1)S/(2n)}^{H}}\rbr{ W_{(2j)S/(2n)}^{H}-W_{(2j-1)S/(2n)}^{H}}  \nonumber  \\
& = n\rbr{\frac{S}{2n}}^{2H} + 2 \sum_{m=1}^{n-1}  \sum_{1\ls i < j \ls n, j-i = m}   \frac{1}{2}\rbr{\rbr{\frac{2j - 2i +1}{2n}S}^{2H}+\rbr{\frac{2j - 2i - 1}{2n}S}^{2H} - 2 \rbr{\frac{2j - 2i }{2n}S}^{2H} }  \nonumber \\
& = n \rbr{\frac{S}{2n}}^{2H} +  \sum_{m=1}^{n-1} (n-m) \rbr{\rbr{\frac{2m+1}{2n}S}^{2H}+\rbr{\frac{2m - 1}{2n}S}^{2H} - 2 \rbr{\frac{2m}{2n}S}^{2H} }  \nonumber \\
& \ge n \rbr{\frac{S}{2n}}^{2H} +  \sum_{m=1}^{n-1} n \rbr{\frac{S}{2n}}^{2H} \rbr{\rbr{\frac{2m+1}{2n}S}^{2H}+\rbr{\frac{2m - 1}{2n}S}^{2H} - 2 \rbr{\frac{2m}{2n}S}^{2H} }  \nonumber
\\
& = n^{1-2H}S^{2H} \sbr{\frac{1}{2^{2H}} +  \sum_{m=1}^{n-1} m^{2H} \rbr{\rbr{1+\frac{1}{2m}}^{2H}+\rbr{1-\frac{1}{2m}}^{2H} - 2 }}. \label{ala}
\end{align} 
Let us consider $h:[-1/2,1/2] \ra \R$, 
$$h(x) := \rbr{1+x}^{2H}+\rbr{1-x}^{2H} - 2.$$
Expanding $g$ in the Taylor series around $0$ and taking into account that $0<2H <1$ and $|x|\le 1/2$, we estimate 
\begin{align*}
h(x) & = \frac{1}{2}\frac{2H(2H-1)}{2} x^2+ 2 \frac{2H(2H-1)(2H-2)(2H-3)}{4!}x^4 + \ldots \nonumber \\
& \gs 2H(2H-1) x^2+ 2 \frac{2H(2H-1)(-2)(-3)}{4!}x^4 + 2 \frac{2H(2H-1)(-2)(-3)(-4)(-5)}{6!}x^6 + \ldots \nonumber \\
& = 2H(2H-1) \sbr{ x^2+ \frac{1}{2}x^4 +\frac{1}{3}x^6  + \ldots } \gs  2H(2H-1) \sbr{ x^2+ \frac{1}{2}x^4\rbr{1 + x^2  + \ldots} }  \nonumber \\
& =  2H(2H-1) \sbr{ x^2+ \frac{1}{2}\frac{x^4}{1-x^2}} \gs  2H(2H-1) \sbr{ x^2+ \frac{1}{2}\frac{x^4}{1-1/4}} =  2H(2H-1) \sbr{ x^2+ \frac{2}{3}x^4}. 
\end{align*}  
Now, using this estimate and $2H(2H-1) \ge -1/4$ we get 
\begin{align*}
\frac{1}{2^{2H}} & +  \sum_{m=1}^{n-1} m^{2H} \rbr{\rbr{1+\frac{1}{2m}}^{2H}+\rbr{1-\frac{1}{2m}}^{2H} - 2 } \gs \frac{1}{2} +  \sum_{m=1}^{\infty} m^{2H} h\rbr{\frac{1}{2m}} \\
& \gs \frac{1}{2} +   2H(2H-1) \sum_{m=1}^{\infty} m^{2H}  \sbr{ \rbr{\frac{1}{2m}}^2+ \frac{2}{3}\rbr{\frac{1}{2m}}^4} \\
& \gs \frac{1}{2} +2H(2H-1) \rbr{ \frac{1}{4}+ \frac{2}{3}\frac{1}{16}+ \int_{1}^{\infty}  \frac{1}{4m^{2-2H}}+ \frac{2}{3}\frac{1}{16m^{4-2H}}\dd m} \\
& \gs  \frac{1}{2} - \frac{1}{4} \frac{1}{4} - \frac{1}{4} \frac{2}{3}\frac{1}{16} + 2H(2H-1) \frac{1}{4(1-2H)} + 2H(2H-1) \frac{2}{3}\frac{1}{16\rbr{3-2H}} \\
& \gs  \frac{1}{2} - \frac{1}{16}  - \frac{1}{96} -  \frac{H}{2}  - \frac{1}{4} \frac{2}{3}\frac{1}{32}  \gs \frac{11}{64} \gs \frac{1}{8}.
\end{align*}  
Substituting this estimate in \eqref{ala} we finally get
\be \label{weakvarhle12}
n^{1-2H} S^{2H} \gs V_n\rbr{W^H,S} \gs \E S_n^2  \gs \frac{1}{8} n^{1-2H} S^{2H}.
\ee
{\bf Upper bound}.
For $H<1/2$ the increments are negatively correlated and we estimate
\begin{align*}
V_n\rbr{W^H,S} & : = \sup_{0 \ls s_1 < t_1 < s_2 < t_2 < \dots < s_n <t_n \ls S} \E \rbr{ \sum_{i=1}^n \rbr{ W_{t_{i}}^{H}-W_{s_{i}}^{H}}}^2 \\
& = \sup_{0 \ls s_1 < t_1 < s_2 < t_2 < \dots < s_n <t_n \ls S}  \sum_{i=1}^n  \sum_{j=1}^n   \E \rbr{ W_{t_{i}}^{H}-W_{s_{i}}^{H}}\rbr{ W_{t_{j}}^{H}-W_{s_{j}}^{H}} \\
& \ls  \sup_{0 \ls s_1 < t_1 < s_2 < t_2 < \dots < s_n <t_n \ls S}  \sum_{i=1}^n \E \rbr{ W_{t_{i}}^{H}-W_{s_{i}}^{H}}^2 \\
& = \sup_{0 \ls s_1 < t_1 < s_2 < t_2 < \dots < s_n <t_n \ls S} n \frac{1}{n}  \sum_{i=1}^n  \rbr{ t_i - s_i}^{2H} \\
& \ls \sup_{0 \ls s_1 < t_1 < s_2 < t_2 < \dots < s_n <t_n \ls S} n \rbr{\frac{1}{n}  \sum_{i=1}^n  \rbr{ t_i - s_i}}^{2H} \\
& \ls  n \rbr{\frac{1}{n} S}^{2H} = n^{1-2H} S^{2H}.
\end{align*}
\subsubsection*{The case $1/2 \ls H  < 1$}
{\bf Upper and lower bound}. For $H \in [1/2,1)$ the increments are non-negatively correlated and we estimate
\begin{align*}
V_n\rbr{W^H,S} & = \sup_{0 \ls s_1 < t_1 < s_2 < t_2 < \dots < s_n <t_n \ls S}  \E \rbr{ \sum_{i=1}^n \rbr{ W_{t_{i}}^{H}-W_{s_{i}}^{H}}}^2 \\
& \gs  \sup_{0 \ls s_1 < t_1 < s_2 < t_2 < \dots < s_n <t_n \ls S}  \sum_{i=1}^n  \E  \rbr{ W_{t_{i}}^{H}-W_{s_{i}}^{H}}^2  \\
& \gs \sup_{0 \ls s_1 < t_1 < s_2 < t_2 < \dots < s_n <t_n \ls S}  \rbr{{t_1} -{s_1}}^{2H} = S^{2H}.
\end{align*}
On the other hand,
\begin{align*}
V_n\rbr{W^H,S} & = \sup_{0 \ls s_1 < t_1 < s_2 < t_2 < \dots < s_n <t_n \ls S}  \E \rbr{ \sum_{i=1}^n \rbr{ W_{t_{i}}^{H}-W_{s_{i}}^{H}}}^2 \\
& \ls  \sup_{0 \ls s_1 < t_1 < s_2 < t_2 < \dots < s_n <t_n \ls S}  \E \rbr{ \sum_{i=1}^n \rbr{ W_{t_{i}}^{H}-W_{s_{i}}^{H}} + \sum_{i=1}^{n-1} \rbr{ W_{s_{i+1}}^{H}-W_{t_{i}}^{H}} }^2  \\
& = \sup_{0 \ls s_1 < t_1 < s_2 < t_2 < \dots < s_n <t_n \ls S}  \rbr{{t_n} -{s_1}}^{2H} = S^{2H}.
\end{align*}
Thus
\be \label{weakvarhge12}
V_n\rbr{W^H,S}  =  S^{2H}.
\ee

\subsection*{Lemmas A and B}

\subsubsection*{Lemma A}

{\bf Lemma A.} {\em Let ${\cal G}$ be a $\sigma$-field. Assume that a non-negative
real random variable $X$ is ${\cal G}$-measurable and 
\begin{equation}
\P\rbr{X\ge Fv}\le Ae^{-Gv^{\alpha}},\quad v\ge1,\label{eq:exp_tails-1}
\end{equation}
for some positive constants $\alpha$, $F$, $A$ and $G$. Then there
exists a positive constant $\tilde{A}$, depending on $\alpha$, $F$,
$A$ and $G$ only, such that for any $\sigma$-field ${\cal H\subseteq G}$
one has
\begin{equation}
\P\rbr{E\rbr{X^{2}|{\cal H}}\ge\rbr{Fv}^{2}}\le\tilde{A}e^{-(G/2)v^{\alpha}},\quad v\ge0,\label{eq:first_ineq}
\end{equation}
\begin{equation}
\P\rbr{E\rbr{X|{\cal H}}\ge Fv}\le\tilde{A}e^{-(G/2)v^{\alpha}},\quad v\ge0.\label{eq:second_ineq}
\end{equation}
which are equivalent to 
\begin{equation}
\P\rbr{E\rbr{X^{2}|{\cal H}}\ge\rbr{2^{1/\alpha}Fv}^{2}}\le\tilde{A}e^{-Gv^{\alpha}},\quad v\ge0,\label{eq:first_ineq-1}
\end{equation}
\begin{equation}
\P\rbr{E\rbr{X|{\cal H}}\ge2^{1/\alpha}Fv}\le\tilde{A}e^{-Gv^{\alpha}},\quad v\ge0.\label{eq:second_ineq-1}
\end{equation}}
{\bf Proof.} Taking $v=G^{-1/\alpha}x$, the estimate (\ref{eq:exp_tails-1})
may be written as 
\begin{equation}
\P\rbr{\frac{G^{1/\alpha}}{F}X\ge x}\le Ae^{-x^{\alpha}},\quad x\ge G^{1/\alpha}.\label{eq:exp_tails1-1}
\end{equation}
For $\lambda\in(0,1)$ we define the convex function
\[
\varphi_{\lambda}(y)=\exp\rbr{\max\cbr{\lambda|y|^{\alpha/2},\frac{2}{\alpha}-1}}.
\]
We notice that for $x>0$
\[
\dd{\varphi_{\lambda}\rbr{x^{2}}}\le\dd{e^{\lambda x^{\alpha}}}=\lambda\alpha x^{\alpha-1}e^{\lambda x^{\alpha}}\dd x.
\]

Now, using (\ref{eq:exp_tails1-1}) we estimate 
\begin{align*}
\E\varphi_{\lambda}((G^{1/\alpha}/F)^{2}X^{2}) & \le\varphi_{\lambda}\rbr{G^{2/\alpha}}\P\rbr{(G^{1/\alpha}/F)^{2}X^{2}\le G^{2/\alpha}}\\
 & \quad+\int_{G^{1/\alpha}}^{+\ns}\varphi_{\lambda}(x^{2})\P\rbr{(G^{1/\alpha}/F)X\in[x,x+\dd x]}\\
 & \le\varphi_{\lambda}\rbr{G^{2/\alpha}}+\int_{G^{1/\alpha}}^{+\ns}\varphi_{\lambda}(x^{2})\P\rbr{(G^{1/\alpha}/F)X\in[x,x+\dd x]}\\
 & =\varphi_{\lambda}\rbr{G^{2/\alpha}}-\int_{G}^{+\ns}\varphi_{\lambda}(x^{2})\dd{\P\rbr{(G^{1/\alpha}/F)X>x}}\\
 & =\varphi_{\lambda}\rbr{G^{2/\alpha}}-\varphi_{\lambda}\rbr{x^{2}}\P\rbr{(G^{1/\alpha}/F)X>x}|_{x=G^{1/\alpha}}^{x=+\ns}\\
 & +\int_{G}^{+\ns}\P\rbr{(G/F)X>x}\dd{\varphi_{\lambda}(x^{2})}\\
 & \le2\varphi_{\lambda}\rbr{G^{2/\alpha}}+\int_{G^{1/\alpha}}^{+\ns}\P\rbr{(G^{1/\alpha}/F)X>x}\lambda\alpha x^{\alpha-1}e^{\lambda x^{\alpha}}\dd x\\
 & \le2\varphi_{\lambda}\rbr{G^{2/\alpha}}+\int_{G^{1/\alpha}}^{+\ns}Ae^{-x^{\alpha}}\lambda\alpha x^{\alpha-1}e^{\lambda x^{\alpha}}\dd x\\
 & \le2\varphi_{\lambda}\rbr{G^{2/\alpha}}+A\lambda\alpha\int_{G^{1/\alpha}}^{+\ns}x^{\alpha-1}e^{-(1-\lambda)x^{\alpha}}\dd x=:A(\lambda).
\end{align*}
By conditional Jensen's inequality, with probability $1$,
\[
\varphi_{1/2}\rbr{E\rbr{(G^{1/\alpha}/F)^{2}X^{2}|{\cal H}}}\le E\rbr{\varphi_{1/2}\rbr{(G^{1/\alpha}/F)^{2}X^{2}}|{\cal H}}
\]
which yields
\[
\E\varphi_{1/2}\rbr{E\rbr{(G^{1/\alpha}/F)^{2}X^{2}|{\cal H}}}\le\E\varphi_{1/2}\rbr{(G^{1/\alpha}/F)^{2}X^{2}}\le A(1/2).
\]
Now Tchebyschev's inequality yields that for any $x\ge0$, 
\begin{align*}
\P\rbr{E\rbr{(G^{1/\alpha}/F)^{2}X^{2}|{\cal H}}\ge x^{2}} & \le\frac{\E\varphi_{1/2}\rbr{E\rbr{(G^{1/\alpha}/F)^{2}X^{2}|{\cal H}}}}{\varphi_{1/2}\rbr{x^{2}}}\\
 & \le A(1/2)\exp\rbr{-\max\cbr{\frac{1}{2}x^{\alpha},\frac{2}{\alpha}-1}}\\
 & \le A(1/2)\exp\rbr{-\frac{1}{2}x^{\alpha}}.
\end{align*}
Now taking any $v\ge0$ and substituting $x=G^{1/\alpha}v$, and we
get (\ref{eq:first_ineq}) with $\tilde{A}=A(1/2)$:
\[
\P\rbr{E\rbr{X^{2}|{\cal H}}\ge F^{2}v^{2}}\le A(1/2)e^{-(1/2)Gv^{\alpha}}.
\]
The second inequality (\ref{eq:second_ineq}) follows immediately
from (\ref{eq:first_ineq}) and the estimate $E\rbr{X|{\cal H}}\le\sqrt{E\rbr{X^{2}|{\cal H}}}$.
\hfill $\blacksquare$

\subsubsection*{Lemma B}

Let $L$ be a positive integer no smaller than $2$ and the sequence
of random variables $M=\rbr{M_{1},M_{2},\ldots,M_{L}}$, such that
$M_{1}\equiv0$, be a martingale with respect to the filtration $\mathbb{G}=\rbr{{\cal G}_{1},{\cal G}_{2},\ldots,{\cal G}_{L}}$.
We introduce the martingale differences 
\[
\Delta M_{l}=M_{l}-M_{l-1},\quad l=2,3,\ldots,L,
\]
the martingale square brackets
\begin{align*}
[M]_{n} & :=\sum_{l=2}^{n}\rbr{\Delta M_{l}}^{2},\quad n=1,2,\ldots,L,
\end{align*}
and the previsible martingale square brackets
\[
\sb M_{n}:=\sum_{l=2}^{n}E\sbr{\rbr{\Delta M_{l}}^{2}|{\cal G}_{l-1}},\quad n=1,2,\ldots,L.
\]
{\bf Lemma B.} {\em Assume that the sequence of real random variables $\rbr{M_{1},M_{2},\ldots,M_{L}}$,
such that $M_{1}\equiv0$, is a martingale with respect to the filtration
$\mathbb{G}=\rbr{{\cal G}_{1},{\cal G}_{2},\ldots,{\cal G}_{L}}$
and the random variables $X_{2},\ldots,X_{L}$ are such that $X_{l}$
is ${\cal G}_{l}$-measurable and
\[
\Delta M_{l}=M_{l}-M_{l-1}=X_{l}-E\rbr{X_{l}|{\cal G}_{l-1}},\quad l=2,3,\ldots,L.
\]
Moreover, assume that there exist some positive constants $\alpha$,
$F$, $A$ and $G$ such that for each $l=2,3,\ldots,L$,
\begin{equation}
\P\rbr{X_{l}\ge Fv}\le Ae^{-Gv^{\alpha}},\quad v\ge1.\label{eq:exp_tails-1-1}
\end{equation}
If $\tilde{A}$ is the same constant as in Lemma A, then for any $x\ge0$
and $v>0$ 
\begin{equation}
\P\rbr{\left|M_{L}\right|\ge x}\le2\exp\rbr{-\frac{x^{2}}{10L\rbr{Fv}^{2}}}+6L\tilde{A}_{H}\exp\rbr{-(G/2)v^{\alpha}}.\label{eq:lemma_B}
\end{equation}}
{\bf Proof.} If ${\cal G}$ and ${\cal H}\subseteq{\cal G}$ are $\sigma$-fields
of subsets of $\Omega$, then for any ${\cal G}$-measurable random
variable $X:\Omega\ra\R$ with finite second moment we have 
\[
E\rbr{\rbr{X-E(X|{\cal H})}^{2}|{\cal H}}=E\rbr{X^{2}|{\cal {\cal H}}}-\rbr{E(X|{\cal {\cal H}})}^{2}
\]
and the following estimates hold almost surely 
\begin{equation}
\rbr{X-E(X|{\cal {\cal H}})}^{2}\le2X^{2}+2\rbr{E(X|{\cal {\cal H}})}^{2}\le2X^{2}+2E(X^{2}|{\cal H})\label{eq:first_cond-1}
\end{equation}
and
\begin{equation}
E\rbr{\rbr{X-E(X|{\cal H})}^{2}|{\cal {\cal H}}}\le E\rbr{X^{2}|{\cal H}}.\label{eq:second_cond-1}
\end{equation}
By (\ref{eq:first_cond-1}), for $l=2,3,\ldots,L$,
\begin{align}
\rbr{\Delta M_{l}}^{2}=\cbr{X_{l}-E\sbr{X_{l}|{\cal G}_{l-1}}}^{2} & \le2X_{l}^{2}+2E\sbr{X_{l}^{2}|{\cal G}_{l-1}}.\label{eq:decomp-1}
\end{align}
Lemma A gives that for any $v\ge0$
\begin{align*}
 & \P\rbr{\rbr{\Delta M_{l}}^{2}\ge4\rbr{Fv}^{2}}\le\P\rbr{2X_{l}^{2}+2E\sbr{X_{l}^{2}|{\cal G}_{l-1}}\ge4\rbr{Fv}^{2}}\\
 & \le\P\rbr{2X_{l}^{2}\ge2\rbr{Fv}^{2}}+\P\rbr{2E\sbr{X_{l}^{2}|{\cal G}_{l-1}}\ge2\rbr{Fv}^{2}}\\
 & \le2\tilde{A}\exp\rbr{-(G/2)v^{\alpha}}.
\end{align*}
Similarly, by (\ref{eq:second_cond-1}) and Lemma A, for any $v\ge0$
\begin{align*}
 & \P\rbr{E\sbr{\rbr{\Delta M_{l}}^{2}|{\cal F}_{l-1}}\ge\rbr{Fv}^{2}}\le\P\rbr{E\sbr{X_{l}^{2}|{\cal G}_{l-1}}\ge\rbr{Fv}^{2}}\\
 & \le\tilde{A}\exp\rbr{-(G/2)v^{\alpha}}.
\end{align*}
Two last estimates give for $v>0$ 
\begin{align}
 & \P\rbr{[M]_{L}+\sb M_{L}\ge5L(Fv)^{2}}=\P\rbr{\sum_{l=2}^{L}\rbr{\Delta M_{l}}^{2}+\sum_{l=2}^{L}E\sbr{\rbr{\Delta M_{l}}^{2}|{\cal F}_{l-1}}\ge5L(Fv)^{2}}\nonumber \\
 & \le\sum_{l=2}^{L}\P\rbr{\rbr{\Delta M_{l}}^{2}\ge4(Fv)^{2}}+\sum_{l=2}^{L}\P\rbr{E\sbr{\rbr{\Delta M_{l}}^{2}|{\cal F}_{l-1}}\ge(Fv)^{2}}\nonumber \\
 & \le3L\tilde{A}\exp\rbr{-(G/2)v^{\alpha}}.\label{eq:rough_square_estim-1}
\end{align}
Next, we will use \cite[Theorem 3.26]{Rio}, which states that
for any $x\ge0$, $y>0$
\begin{equation}
\P\rbr{M_{L}\ge x,\ [M]_{L}+\sb M_{L}\le y}\le\exp\rbr{-\frac{x^{2}}{2y}}.\label{eq:Rio-1}
\end{equation}
This yields for $x\ge0$ and $v>0$
\begin{align}
\P\rbr{M_{L}\ge x}= & \P\rbr{M_{L}\ge x,\ [M]_{L}+\sb M_{L}\le5L(Fv)^{2}}+\P\rbr{[M]_{L}+\sb M_{L}\ge5L(Fv)^{2}}\nonumber \\
\le & \exp\rbr{-\frac{x^{2}}{10L\rbr{Fv}^{2}}}+3L\tilde{A}_{H}\exp\rbr{-(G/2)v^{\alpha}}.\label{eq:last-1}
\end{align}
A symmetric estimate holds for $-M_{L}$, which together with (\ref{eq:last-1})
yields (\ref{eq:lemma_B}).

\hfill $\blacksquare$

\bibliographystyle{alpha} % Style BST file
\bibliography{/Users/rafallochowski/biblio/biblio}

\end{document}